\documentclass[12pt]{article} 
\usepackage{latexsym, 
amscd, 
graphicx, epsfig, amssymb}

\newtheorem{thm}{Theorem}[section]
\newtheorem{defn}[thm]{Definition}
\newtheorem{prop}[thm]{Proposition}
\newtheorem{cor}[thm]{Corollary}
\newtheorem{lemma}[thm]{Lemma}

\newtheorem{cond}[thm]{Condition}
\newcommand{\halmos}{\rule{1ex}{1.4ex}}

\newcommand{\binom}[2]{{{#1}\choose {#2}}}
\newcommand{\text}[1]{\mbox{\rm #1}}
\newcommand{\mod}{\;\;\mbox{\rm mod}\;}

\newcommand{\nn}{\nonumber \\}

 \newcommand{\res}{\mbox{\rm Res}}

 \newcommand{\pf}{{\it Proof.}\hspace{2ex}}
 
 \newcommand{\epfv}{\hspace*{\fill}\mbox{$\halmos$}\vspace{1em}}

\newcommand{\tr}{\mbox{\rm Tr}}
\newcommand{\A}{\mathcal{A}}
\newcommand{\Y}{\mathcal{Y}}
\newcommand{\C}{\mathbb{C}}
\newcommand{\Z}{\mathbb{Z}}
\newcommand{\R}{\mathbb{R}}
\newcommand{\Q}{\mathbb{Q}}

\newcommand{\F}{\mathcal{F}}

\newcommand{\V}{\mathcal{V}}
\newcommand{\one}{\mathbf{1}}

\title{ {\bf Modular invariance for conformal full field algebras} }
\date{}
\author{Yi-Zhi Huang and Liang Kong}
\begin{document}

\bibliographystyle{alpha}
\maketitle

\begin{abstract} 
Let $V^{L}$ and $V^{R}$ be simple vertex operator algebras satisfying 
certain natural uniqueness-of-vacuum, complete reducibility
and cofiniteness conditions and let $F$ be a conformal full field 
algebra over $V^{L}\otimes V^{R}$. 
We prove that the 
$q_{\tau}$-$\overline{q_{\tau}}$-traces (natural traces involving 
$q_{\tau}=e^{2\pi i\tau}$ and 
$\overline{q_{\tau}}=\overline{e^{2\pi i\tau}}$)
of geometrically 
modified genus-zero correlation functions for $F$
are convergent in suitable regions
and can be extended to doubly periodic functions with periods 
$1$ and $\tau$. 
We obtain necessary and 
sufficient conditions for these functions to be modular invariant. 
In the case that $V^{L}=V^{R}$ and $F$ is  one of those constructed 
by the authors in \cite{HK}, we prove that all these functions 
are modular invariant. 
\end{abstract}

\renewcommand{\theequation}{\thesection.\arabic{equation}}
\renewcommand{\thethm}{\thesection.\arabic{thm}}
\setcounter{equation}{0}
\setcounter{thm}{0}
\setcounter{section}{-1}

\section{Introduction}

In this paper, we construct  genus-one full conformal field theories
(genus-one conformal field theories 
with both chiral and anti-chiral parts)
from genus-zero full conformal field theories. 
More precisely, we construct genus-one correlation functions 
from genus-zero correlation functions for a conformal full field algebra
over $V^{L}\otimes V^{R}$ in the sense of \cite{HK}, where 
$V^{L}$ and $V^{R}$ are vertex operator algebras satisfying 
certain natural uniqueness-of-vacuum, complete reducibility
and cofiniteness conditions (see below).

Since Kontsevich and Segal (see \cite{S1}, \cite{S2}
and \cite{S3}) gave a geometric definition of 
(two-dimensional) conformal field theory in 1987 by axiomatizing 
the properties of path integrals used by physicists, 
constructing conformal field theories satisfying 
this definition became an important unsolved mathematical problem. 
Around the same time, E. Verlinde \cite{V} and Moore-Seiberg 
\cite{MS1} \cite{MS2} (see also \cite{MS3}) made a
major breakthrough by discovering the relation between 
fusion rules and modular transformations and 
the modular tensor category structures associated to
rational conformal field theories, assuming 
that these theories had been constructed. 

For some of the applications of conformal 
field theories (for example, for the construction and study of 
knot and three-manifold invariants and the proof of the Verlinde 
formula in algebraic geometry), a construction of modular functors
introduced by Segal \cite{S2} \cite{S3} is enough. See, for example,
\cite{T} and \cite{BK} for the construction of some 
examples of modular functors. 
In \cite{H11} (see also \cite{H8} for an announcement
and \cite{H9} for an exposition), 
the first author constructed modular tensor categories
from representations of vertex operator algebras satisfying 
the conditions mentioned above. Combining with the results in 
\cite{T} and \cite{BK}, this result of \cite{H11} 
in fact gives a construction 
modular functors from representations of vertex operator algebras
satisfying the conditions alluded above. 

However, for some other applications, modular functors are 
far from enough.
One extreme example is the conformal-field-theoretic construction 
used in Frenkel-Lepowsky-Meurman's 
proof of the McKay-Thompson conjecture
on the existence of the moonshine module 
\cite{FLM} and in Borcherds' proof of the Monstrous moonshine conjecture
\cite{B}. In this example, the modular functor associated to the 
moonshine module vertex operator algebra is actually trivial
and thus does not play any role. On the other hand, we know that
the construction of the moonshine module vertex operator algebra
by Frenkel-Lepowsky-Meurman
and the proof of the Monstrous moonshine conjecture
by Borcherds can be interpreted as one of the deepest applications of 
conformal field theory. 

A program to construct conformal field theories was first initiated 
by I.~Frenkel even before Kontsevich and Segal gave their geometric
definition of conformal field theory. Under the direction of 
I. Frenkel, He \cite{Ts} constructed genus-zero and genus-one parts of 
conformal field theories associated to tori. 
Tsukada constructed not only chiral theories but also  full theories
which contain both chiral and antichiral parts. 
In \cite{H1}, \cite{H2}, \cite{H3}, \cite{H4},  \cite{H5},  \cite{H6} and 
\cite{H7}, the first author constructed 
chiral genus-zero and genus-one conformal field theories 
from representations of vertex operator algebras
satisfying suitable conditions. In \cite{HK}, the authors constructed 
genus-zero full conformal field theories from representations 
of vertex operator algebras satisfying the conditions mentioned above
and given precisely below and in Section 3
(see also \cite{K} and \cite{K2}). 

The present paper is a continuation of
the paper \cite{HK} and the results obtained in 
the present paper can be viewed as one step in a program 
of constructing conformal field theories from representations 
of vertex operator algebras. As is mentioned above, 
we construct genus-one full conformal field theories in this paper. 
More precisely, 
let $V^{L}$ and $V^{R}$ be simple vertex operator algebras satisfying 
the following conditions for a vertex operator algebra $V$:
(i) For $n<0$, $V_{(n)}=0$, $V_{(0)}=\mathbb{C}\mathbf{1}$ and 
$W_{(0)}=0$
for irreducible $V$-module $W$ not isomorphic to $V$.
(ii) Every $\mathbb{N}$-gradable weak $V$-module is completely reducible.
(iii) $V$ is $C_{2}$-cofinite, that is, $\dim V/C_{2}(V)<\infty$, 
where $C_{2}(V)$ is the subspace of $V$ spanned by elements of 
the form $u_{-2}v$ for $u, v\in V$.
Let $F$ be a conformal full field algebra over $V^{L}\otimes V^{R}$ (see 
\cite{HK} or Section 1 for the definition and basic properties).
We construct genus-one correlation functions using 
$q_{\tau}$-$\overline{q_{\tau}}$-traces 
(natural traces involving 
$q_{\tau}=e^{2\pi i\tau}$ and $\overline{q_{\tau}}=e^{2\pi i\tau}$) 
of geometrically 
modified genus-zero correlation functions for $F$.
We prove that these functions are  doubly periodic and
we obtain conditions which are equivalent to 
the modular invariance of these functions. 
When $V^{L}=V^{R}=V$ and $F$ is the conformal full field algebras 
over $V\otimes V$ constructed 
by the authors in \cite{HK}, we prove that these conditions are 
satisfied and thus all these functions 
are modular invariant. 

We note that based on the
existence of the structure of a modular tensor category on the category
of modules for a vertex operator algebra, the existence of conformal
blocks with monodromies compatible with the modular tensor category and
all the necessary convergence properties,
Felder, Fr\"{o}hlich,
Fuchs and Schweigert \cite{FFFS},  Fuchs, Runkel, Schweigert and
Fjelstad \cite{FRS1} \cite{FRS2} \cite{FRS3} \cite{FRS4}
\cite{FRS5} \cite{FFRS} \cite{FjFRS2}, and 
Fr\"{o}hlich, Fuchs, Runkel and Schweigert \cite{FrFRS} studied
open-closed conformal field theories (in particular full (closed)
conformal field theories) using the theory of tensor categories and
three-dimensional topological field theories. In particular, they
constructed correlation functions as states in some three-dimensional
topological field theories. 
In \cite{HK} and the present paper, what we need in our work are theorems 
proved by the first author in \cite{H6}, \cite{H7}, \cite{H10}
and \cite{H11}
when the vertex operator algebras we start with
satisfy some natural conditions.
Our work in
\cite{HK} and the present paper not only replaced these fundamental 
but hard-to-verify
assumptions by natural, purely algebraic and easy-to-verify conditions
on vertex operator algebras, but also constructed explicitly genus-zero
and genus-one correlation functions from intertwining operators for the
vertex operator algebras.

The present paper depends heavily on the results obtained in \cite{H7}, 
even more heavily than the papers \cite{H10} and \cite{H11}. 
In \cite{H10} and \cite{H11}, we need only
the properties of certain special two-point
genus-one chiral correlation functions obtained in \cite{H7},
namely, those obtained from iterates of intertwining operators with the
intermediate modules being the vertex operator algebra itself.
In this paper, we need the full strength of the results obtained 
in \cite{H7}. In particular, using the results of \cite{H7}, 
we prove a symmetry property of the matrix elements of the 
actions of the modular transformation $\tau\mapsto -\frac{1}{\tau}$
on the space of one-point genus-one correlation functions
(Theorem \ref{s-iden-thm}), which is 
a generalization of the symmetry property of the matrix elements 
of the action of the same modular transformation on the space
of vacuum characters (Theorem 5.6 in \cite{H10}). This generalization
is exactly what we need in Section 5 to prove the modular invariance of the 
genus-one correlation functions for the conformal full field algebras 
constructed in \cite{HK}. To prove Theorem \ref{s-iden-thm},
we prove an explicit formula (\ref{s-formula}) for these matrix elements. 
In the case that $a=e$, we recover the 
formula for the matrix elements of the action of the 
same modular transformation on the space of vacuum characters, 
obtained first by Moore-Seiberg from the Verlinde formula and 
proved for vertex operator algebras satisfying the conditions above 
by the first author in \cite{H7}. Note that  in the 
case $a=e$ this formula 
(\ref{s-formula}) was shown 
in \cite{H11} to 
be equivalent to the nondegeneracy property of the 
semisimple ribbon (tensor) category of modules for the vertex operator algebra. 
In the special case of discrete series,
a formula for these matrix elements was given in \cite{MS3}. 

The present paper is organized as follows: In Section 1,
we recall the basic definition and constructions in the 
theory of conformal full field algebras over $V^{L}\otimes V^{R}$
given first in \cite{HK}. 
In Section 2, we recall the chiral genus-one theory 
constructed from intertwining 
operator algebras in \cite{H7}. These two sections 
are given here for the convenience of the reader. 
The reader is referred to \cite{HK} and \cite{H7} for more details. 
In Section 3, we prove the convergence of 
$q_{\tau}$-$\overline{q_{\tau}}$-traces of genus-zero correlation 
functions in suitable regions and show that these can be 
extended to doubly periodic functions with periods 
$1$ and $\tau$. We also give conditions which are 
equivalent to the modular invariance of these functions in 
this section. In Section 4, we study the matrix elements of 
the action of the modular transformation $S: \tau \mapsto -\frac{1}{\tau}$
on chiral genus-one correlation functions. 
In particular we obtain an explicit formula  for
these matrix elements. This formula allows us to derive a 
a symmetry property of these matrix elements.
In Section 5, for the conformal full 
field algebras constructed in \cite{HK},
we prove the modular invariance  of the correlation
functions using the results obtained in Sections 3 and 4.

\paragraph{Acknowledgment} 
The first author is partially supported 
by NSF grant DMS-0401302.

\renewcommand{\theequation}{\thesection.\arabic{equation}}
\renewcommand{\thethm}{\thesection.\arabic{thm}}
\setcounter{equation}{0}
\setcounter{thm}{0}

\section{Conformal full field algebras}

In this section, we recall the notion of conformal full field algebra
over $V^{L}\otimes V^{LR}$ introduced in \cite{HK} and review the 
construction 
of conformal full field algebras over $V^{L}\otimes V^{LR}$
when $V^{L}=V^{R}$ satisfying suitable conditions
given in the same paper. See \cite{HK} for more details and 
other variants of full field algebras. 

Let $\mathbb{F}_{n}(\C) =\{ (z_{1}, \ldots, z_{n})\in \C^n \;|\; z_i\neq z_j
\mbox{ if } i\neq j \}$. 
For an $\R\times \R$-graded vector space $F=\coprod_{r, s\in \R}F_{(r, s)}$,
let $F'=\coprod_{r, s\in \R}F_{(r, s)}^{*}$ and 
$\overline{F}=\prod_{r, s\in \R}F_{(r, s)}$ be the graded dual 
and algebraic completion
of $F$, respectively. For $r, s\in \R$, let 
$P_{r, s}$ be the projection from $F$ or $\overline{F}$ to $F_{(r, s)}$.
A series $\sum f_{n}$ in $\overline{F}$ is said to be {\it
absolutely convergent} if for any $f'\in F'$, $\sum|\langle f', f_{n}\rangle|$
is convergent. The sums $\sum|\langle f', f_{n}\rangle|$ for $f'\in 
F'$ define a linear functional on $F'$. We call this linear functional 
the {\it sum} of the series and denote it by the same notation 
$\sum f_{n}$. If the homogeneous subspaces of $F$ are all finite-dimensional, 
then $\overline{F}=(F')^{*}$ and, in this case, the sum 
of an absolutely convergent series is always in $\overline{F}$. 
When the sum is in $\overline{F}$,
we say that the series is {\it absolutely convergent in $\overline{F}$}.

\begin{defn}\label{ffa}
{\rm Let $(V^{L}, Y^{L}, \one^{L}, \omega^{L})$ 
and $(V^{R}, Y^{R}, \one^{R}, \omega^{R})$ be vertex operator algebras. 
A {\it conformal full field algebra over $V^{L}\otimes V^{R}$} 
is an $\R\times \R$-graded vector space $F=\coprod_{r, s\in \R}F_{(r, s)}$
(graded by {\it left conformal weight}
or simply {\it left weight} and by {\it right conformal weight}
or simply {\it right weight}),  
equipped with {\it correlation function maps}
$$\begin{array}{rrcl}
m_{n}: & F^{\otimes n}\times
\mathbb{F}_{n}(\C)&\to& \overline{F}\\
&(u_{1}\otimes \cdots \otimes u_{n}, (z_{1}, \dots, z_{n}))&\mapsto&
m_{n}(u_{1}, \dots, u_{n}; z_{1}, \bar{z}_{1}, \dots, z_{n}, \bar{z}_n),
\end{array}$$
for $n\in \Z_+$, an injective grading-preserving linear map $\rho$ from 
the  $V^{L}\otimes V^{R}$ to $F$
satisfying the following axioms: 

\begin{enumerate}

\item There exists $M\in \R$ such that 
$F_{(r,s)}=0$ if $n<M$ or $m<M$.

\item $\dim F_{(r,s)}<\infty$ for $m, n\in \R$. 

\item For $n\in \Z_+$, 
$m_n(u_{1}, \dots, u_{n}; z_{1}, \bar{z}_{1}, \dots, z_{n}, \bar{z}_n)$ is 
linear in $u_{1}, \dots, u_{n}$ and smooth in the real and imaginary 
parts of $z_{1}, \dots, z_{n}$.

\item For $u\in F$, $m_{1}(u; 0, 0)=u$. 

\item For $n\in \Z_{+}$, $u_{1}, \dots, u_{n}\in F$, 
\begin{eqnarray*}
&m_{n+1}(u_{1}, \dots, u_{n}, \one; z_{1}, \bar{z}_{1}, \dots, z_{n}, \bar{z}_{n}, 
z_{n+1}, \bar{z}_{n+1})&\nn
&=m_n(u_{1}, \dots, u_{n}; z_{1}, \bar{z}_{1}, \dots, z_{n}, \bar{z}_{n}),&
\end{eqnarray*}
where $\one =\rho(\one^{L}\otimes \one^{R})$.

\item The {\it  convergence property}:   
For $k, l_{1}, \ldots, l_{k} \in \Z_+$ and 
$u_{1}^{(1)},\ldots, u_{l_{1}}^{(1)},\ldots, u_{1}^{(k)}$, $\dots,
u_{l_{k}}^{(k)} \in F$,  the series 
\begin{eqnarray}\label{ffa-conv-axiom}
\lefteqn{\sum_{r_{1}, s_{1}, \ldots, r_k, s_{k}} 
m_{k}(P_{r_{1}, s_{1}} m_{l_{1}}(u_{1}^{(1)}, 
\ldots, u_{l_{1}}^{(1)}; z_{1}^{(1)}, \bar{z}_{1}^{(1)}, \dots, 
z_{l_{1}}^{(1)}, \bar{z}_{l_{1}}^{(1)} ), \dots,}  \nn
&& 
P_{r_k, s_{k}} m_{l_{k}} (u_{1}^{(k)}, \ldots, u_{l_{k}}^{(k)};
z_{1}^{(k)}, \bar{z}_{1}^{(k)}, \ldots, z_{l_{k}}^{(k)}, 
\bar{z}_{l_{k}}^{(k)}); z^{(0)}_{1}, \bar{z}_{1}^{(0)}, \ldots, 
z^{(0)}_{k}, \bar{z}^{(0)}_k)\nn
&&
\end{eqnarray}
converges absolutely to  
\begin{eqnarray}\label{sum}
\lefteqn{m_{l_{1}+ \cdots + l_k}(u_{1}^{(1)}, \dots, u_{l_k}^{(k)};
z_{1}^{(1)}+z^{(0)}_{1}, \bar{z}_{1}^{(1)}+ \bar{z}^{(0)}_{1}, \dots,
z_{l_{1}}^{(1)}+z^{(0)}_{1}, }  \nn
&&\quad
\bar{z}_{l_{1}}^{(1)}+\bar{z}^{(0)}_{1},  
\dots, 
z_{1}^{(k)}+z^{(0)}_k,  
\bar{z}_{1}^{(k)}+ \bar{z}^{(0)}_k, \dots,
z_{l_k}^{(k)}+z^{(0)}_k, \bar{z}_{l_k}^{(k)}+\bar{z}^{(0)}_k)  \nn
&&
\end{eqnarray}
when $|z_p^{(i)}| + |z_q^{(j)}|< |z^{(0)}_i
-z^{(0)}_j|$ for $i,j=1, \ldots, k$, $i\ne j$ and for $p=1, 
\dots,  l_i$ and $q=1, \dots, l_j$.

\item The {\it permutation property}: For any $n\in \Z_{+}$ and 
any $\sigma\in S_{n}$,
we have 
\begin{eqnarray}
\lefteqn{m_{n}(u_{1}, \dots, u_{n}; z_{1}, \bar{z}_{1}, 
\dots, z_{n}, \bar{z}_{n})} \nn
&& = m_{n}(u_{\sigma(1)}, \dots, u_{\sigma(n)}; 
z_{\sigma(1)}, \bar{z}_{\sigma(1)}, \dots, 
z_{\sigma(n)}, \bar{z}_{\sigma(n)})   \label{ffa-perm-axiom}
\end{eqnarray}
for $u_{1}, \dots, u_{n}\in F$ and $(z_{1}, \dots, z_{n})\in 
\mathbb{F}_{n}(\C)$.

\item Let 
$$\begin{array}{rrcl}
\mathbb{Y}:& F^{\otimes 2} \times \C^{\times} &\to&
\overline{F}\\
&(u\otimes v, z, \bar{z})&\mapsto& \mathbb{Y}(u; z, \bar{z})v
\end{array}$$ 
be given by
$$\mathbb{Y}(u; z, \bar{z})v = m_{2}(u \otimes v; z,\bar{z}, 0, 0)$$
for $u, v\in F$. Then 
$$\mathbb{Y}(\rho(u^{L}\otimes u^{R}); z, \bar{z})\rho(v^{L}\otimes v^{R})
=\rho(Y^{L}(u^{L}, z)u^{R}\otimes Y^{R}(u^{R}, \overline{z})v^{R})$$
for $u^{L}, v^{L}\in V^{L}$, $u^{R}, v^{R}\in V^{R}$, and
there exist operators 
$L^{L}(n)$ and $L^{R}(n)$ for $n\in \Z$ such that
\begin{eqnarray*}
\mathbb{Y}(\rho(\omega^{L}\otimes \mathbf{1}^{R}); z, \bar{z})
&=&\sum_{n\in \Z}L^{L}(n)z^{-n-2},\\
\mathbb{Y}(\rho(\mathbf{1}^{L}\otimes \omega^{R}); z, \bar{z})
&=&\sum_{n\in \Z}L^{R}(n)\bar{z}^{-n-2}.
\end{eqnarray*}

%\item For $n\in \Z_{+}$, $a\in \R$, $u_{1}, 
%\dots, u_{n}\in F$,
%\begin{eqnarray*}
%\lefteqn{e^{a(L^{L}(0)+L^{R}(0))}m_{n}(u_{1}, 
%\dots, u_{n}; z_{1}, \bar{z}_{1},
%\dots, z_{n}, \bar{z}_{n})}\nn
%&&=m_{n}(e^{a(L^{L}(0)+L^{R}(0))}u_{1}, \dots, e^{a(L^{L}(0)+L^{R}(0))}u_{n}; 
%e^{a}z_{1}, e^{a}\bar{z}_{1},
%\dots, e^{a}z_{n}, e^{a}\bar{z}_{n}).
%\end{eqnarray*}

\item The {\it single-valuedness property}: 
$e^{2\pi i (L^L(0)-L^R(0))} = I_F$. 

\item {\it The $L^L(0)$- and $L^R(0)$-grading properties}:
For $r, s\in \R$ and $u\in F_{(r, s)}$,
\begin{eqnarray*}
L^{L}(0)u&=&ru,\\
L^{R}(0)u&=&su,
\end{eqnarray*}

\item {\it The $L^L(0)$- and $L^R(0)$-bracket properties}: For $u\in F$, 
\begin{eqnarray}
\left[L^L(0), \mathbb{Y}(u; z,\bar{z})\right] 
&=& z\frac{\partial}{\partial z} \mathbb{Y}(u; z, \bar{z})
+\mathbb{Y}(L^L(0) u; z,\bar{z})     \label{d-l} \\
\left[ L^R(0), \mathbb{Y}(u; z, \bar{z})\right]
&=& \bar{z}\frac{\partial}{\partial \bar{z}} \mathbb{Y}(u; z,\bar{z})
+ \mathbb{Y}(L^R(0) u; z, \bar{z}). \label{d-r}
\end{eqnarray}

\item The {\it $L^{L}(-1)$- and $L^{R}(-1)$-derivative property}: For $u\in F$,
\begin{eqnarray}
&&\left[ L^L(-1), \mathbb{Y}(u; z,\bar{z})\right] = 
\mathbb{Y}(L^{L}(-1)u; z,\bar{z}) =
\frac{\partial}{\partial z}\mathbb{Y}(u; z, \bar{z}),\\
&&\left[ L^R(-1), \mathbb{Y}(u; z, \bar{z})\right] = 
\mathbb{Y}(L^{R}(-1)u; z, \bar{z}) =
\frac{\partial}{\partial \bar{z}}\mathbb{Y}(u; z, \bar{z}). \label{D-L-R-b-d}
\end{eqnarray}

\end{enumerate}
}

\end{defn}

We denote the conformal full field algebra over $V^{L}\otimes V^{R}$
defined above by 
$(F, m, \rho)$ or simply by $F$.
In the definition above, we use the notations 
$$m_{n}(u_{1}, \dots, u_{n}; z_{1},
\bar{z}_{1} \dots, z_{n}, \bar{z}_n)$$ instead of
$$m_{n}(u_{1}, \dots,
u_{n}; z_{1},\dots, z_{n})$$
to emphasis that 
these are in general not holomorphic in $z_{1},\dots, z_{n}$.
For $u'\in F'$, $u_{1}, \dots, u_{n}\in F$, 
$$\langle u', m_{n}(u_{1}, \dots, u_{n}; z_{1},
\bar{z}_{1} \dots, z_{n}, \bar{z}_n)\rangle$$
as a function of $z_{1}, \dots, z_{n}$ 
is called a {\it correlation function}.
The map $\mathbb{Y}$ is called the {\it full 
vertex operator map}
and for $u\in F$, $\mathbb{Y}(u; z, \bar{z})$ is called the {\it 
full vertex operator} associated to $u$. The element $\rho(\one \otimes 
\one)$ is called the {\it vacuum} of $F$. The elements
$\rho(\omega^{L}\otimes \mathbf{1}^{R})$ and 
$\rho(\mathbf{1}^{L}\otimes\omega^{R})$ are called the 
{\it left conformal element} $\rho(\omega^{L}\otimes \mathbf{1}^{R})$
and {\it right conformal element}, respectively.

{\it Homomorphisms} and {\it isomorphisms} for conformal 
full field algebras over $V^{L}\otimes V^{R}$ are
defined in the obvious way. 

For a conformal full field algebra over $V^{L}\otimes V^{R}$, 
a formal full vertex operator map 
\begin{eqnarray*}
\mathbb{Y}_{f}: F\otimes F&\to& F\{x, \bar{x}\}\nn
u\otimes v&\mapsto& \Y_{f}(u; x, \overline{x})v
\end{eqnarray*}
was obtained in \cite{HK} such that 
$$\mathbb{Y}(u; z, \overline{z})=\mathbb{Y}_{f}(u; x, 
\overline{x})|_{x^{r}=e^{r\log z}, 
\overline{x}^{s}=e^{s\overline{\log z}}, r, s\in \R}$$
for $u\in F$ and $z\in \C^{\times}$. We can also substitute $e^{r\log z}$
and $e^{s\overline{\log \zeta}}$ for $x^{r}$ and $\overline{x}^{s}$
to obtain 
$\mathbb{Y}(u; z, \zeta)$ for $u\in F$.

For the operators
$L^{L}(n)$ and $L^{R}(n)$ for $n\in \Z$, we have the 
following bracket formulas:
For $m, n\in \Z$, 
\begin{eqnarray*}
[L^{L}(m), L^{L}(n)]
&=&(m-n)L^{L}(m+n)+\frac{c^{L}}{12}(m^{3}-m)\delta_{m+n, 0},\\
{[L^{R}(m), L^{R}(n)]}
&=&(m-n)L^{R}(m+n)+\frac{c^{R}}{12}(m^{3}-m)\delta_{m+n, 0},\\
{[L^{L}(m), L^{R}(n)]}&=&0.
\end{eqnarray*}

Let $F$ be a module for the vertex operator algebra
$V^L\otimes V^R$ and  $\Y$ an intertwining operator 
of type $\binom{F}{FF}$. In \cite{HK}, a splitting 
$\mathbb{Y}^{\Y}: (F\otimes F)\times \C^{\times}\to \overline{F}$
and a formal splitting 
$\mathbb{Y}_{f}^{\Y}: F\otimes F\to \overline{F}\{x, \overline{x}\}$
of $\Y$ are constructed. Substituting $e^{r\log z}$ 
and $e^{s\overline{\log \zeta}}$
for $x^{r}$ and $\overline{x}^{s}$ in the images of $F\otimes F$
under $\mathbb{Y}_{f}^{\Y}$, we obtain an analytic splitting 
$\mathbb{Y}_{\rm an}^{\Y}: (F\otimes F)\times (\C^{\times}\times \C^{\times})
\to \overline{F}$.

One of the main result of \cite{HK} is the following theorem:

\begin{thm}  \label{r-conf-alg-2-prop}
Let $V^{L}$ and $V^{R}$ be vertex operator algebras satisfying the 
following conditions for a vertex operator algebra $V$:
(i) Every $\C$-graded $L(0)$-semisimple generalized $V$-module 
is a direct sum of $\C$-graded irreducible $V$-modules.
(ii) There are only finitely many inequivalent $\C$-graded irreducible $V$-modules
and they are all $\R$-graded. (iii) Every $\R$-graded irreducible $V$-module
$W$ satisfies the $C_{1}$-cofiniteness condition, that is, 
$\dim W/C_{1}(W)<\infty$, where $C_{1}(V)$ is the subspace of $V$ spanned 
by elements of the form $u_{1}w$ for $u\in V_{+}=\coprod_{n>0}V_{(n)}$ and $w\in W$.
Then a conformal full field algebra over 
$V^L\otimes V^R$ is equivalent to a module $F$ 
for the vertex operator algebra $V^L\otimes V^R$ equipped with
an intertwining operator
$\Y$ of type $\binom{F}{FF}$ and an injective 
module map $\rho: V^L\otimes V^R \to F$, 
satisfying the following conditions:
\begin{enumerate}

\item The {\it identity property}: 
$\Y(\rho(\one^L\otimes \one^R), x)=I_{F}$.

\item The {\it creation property}: For $u\in F$, 
$\lim_{x\rightarrow 0} \Y(u, x)\rho(\one^L\otimes \one^R)=u$.

\item The {\it associativity}: For
$u,v,w\in F$ and $w'\in F'$, 
\begin{eqnarray}  \label{asso-z-zeta}
\lefteqn{\langle w', \mathbb{Y}^{\Y}_{\rm an}(u; z_{1},\zeta_{1})
\mathbb{Y}^{\Y}_{\rm an}(v; z_{2}, \zeta_{2})w\rangle} \nn
&&= \langle w', 
\mathbb{Y}^{\Y}_{\rm an}(\mathbb{Y}^{\Y}_{\rm an}(u; z_{1}-z_{2}, \zeta_{1}-\zeta_{2})v;
 z_{2}, \zeta_{2})w\rangle
\end{eqnarray}
holds when
$|z_{1}|>|z_{2}|>0$ and $|\zeta_{1}|>|\zeta_{2}|>0$.

\item The {\it single-valuedness property}: 
\begin{equation} \label{sing-val-1}
e^{2\pi i (L^L(0)-L^R(0))} = I_{F}.
\end{equation}   

\item The {\it skew symmetry}: 
\begin{equation} \label{skew-1-1}
\mathbb{Y}^{\Y}(u; 1, 1)v = e^{L^L(-1)+L^R(-1)}
\mathbb{Y}^{\Y}(v; e^{\pi i}, e^{-\pi i})u.
\end{equation}
\end{enumerate}
\end{thm}

Let $V$ be a simple vertex operator algebra satisfying the following:

\begin{cond}[uniqueness of vacuum]\label{u-o-v}
{\rm For $n<0$, $V_{(n)}=0$ and $V_{(0)}=\mathbb{C}\mathbf{1}$
and 
for irreducible $V$-module $W$ not isomorphic to $V$, $W_{(0)}=0$.}
\end{cond}

\begin{cond}[complete reducibility]\label{c-red}
{\rm Every $\mathbb{N}$-gradable weak 
$V$-module is completely reducible.}
\end{cond}

\begin{cond}[$C_{2}$-cofiniteness]\label{c-2-c-f}
{\rm $V$ is $C_{2}$-cofinite, that is, $\dim V/C_{2}(V)<\infty$, 
where $C_{2}(V)$ is the subspace of $V$ spanned by elements of 
the form $u_{-2}v$ for $u, v\in V$.}
\end{cond}

We now recall the construction of conformal 
full field algebras over $V\otimes V$
in \cite{HK}. 

Let $\A$ be the set of equivalence classes of irreducible $V$-modules. For 
any $a\in \A$, we choose a representative $W^{a}$ from $a$. 
Then there exists $h_{a}\in \R$ such that $W^{a}=\coprod_{n\in \Z}
W^{a}_{(n+h_{a})}$.

For a single-valued
branch $f_{1}(z_{1}, z_{2})$ of a multivalued
analytic function in a region $A$, we use $E(f_{1}(z_{1}, z_{2}))$ 
to denote the multivalued analytic extension together with the 
preferred branch $f_{1}(z_{1}, z_{2})$. 
Let $w_{1}=w_{1}(z_{1}, z_{2})$ and $w_{2}=w_{2}(z_{1}, z_{2})$
be a change of variables and $f_{2}(z_{1}, z_{2})$ 
a branch 
of $E(f_{1}(z_{1}, z_{2}))$ in a region $B$ containing
$w_{1}(z_{1}, z_{2})=0$
and $w_{2}(z_{1}, z_{2})=0$ such that $A\cap B\ne \emptyset$ and 
$f_{1}(z_{1}, z_{2})=f(z_{1}, z_{2})$ for $(z_{1}, z_{2})\in A\cap B$.
Then we use 
$$\res_{w_{1}=0\;|\;w_{2}}E(f_{1}(z_{1}, z_{2}))$$
to denote the coefficient of $w_{1}^{-1}$ in the expansion of 
$f_{2}(z_{1}, z_{2})$ as a series in powers of $w_{1}$  whose coefficients
are analytic functions of $w_{2}$.

For $a_{1}, a_{2}, a_{3}\in \A$, $w_{a_{1}}\in W^{a_{1}}$, 
$w_{a_{2}}\in W^{a_{2}}$, 
$w_{a'_{1}}\in (W^{a_{1}})'$, $w'_{a_{2}}\in (W^{a_{1}})'$, 
$\Y_{1}\in 
\V_{a_{1}a_{2}}^{a_{3}}$ and $\Y_{2}\in \V_{a'_{1}a'_{2}}^{a'_{3}}$, 
Let $\langle \Y_{1},
\Y_{2}\rangle_{\V_{a_{1}a_{2}}^{a_{3}}}\in \C$ be given by
\begin{eqnarray}\label{inner}
\lefteqn{\res_{1-z_{1}-z_{2}=0\;|\;z_{2}}(1-z_{1}-z_{2})^{-1}
E(\langle e^{L(1)}\Y_{2}((1-z_{1}-z_{2})^{L(0)}
\tilde{w}'_{a_{1}}, z_{1})
\tilde{w}'_{a_{2}}, }\nn
&&\quad\quad\quad\quad\quad\quad\quad\quad\quad\quad
e^{L(1)}\Y_{1}((1-z_{1}-z_{2})^{L(0)}\tilde{w}_{a_{1}}, z_{2})
\tilde{w}_{a_{2}}\rangle)\nn
&&=\langle w'_{a_{1}}, w_{a_{1}}\rangle\langle w'_{a_{2}}, w_{a_{2}}\rangle
\langle \Y_{1},
\Y_{2}\rangle_{\V_{a_{1}a_{2}}^{a_{3}}}.
\end{eqnarray}
It was shown in \cite{HK} that $\langle \Y_{1},
\Y_{2}\rangle_{\V_{a_{1}a_{2}}^{a_{3}}}$ indeed exists. 
Clearly, $\langle \Y_{1},
\Y_{2}\rangle_{\V_{a_{1}a_{2}}^{a_{3}}}$ is bilinear in $\Y_{1}$ and $\Y_{2}$.
Thus we have a pairing $\langle \cdot, \cdot\rangle_{\V_{a_{1}a_{2}}^{a_{3}}}: 
\V_{a_{1}a_{2}}^{a_{3}}\otimes \V_{a'_{1}a'_{2}}^{a'_{3}}
\to \C$. The following is another main result of \cite{HK}:

\begin{thm}
The pairing $\langle \cdot, \cdot\rangle_{\V_{a_{1}a_{2}}^{a_{3}}}: 
\V_{a_{1}a_{2}}^{a_{3}}\otimes \V_{a'_{1}a'_{2}}^{a'_{3}}
\to \C$ is nondegenerate. In particular, 
$N_{a'_{1}a'_{2}}^{a'_{3}}
=N_{a_{1}a_{2}}^{a_{3}}.$
\end{thm}

Recall that in \cite{H10},  an 
action of $S_{3}$ on the space 
of intertwining operators for a vertex operator algebras was given.
We choose a canonical 
basis of $\V_{a_{1}a_{2}}^{a_{3}}$ for $a_{1}, a_{2}, a_{3}\in
\A$ when one of $a_{1}, a_{2}, a_{3}$ is $e$: For $a\in \A$, we 
choose $\Y_{ea; 1}^{a}$ to 
be the vertex operator 
$Y_{W^{a}}$ defining the module structure on $W^{a}$ and we choose 
$\Y_{ae; 1}^{a}$ to be the intertwining operator defined using 
the action of $\sigma_{12}$ on choose $\Y_{ea; 1}^{a}$, 
or equivalently, using the skew-symmetry 
in this case,
\begin{eqnarray*}
\Y_{ae; 1}^{a}(w_{a}, x)u&=&\sigma_{12}(\Y_{ea; 1}^{a})(w_{a}, x)u\nn
&=&e^{xL(-1)}\Y_{ea; 1}^{a}(u, -x)w_{a}\nn
&=&e^{xL(-1)}Y_{W^{a}}(u, -x)w_{a}
\end{eqnarray*}
for $u\in V$ and $w_{a}\in W^{a}$. 
Since $V'$ as the contragredient of the irreducible adjoint module 
$V$ is an irreducible $V$-module (see \cite{FHL} 
and has a nonzero homogeneous 
subspace of weight $0$, 
as a $V$-module it must be isomorphic to $V$. So we have 
$e'=e$. From \cite{FHL}, we know that there is a nondegenerate
invariant 
bilinear form $(\cdot, \cdot)$ on $V$ such that $(\mathbf{1}, 
\mathbf{1})=1$. 
We choose $\Y_{aa'; 1}^{e}=\Y_{aa'; 1}^{e'}$
to be the intertwining operator defined using the action of 
$\sigma_{23}$ by
$$\Y_{aa'; 1}^{e'}=\sigma_{23}(\Y_{ae; 1}^{a}),$$
that is,
$$(u, \Y_{aa'; 1}^{e}(w_{a}, x)w_{a'})
=e^{\pi i h_{a}}\langle \Y_{ae; 1}^{a}(e^{xL(1)}(e^{-\pi i}x^{-2})^{L(0)}w_{a}, x^{-1})u, 
w_{a'}\rangle$$
for $u\in V$, $w_{a}\in W^{a}$ and $w_{a'}\in W^{a'}$. Since the actions of
$\sigma_{12}$
and $\sigma_{23}$ generate the action of $S_{3}$ on $\mathcal{V}$, we have
$$\Y_{a'a; 1}^{e}=\sigma_{12}(\Y_{aa'; 1}^{e})$$
for any $a\in \mathcal{A}$.

As in \cite{HK}, 
for $a\in \A$, let 
$$F_{a}=F(\Y_{ae; 1}^{a} \otimes \Y_{a'a; 1}^{e};
\Y_{ea; 1}^{a}\otimes \Y_{aa'; 1}^{e})\ne 0$$
and we use $\sqrt{F_{a}}$ to denote the 
square root $\sqrt{|F_{a}|}e^{\frac{i\arg F_{a}}{2}}$
of $F_{a}$. 
For $a_{1}, a_{2}, a_{3}\in \A$,
consider the modified pairings 
$$\frac{\sqrt{F_{a_{3}}}}{\sqrt{F_{a_{1}}}\sqrt{F_{a_{2}}}}
\langle \cdot, \cdot\rangle_{\V_{a_{1}a_{2}}^{a_{3}}}.$$
These pairings give a nondegenerate bilinear form 
$(\cdot, \cdot)_{\V}$ on $\V$. 
For any basis $\{\Y_{a_{1}a_{2}; i}^{a_{3}}\;
|\;i=\dots, N_{a_{1}a_{2}}^{a_{3}}\}$ of 
$\V_{a_{1}a_{2}}^{a_{3}}$ and 
any $\sigma\in S_{3}$, $\{\sigma(\Y_{a_{1}a_{2}; i}^{a_{3}})\;
|\;i=\dots, N_{a_{1}a_{2}}^{a_{3}}\}$ is a basis of 
$\sigma(\V_{a_{1}a_{2}}^{a_{3}})$.

We have the following result from \cite{HK}:

\begin{prop}\label{skew}
The nondegenerate bilinear form $(\cdot, \cdot)_{\V}$
is invariant with respect to the action of $S_{3}$ on $\V$, that is,
for $a_{1}, a_{2}, a_{3}\in \A$, $\sigma\in S_{3}$, $\Y_{1}\in 
\V_{a_{1}a_{2}}^{a_{3}}$
and $\Y_{2}\in \V_{a'_{1}a'_{2}}^{a'_{3}}$, 
$$(\sigma(\Y_{1}), \sigma(\Y_{2}))_{\V}= 
(\Y_{1}, \Y_{2})_{\V}.$$
Equivalently, for $a_{1}, a_{2}, a_{3}\in \A$,
$$\left\{\left.\frac{\sqrt{F_{a_{1}}}\sqrt{F_{a_{2}}}
\sqrt{F_{a_{\sigma^{-1}(3)}}}}
{\sqrt{F_{a_{\sigma^{-1}(1)}}}\sqrt{F_{a_{\sigma^{-1}(2)}}}\sqrt{F_{a_{3}}}}
\sigma(\Y_{a'_{1}a'_{2}; i}^{\prime; a'_{3}})\;\right|\;i=\dots, 
N_{a_{1}a_{2}}^{a_{3}}\right\}$$
is the dual basis of $\{\sigma(\Y_{a_{1}a_{2}; i}^{a_{3}})\;
|\;i=\dots, N_{a_{1}a_{2}}^{a_{3}}\}$, where 
$\{\Y_{a_{1}a_{2}; i}^{a_{3}}\;
|\;i=\dots, N_{a_{1}a_{2}}^{a_{3}}\}$ is a basis of 
$\V_{a_{1}a_{2}}^{a_{3}}$ and $\{\Y_{a'_{1}a'_{2}; i}^{\prime;a'_{3}}\;
|\;i=\dots, N_{a_{1}a_{2}}^{a_{3}}\}$ is its dual basis 
with respect to the pairing 
$\langle \cdot, \cdot\rangle_{\V_{a_{1}a_{2}}^{a_{3}}}$.
\end{prop}

Let 
$$F=\oplus_{a\in \A} W^{a}\otimes W^{a'}.$$
For $w_{a_{1}}\in W^{a_{1}}$, $w_{a_{2}}\in W^{a_{2}}$, 
$w_{a'_{1}}\in W^{a'_{1}}$ and $w_{a'_{2}}\in W^{a'_{2}}$, 
we define 
\begin{eqnarray*}
\lefteqn{\mathbb{Y}((w_{a_{1}}\otimes w_{a'_{1}}); z, \overline{z})
(w_{a_{2}}\otimes w_{a'_{2}})}\nn
&&=\sum_{a_{3}\in \A}\sum_{p=1}^{N_{a_{1}a_{2}}^{a_{3}}}
\Y_{a_{1}a_{2}; p}^{a_{3}}(w_{a_{1}}, z) w_{a_{2}}
\otimes \Y_{a'_{1}a'_{2}; p}^{\prime; a'_{3}}
(w_{a'_{1}}, \overline{z})w_{a'_{2}}.
\end{eqnarray*}

In \cite{HK}, we proved the following result:

\begin{thm}\label{full-const}
The quadruple $(F, \mathbb{Y}, \one\otimes \one, \omega\otimes \one,
\one \otimes \omega)$ is a conformal full field algebra over 
$V\otimes V$. 
\end{thm}

\renewcommand{\theequation}{\thesection.\arabic{equation}}
\renewcommand{\thethm}{\thesection.\arabic{thm}}
\setcounter{equation}{0}
\setcounter{thm}{0}

\section{Modular invariance for intertwining operator algebras}

In this section, we review the modular invariance of intertwining 
operator algebras proved in \cite{H7}. We assume in this section that 
$V$ is  a simple vertex operator algebra satisfying Conditions 
\ref{c-red} and \ref{c-2-c-f} and the condition
that for $n<0$, $V_{(n)}=0$ and $V_{(0)}=\mathbb{C}\mathbf{1}$.
Note that the last condition is weaker than Condition \ref{u-o-v}.
We shall 
use $Y$ to denote the vertex operator maps for the algebra
$V$ and for $V$-modules.

Let $A_{j}$, $j\in \mathbb{Z}_{+}$, be the 
complex numbers defined 
by 
\begin{eqnarray*}
\frac{1}{2\pi i}\log(1+2\pi i y)=\left(\exp\left(\sum_{j\in \mathbb{Z}_{+}}
A_{j}y^{j+1}\frac{\partial}{\partial y}\right)\right)y.
\end{eqnarray*} 
For any $V$-module 
$W$, we shall denote the operator 
$\sum_{j\in \mathbb{Z}_{+}}
A_{j}L(j)$
on $W$ by $L_{+}(A)$. Then 
$$e^{\displaystyle -\sum_{j\in \mathbb{Z}_{+}}
A_{j}L(j)}=e^{-L_{+}(A)}.$$
Let
$$\mathcal{U}(x)=(2\pi i)^{L(0)}x^{L(0)}e^{-L^{+}(A)}
\in (\mbox{\rm End}\;W)\{x\}$$
where $(2\pi i)^{L(0)}=e^{(\log 2\pi +i \frac{\pi}{2})L(0)}$.
Let $B_{j}\in \Q$ for $j\in \Z_{+}$ be defined by 
\begin{eqnarray*}
\log(1+y)=\left(\exp\left(\sum_{j\in \mathbb{Z}_{+}}
B_{j}y^{j+1}\frac{\partial}{\partial y}\right)\right)y.
\end{eqnarray*} 
Then it is easy to see that 
$$\mathcal{U}(x)=x^{L(0)}e^{-L^{+}(B)}(2\pi i)^{L(0)}.$$

For any $z\in \C$, we shall denote $e^{2\pi iz}$ by $q_{z}$
and we shall also use $\mathcal{U}(q_{z})$ to denote the 
the map obtained by substituting $e^{2\pi izL(0)}$ for 
$x^{L(0)}$ in $\mathcal{U}(x)$, that is, 
\begin{eqnarray}\label{ux-1}
\mathcal{U}(q_{z})&=&(2\pi i)^{L(0)}e^{2\pi izL(0)}e^{-L^{+}(A)}\nn
&=&e^{2\pi izL(0)}e^{-L^{+}(B)}(2\pi i)^{L(0)}.
\end{eqnarray}

For $V$-modules $W_{i}$
and $\tilde{W}_{i}$, $i=1, \dots, n$,  intertwining operators
$\mathcal{Y}_{i}$, $i=1, \dots, n$,  of 
types $\binom{\tilde{W}_{i-1}}{ W_{i}\tilde{W}_{i}}$, respectively,
where we use the convention $\tilde{W}_{0}=
\tilde{W}_{n}$, and  $w_{i}\in W_{i}$, $i=1, \dots, n$,
we shall consider the element
\begin{eqnarray}\label{correl-fn}
\lefteqn{(F_{\mathcal{Y}_{1}, \dots, \mathcal{Y}_{n})}(w_{1}, \dots, w_{n};
z_{1}, \dots, z_{n}; q)}\nn
&&=\tr_{\tilde{W}_{n}}\mathcal{Y}_{n}(\mathcal{U}(q_{z_{1}})w_{1}, q_{z_{1}})
\cdots
\mathcal{Y}_{1}(\mathcal{U}(q_{z_{n}})w_{n}, q_{z_{n}})
q^{L(0)-\frac{c}{24}}
\end{eqnarray}
of $\mathbb{G}_{|q_{z_{1}}|>\cdots>|q_{z_{n}}|>0}((q))$,
where for complex variables $\xi_{1}, \dots, \xi_{n}$,
$\mathbb{G}_{|\xi_{1}|>\cdots>|\xi_{n}|>0}$ is the space
of all multivalued analytic functions in $\xi_{1}, \dots, \xi_{n}$
defined on the region
$|\xi_{1}|>\cdots>|\xi_{n}|>0$ with preferred branches 
in the simply-connected region $|\xi_{1}|>\cdots>|\xi_{n}|>0$,
$0\le \arg \xi_{i}<2\pi$,
$i=1, \dots, n$.

In \cite{H7}, the following result was proved:

\begin{thm}\label{conv-ext}
In the region $1>|q_{z_{1}}|>\cdots 
>|q_{z_{n}}|>|q_{\tau}|>0$, the series
\begin{equation}\label{correl-fn-tau}
F_{\mathcal{Y}_{1}, \dots, \mathcal{Y}_{n}}(w_{1}, 
\dots, w_{n};
z_{1}, \dots, z_{n}; q_{\tau})
\end{equation}
is absolutely convergent and can be analytically extended
to a (multivalued) analytic function in the region given by
$\Im(\tau)>0$ (here $\Im(\tau)$ is the imaginary part of $\tau$), 
$z_{i}\ne z_{j}+k\tau+l$ for $i, j=1, \dots, n$,
$i\ne j$, $k, l\in \mathbb{Z}$.
\end{thm}

We shall denote the (multivalued) analytic extension given in Theorem 
\ref{conv-ext} by 
$$(\Phi(\mathcal{Y}_{1}\otimes \cdots\otimes \mathcal{Y}_{n}))(w_{1}, 
\dots, w_{n};
z_{1}, \dots, z_{n}; q_{\tau}).$$
Note that in \cite{H7}, this function is denoted by 
$$\overline{F}_{\mathcal{Y}_{1}, \dots, \mathcal{Y}_{n}}(w_{1}, 
\dots, w_{n};
z_{1}, \dots, z_{n}; q_{\tau}).$$
Here we change the notation
to avoid confusions with  notations related to algebraic extensions or
to complex conjugations and also for convenience in later sections
and for future use. 

In \cite{H7}, the following genus-one duality results were proved:

\begin{thm}[Genus-one commutativity]
Let $W_{i}$ and $\tilde{W}_{i}$ 
be $V$-modules and
$\mathcal{Y}_{i}$ intertwining operators  of 
types $\binom{\tilde{W}_{i-1}}{W_{i}\tilde{W}_{i}}$ ($i=1, \dots, n$,
$\tilde{W}_{0}=\tilde{W}_{n}$), 
respectively. Then for any $1\le k\le n-1$, there exist
$V$-modules $\check{W}_{k}$ and intertwining operators
$\check{\mathcal{Y}}_{k}$ and $\check{\mathcal{Y}}_{k+1}$
of types $\binom{\check{W}_{k}}{W_{k}\tilde{W}_{k+1}}$
and $\binom{\tilde{W}_{k-1}}{W_{k+1}\check{W}_{k}}$, respectively,
such that 
$$F_{\mathcal{Y}_{1}, \dots, \mathcal{Y}_{n}}(w_{1}, 
\dots, w_{n};
z_{1}, \dots, z_{n}; q_{\tau})$$
and 
\begin{eqnarray*}
\lefteqn{F_{\mathcal{Y}_{1}, \dots, \mathcal{Y}_{k-1},
\check{\mathcal{Y}}_{k+1}, \check{\mathcal{Y}}_{k}, \mathcal{Y}_{k+2}
\dots, \mathcal{Y}_{n}}(w_{1}, 
\dots, w_{k-1}, w_{k+1}, w_{k}, w_{k+2}, \dots, w_{n};}\nn
&&\quad\quad\quad\quad\quad\quad\quad\quad\quad\quad\quad\quad
\quad\quad\quad
z_{1}, \dots, z_{k-1}, z_{k+1}, z_{k}, z_{k+2}, 
\dots, z_{n}; q_{\tau})
\end{eqnarray*}
are analytic extensions of each other, or equivalently, 
\begin{eqnarray*}
\lefteqn{(\Phi(\mathcal{Y}_{1}\otimes \cdots\otimes \mathcal{Y}_{n}))(w_{1}, 
\dots, w_{n};
z_{1}, \dots, z_{n}; \tau)}\nn
&&=(\Phi(\mathcal{Y}_{1}\otimes\cdots \otimes\mathcal{Y}_{k-1}\otimes
\check{\mathcal{Y}}_{k+1}\otimes \check{\mathcal{Y}}_{k}\otimes \mathcal{Y}_{k+2}
\otimes\cdots\otimes \mathcal{Y}_{n}))(w_{1},
\dots, w_{k-1}, \nn
&&\quad\quad\quad\quad\quad w_{k+1}, w_{k}, w_{k+2}, \dots, w_{n};
z_{1}, \dots, z_{k-1}, z_{k+1}, z_{k}, z_{k+2}, 
\dots, z_{n}; \tau).
\end{eqnarray*}
More generally, for any $\sigma\in S_{n}$, 
there exist $V$-modules $\check{W}_{i}$ ($i=1, \dots, n$) and 
intertwining operators $\check{\mathcal{Y}}_{i}$ of 
types $\binom{\check{W}_{i-1}}{W_{\sigma(i)}\check{W}_{i}}$ ($i=1, \dots, n$,
$\check{W}_{0}=\check{W}_{n}=\tilde{W}_{n}$), respectively, 
such that 
$$F_{\mathcal{Y}_{1}, \dots, \mathcal{Y}_{n}}(w_{1}, 
\dots, w_{n};
z_{1}, \dots, z_{n}; q_{\tau})$$
and 
$$F_{\check{\mathcal{Y}}_{1}, \dots, \check{\mathcal{Y}}_{n}}(w_{\sigma(1)}, 
\dots, w_{\sigma(n)};
z_{\sigma(1)}, \dots, z_{\sigma(n)}; q_{\tau})$$ 
are analytic extensions of each other, or equivalently,
\begin{eqnarray*}
\lefteqn{(\Phi(\mathcal{Y}_{1}\otimes \cdots\otimes \mathcal{Y}_{n}))(w_{1}, 
\dots, w_{n};
z_{1}, \dots, z_{n}; \tau)}\nn
&&=(\Phi(\check{\mathcal{Y}}_{1}\otimes \cdots\otimes
\check{\mathcal{Y}}_{n}))(w_{\sigma(1)}, 
\dots, w_{\sigma(n)};
z_{\sigma(1)}, \dots, z_{\sigma(n)}; \tau).
\end{eqnarray*}
\end{thm}

\begin{thm}[Genus-one associativity]\label{g1-asso}
Let $W_{i}$ and $\tilde{W}_{i}$ for $i=1, \dots, n$ be 
$V$-modules and
$\mathcal{Y}_{i}$ intertwining operators  of 
types $\binom{\tilde{W}_{i-1}}{W_{i}\tilde{W}_{i}}$ ($i=1, \dots, n$,
$\tilde{W}_{0}=\tilde{W}_{n}$), 
respectively. Then for any $1\le k\le n-1$, 
there exist a $V$-module $\check{W}_{k}$ and 
intertwining operators $\check{\mathcal{Y}}_{k}$ and 
$\check{\mathcal{Y}}_{k+1}$ of 
types $\binom{\check{W}_{k}}{W_{k} W_{k+1}}$ and 
$\binom{\tilde{W}_{k-1}}{\check{W}_{k}\tilde{W}_{k+1}}$, respectively, 
such that 
\begin{eqnarray}\label{g-1-iter}
\lefteqn{(\Phi(\mathcal{Y}_{1}\otimes \cdots\otimes \mathcal{Y}_{k-1}\otimes
\check{\mathcal{Y}}_{k+1}\otimes \mathcal{Y}_{k+2}\otimes \cdots\otimes
\mathcal{Y}_{n}))(w_{1}, 
\dots, w_{k-1}, }\nn
&&\quad\quad\quad
\check{\mathcal{Y}}(w_{k}, z_{k}-z_{k+1})
w_{k+1}, 
w_{k+2}, \dots, w_{n};
z_{1}, \dots, z_{k-1}, z_{k+1}, \dots, z_{n}; \tau)\nn
&&=\sum_{r\in \mathbb{R}}
(\Phi(\mathcal{Y}_{1}\otimes \cdots\otimes \mathcal{Y}_{k-1}\otimes
\check{\mathcal{Y}}_{k+1}\otimes \mathcal{Y}_{k+2}\otimes \cdots\otimes
\mathcal{Y}_{n}))(w_{1}, 
\dots, w_{k-1}, \nn
&&\quad\quad
P_{r}(\check{\mathcal{Y}}(w_{k}, z_{k}-z_{k+1})
w_{k+1}), 
w_{k+2}, \dots, w_{n};
z_{1}, \dots, z_{k-1}, z_{k+1}, \dots, z_{n}; \tau)\nn
&&
\end{eqnarray}
is absolutely convergent when $1>|q_{z_{1}}|>\cdots 
>|q_{z_{k-1}}|>|q_{z_{k+1}}|>
\dots >|q_{z_{n}}|>|q_{\tau}|>0$ and 
$1>|q_{(z_{k}-z_{k+1})}-
1|>0$
and is convergent to 
$$(\Phi(\mathcal{Y}_{1}\otimes \cdots\otimes \mathcal{Y}_{n}))(w_{1}, 
\dots, w_{n};
z_{1}, \dots, z_{n}; \tau)$$ 
when $1>|q_{ z_{1}}|>\cdots 
>|q_{z_{n}}|>|q_{\tau}|>0$ and 
$|q_{(z_{k}-z_{k+1})}|>1>|q_{(z_{k}-z_{k+1})}-
1|>0$.
\end{thm}

Let $W_{i}$ be $V$-modules and $w_{i}\in W_{i}$ for $i=1, \dots, n$.
For  any $V$-modules
$\tilde{W}_{i}$ and any intertwining operators
$\mathcal{Y}_{i}$, $i=1, \dots, n$,  of 
types $\binom{\tilde{W}_{i-1}}{W_{i}\tilde{W}_{i}}$, respectively,
we have a genus-one correlation function
$$(\Phi(\mathcal{Y}_{1}\otimes \cdots\otimes \mathcal{Y}_{n}))(w_{1}, \dots, w_{n};
z_{1}, \dots, z_{n}; \tau).$$
Note that these multivalued functions actually have preferred
branches in the region $1>|q_{z_{1}}|>\cdots >|q_{z_{n}}|>|q_{\tau}|>0$
given by the intertwining operators $\mathcal{Y}_{1}, 
\dots, \mathcal{Y}_{n}$. Thus linear 
combinations of these functions make sense. 
For fixed $V$-modules $W_{i}$ and $w_{i}\in W_{i}$
for $i=1, \dots, n$, denote the vector space
spanned by all such functions by  $\mathcal{F}_{w_{1}, \dots, w_{n}}$. 
For any single valued analytic function $f(\tau)$ of $\tau$ 
and any $r\in \R$, 
we choose a branch of the multivalued 
analytic function $f(\tau)^{r}$
to be $e^{r\log f(\tau)}$. 
The following theorem is one of the main result of \cite{H7}:

\begin{thm}\label{mod-inv}
For any $V$-modules $\tilde{W}_{i}$
and any
intertwining operators $\mathcal{Y}_{i}$ ($i=1, \dots, n$)  of 
types $\binom{\tilde{W}_{i-1}}{W_{i}\tilde{W}_{i}}$, respectively,
and any 
$$\left(\begin{array}{cc}
\alpha&\beta\\
\gamma&\delta
\end{array}\right)\in SL(2, \mathbb{Z}),$$
\begin{eqnarray*}
\lefteqn{(\Phi(\mathcal{Y}_{1}\otimes \cdots\otimes \mathcal{Y}_{n}))
\Biggl(\left(\frac{1}{\gamma\tau+\delta}\right)^{L(0)}w_{1}, \dots,
\left(\frac{1}{\gamma\tau+\delta}\right)^{L(0)}w_{n};}\nn
&&\quad\quad\quad\quad\quad\quad\quad\quad\quad\quad\quad\quad
\quad\quad\quad
\frac{z_{1}}{\gamma\tau+\delta}, \dots, \frac{z_{n}}{\gamma\tau+\delta}; 
\frac{\alpha\tau+\beta}{\gamma\tau+\delta}\Biggr)
\end{eqnarray*}
is in $\mathcal{F}_{w_{1}, \dots, w_{n}}$. 
\end{thm}

\renewcommand{\theequation}{\thesection.\arabic{equation}}
\renewcommand{\thethm}{\thesection.\arabic{thm}}
\setcounter{equation}{0}
\setcounter{thm}{0}

\section{Genus-one correlation functions and modular invariance 
for conformal full field algebras}

Let $V^{L}$ and $V^{R}$ be simple vertex operator algebras satisfying 
Conditions \ref{u-o-v}, \ref{c-red} and \ref{c-2-c-f}.
Let $F$ be a conformal
full field algebra over $V^L\otimes V^R$. In particular, $F$ 
is $\R\times \R$-graded, that is,
$F=\coprod_{r, s\in \R}F_{(r, s)}$ where $F_{(r, s)}$ for 
$r, s\in \R$ are eigenspaces for the operators $L^{L}(0)$ and
$L^{R}(0)$ with eigenvalues $r$ and $s$, respectively. 
For a linear map $f: F\to \overline{F}$, we define the 
$q$-$\overline{q}$-trace 
of $f$ to be 
$$\tr_{F}fq^{L^{L}(0)-\frac{c^{L}}{24}}
\overline{q}^{L^{R}(0)-\frac{c^{R}}{24}}
=\sum_{r, s\in \R}\tr_{F_{(m, n)}}
fq^{r-\frac{c^{L}}{24}}\overline{q}^{s-\frac{c^{R}}{24}}.$$

As in \cite{HK}, we choose the branches of the functions 
$z^{r}$ and $\overline{z}^{s}$ for $r, s\in \R$ to be 
$e^{r\log z}$ and $e^{s\overline{\log z}}$, respectively. 
On the other hand, for the functions 
$(e^{2\pi i z})^{r}$ and $\left(\overline{e^{2\pi i z}}\right)^{s}
=(e^{-2\pi i \overline{z}})^{s}$ for $r, s\in \R$, we choose
their branches to be 
$e^{2\pi irz}$ and $e^{-2\pi is\overline{z}}$, respectively. 
For multivalued analytic functions obtained from these functions
using products and sums, we choose their
branches to be the ones obtained from the branches we choose above 
using the same operations.

Now let 
\begin{eqnarray*}
\mathcal{U}^{L}(e^{2\pi iz})
&=&(2\pi i)^{L^{L}(0)}(e^{2\pi iz})^{L^{L}(0)}e^{-L_{+}^{L}(A)}\nn
&=&(2\pi i)^{L^{L}(0)}e^{2\pi izL^{L}(0)}e^{-L_{+}^{L}(A)},\nn
\mathcal{U}^{R}(e^{2\pi iz})
&=&(2\pi i)^{L^{R}(0)}(e^{2\pi iz})^{L^{R}(0)}e^{-L_{+}^{R}(A)}\nn
&&(2\pi i)^{L^{R}(0)}e^{2\pi izL^{R}(0)}e^{-L_{+}^{R}(A)},
\end{eqnarray*}
where $L_+^L(A)=\sum_{j\in \Z_+} A_j L^L(j)$ and 
$L_+^R(A)=\sum_{j\in \Z_+} A_j L^R(j)$.
By (\ref{ux-1}) and our choice of the branches of the functions 
$(e^{2\pi i z})^{r}$ and $\left(\overline{e^{2\pi i z}}\right)^{s}
=(e^{-2\pi i \overline{z}})^{s}$ above, we have
\begin{eqnarray*}
\overline{\mathcal{U}^{R}(e^{2\pi iz})}
&=&\overline{(e^{2\pi iz})^{L^{R}(0)}e^{-L_{+}^{L}(B)}(2\pi i)^{L^{R}(0)}}\nn
&=&\overline{(e^{2\pi iz})^{L^{R}(0)}e^{-L_{+}^{L}(B)}(2\pi)^{L^{R}(0)}
e^{\frac{\pi i}{2}L^{R}(0)}}\nn
&=&(\overline{e^{-2\pi iz}})^{L^{R}(0)}e^{-L_{+}^{L}(B)}(2\pi)^{L^{R}(0)}
e^{-\frac{\pi i}{2}L^{R}(0)}\nn
&=&(e^{-2\pi i\overline{z}})^{L^{R}(0)}e^{-L_{+}^{L}(B)}(2\pi)^{L^{R}(0)}
e^{\frac{\pi i}{2}L^{R}(0)}e^{-\pi iL^{R}(0)}\nn
&=&\mathcal{U}^{R}(e^{-2\pi i\overline{z}})e^{-\pi iL^{R}(0)}.
\end{eqnarray*}

For any $u\in F$ and $z\in \C^{\times}$, we shall call the operator
$$\mathbb{Y}(\mathcal{U}^{L}(e^{2\pi iz})
\overline{\mathcal{U}^{R}(e^{2\pi iz})}
u; e^{2\pi iz_{1}}, \overline{e^{2\pi iz}}): F\to \overline{F}$$
a {\it geometrically-modified full vertex operator}.
For $u_{1}, \dots, u_{k}\in F$ and $z_{1}, \dots, z_{k}\in 
\C$ satisfying $|e^{2\pi iz_{1}}|>\cdots>|e^{2\pi iz_{k}} |>0$, 
the product 
\begin{eqnarray*}
\lefteqn{\mathbb{Y}(\mathcal{U}^{L}(e^{2\pi iz_{1}})
\overline{\mathcal{U}^{R}(e^{2\pi iz_{1}})}
u_{1}; e^{2\pi iz_{1}}, \overline{e^{2\pi iz_{1}}})\cdot}\nn
&&\quad\quad\quad\cdots \mathbb{Y}(\mathcal{U}^{L}(e^{2\pi iz_{k}})
\overline{\mathcal{U}^{R}(e^{2\pi iz_{k}})}
u_{1}; e^{2\pi iz_{k}}, \overline{e^{2\pi iz_{k}}})\nn
&&=\mathbb{Y}(\mathcal{U}^{L}(e^{2\pi iz_{1}})
\mathcal{U}^{R}(e^{-2\pi i\overline{z}_{1}})e^{-\pi iL^{R}(0)}
u_{1}; e^{2\pi iz_{1}}, e^{-2\pi i\overline{z}_{1}})\cdot\nn
&&\quad\quad\quad\cdots \mathbb{Y}(\mathcal{U}^{L}(e^{2\pi iz_{k}})
\mathcal{U}^{R}(e^{-2\pi i\overline{z}_{k}})e^{-\pi iL^{R}(0)}
u_{1}; e^{2\pi iz_{k}}, e^{-2\pi i\overline{z}_{k}})\nn
\end{eqnarray*}
of geometrically-modified full vertex operators
is a linear map from $F$ to $\overline{F}$. So
we have the $q$-$\overline{q}$-trace
\begin{eqnarray}\label{q-tr-fvo}
\lefteqn{\tr_{F}\mathbb{Y}(\mathcal{U}^{L}(e^{2\pi iz_{1}})
\overline{\mathcal{U}^{R}(e^{2\pi iz_{1}})}
u_{1}; e^{2\pi iz_{1}}, \overline{e^{2\pi iz_{1}}})\cdot}\nn
&&\quad\quad\quad\cdots \mathbb{Y}(\mathcal{U}^{L}(e^{2\pi iz_{1}}) 
\overline{\mathcal{U}^{R}(e^{2\pi iz_{k}})}
u_{1}; e^{2\pi iz_{1}}, \overline{e^{2\pi iz_{k}}})
q^{L^{L}(0)-\frac{c^{L}}{24}}\overline{q}^{L^{R}(0)-\frac{c^{R}}{24}}\nn
&&=\tr_{F}\mathbb{Y}(\mathcal{U}^{L}(e^{2\pi iz_{1}})
\mathcal{U}^{R}(e^{-2\pi i\overline{z}_{1}})e^{-\pi iL^{R}(0)}
u_{1}; e^{2\pi iz_{1}}, e^{-2\pi i\overline{z}_{1}})\cdot\nn
&&\quad\quad\quad\cdots \mathbb{Y}(\mathcal{U}^{L}(e^{2\pi iz_{1}}) 
\mathcal{U}^{R}(e^{-2\pi i\overline{z}_{k}})e^{-\pi iL^{R}(0)}
u_{1}; e^{2\pi iz_{1}}, e^{-2\pi i\overline{z}_{k}})\cdot\nn
&&\quad\quad\quad\quad\quad\quad\cdot
q^{L^{L}(0)-\frac{c^{L}}{24}}\overline{q}^{L^{R}(0)-\frac{c^{R}}{24}}.
\end{eqnarray}

As a module 
for the vertex operator algebra $V^{L}\otimes V^{R}$, $F$ is a 
direct sum of irreducible modules for $V^{L}\otimes V^{R}$. Let 
$\A^{L}$ ($\A^{R}$) be the set of equivalence classes of irreducible 
modules for $V^{L}$ ($V^{R}$). For each $a^{L}\in \A^{L}$ 
($a^{R}\in \A^{R}$), choose a representative $W^{a^{L}}$ 
($W^{a^{R}}$) of $a^{L}$ ($a^{R}$). Then there 
exist a positive integer $N$ and maps $r^{L}: \{1, \dots, N\}
\to \A^{L}$, $r^{R}: \{1, \dots, N\}
\to \A^{R}$ such that $F$ is isomorphic as a $V^{L}\otimes V^{R}$-module
to $\coprod_{n=1}^{N}W^{r^{L}(n)}\otimes W^{r^{R}(n)}$. 
We now shall identify the vector space 
$F$ with this $V^{L}\otimes V^{R}$-module. 
Then the full vertex operator map can be written as 
\begin{equation}\label{full-vo-decomp}
\mathbb{Y}=\sum_{l, m, n=1}^{N}
\sum_{i=1}^{N_{r^{L}(m)r^{L}(n)}^{r^{L}(l)}}
\sum_{j=1}^{N_{r^{R}(m)r^{R}(n)}^{r^{R}(l)}}
d_{mn; ij}^{l; (1, 2)}\Y_{r^{L}(m)r^{L}(n); i}^{r^{L}(l); (1)}\otimes 
\Y_{r^{R}(m)r^{R}(n); j}^{r^{R}(l); (2)}.
\end{equation}
where $d_{mn; ij}^{l; (1, 2)}\in \C$ for $l, m, n=1, \dots, N$,
$i=1, \dots, N_{r^{L}(m)r^{L}(n)}^{r^{L}(l)}$, 
$j=1, \dots, N_{r^{R}(m)r^{R}(n)}^{r^{R}(l)}$, and
$p, q$ any indices, and
$$\{\Y_{r^{L}(m)r^{L}(n); i}^{r^{L}(l); (1)}\;|\;i=1, \dots, 
N_{r^{L}(m)r^{L}(n)}^{r^{L}(l)}\}$$
and 
$$\{\Y_{r^{R}(m)r^{R}(n); j}^{r^{R}(l); (2)}\;|\;i=1, \dots, 
N_{r^{R}(m)r^{R}(n)}^{r^{R}(l)}\}$$
for $l, m, n=1, \dots, N$
are basis of 
$\V_{r^{L}(m)r^{L}(n)}^{r^{L}(l)}$ and 
$\V_{r^{R}(m)r^{R}(n)}^{r^{R}(l)}$, respectively. 

For $\tau\in \mathbb{H}$, let $q_{\tau}=e^{2\pi i\tau}$.
we have:

\begin{prop}\label{conv-tr-fvo}
For $u_{1}, \dots, u_{k}\in F$, the 
$q$-$\overline{q}$-trace (\ref{q-tr-fvo}) with $q=q_{\tau}$ is absolutely convergent 
when $1>|e^{2\pi iz_{1}}|>\cdots>|e^{2\pi iz_{k}} |>|q_{\tau}|>0$. Moreover, 
the sum of this 
$q$-$\overline{q}$-trace can be analytically extended to 
a multi-valued analytic function of $z_{1}, \xi_{1},
\dots, z_{k}, \xi_{k}$, $\tau$ and $\sigma$ in the region 
given by $z_{i}\ne z_{j}+m+n\tau$ for $i\ne j$ and $m, n\in \Z$
and $\xi_{i}\ne \xi_{j}+m+n\sigma$ for $i\ne j$ and $m, n\in \Z$,
with the sum 
of the $q$-$\overline{q}$-trace (\ref{q-tr-fvo}) with $q=q_{\tau}$
as a preferred value 
at the special points 
$$(z_{1}, \xi_{1}=\overline{z}_{1}, \dots, z_{k},
\xi_{k}=\overline{z}_{k}; \tau, \sigma=-\overline{\tau})$$
satisfying 
$1>|e^{2\pi iz_{1}}|>\cdots>|e^{2\pi iz_{k}} |>|q_{\tau}|>0$.
\end{prop}
\pf
This result follows immediately from (\ref{full-vo-decomp})
and the convergence and analytic extension properties for 
$q$-$\overline{q}$-traces of
intertwining operators for the vertex operator algebras $V^{L}$ 
and $V^{R}$. 
\epfv

Recall from \cite{HK} the (genus-zero) correlation functions
$$m_{k}(u_{1}, \dots, u_{k}; z_{1}, \overline{z}_{1},
\dots, z_{k}, \overline{z}_{k})$$
for $k\in \Z_{+}$ and $u_{1}, \dots, u_{k}\in F$
and their analytic extensions
$$E(m)_{k}(u_{1}, \dots, u_{k}; z_{1}, \zeta_{1},
\dots, z_{k}, \zeta_{k}).$$
The functions 
\begin{eqnarray*}
\lefteqn{m_{k}(\mathcal{U}^{L}(e^{2\pi iz_{1}})
\mathcal{U}^{R}(e^{-2\pi i\overline{z}_{1}})e^{-\pi iL^{R}(0)}u_{1}, }\nn
&&\dots, 
\mathcal{U}^{L}(e^{2\pi iz_{k}})
\mathcal{U}^{R}(e^{-2\pi i\overline{z}_{k}})e^{-\pi iL^{R}(0)} u_{k}; 
e^{2\pi iz_{1}}, e^{-2\pi i\overline{z}_{1}},
\dots, e^{2\pi iz_{k}}, e^{-2\pi i\overline{z}_{k}})
\end{eqnarray*}
are called {\it geometrically-modified (genus-zero) correlation functions}.

\begin{cor}
For  $u_{1}, \dots, u_{k}\in F$, 
\begin{eqnarray}\label{general-tr}
\lefteqn{\tr_{F}E(m)_{k}(\mathcal{U}^{L}(e^{2\pi iz_{1}})
\mathcal{U}^{R}(e^{-2\pi i\xi_{1}})e^{-\pi iL^{R}(0)}u_{1}, }\nn
&&\quad\quad \quad\quad\quad\dots, 
\mathcal{U}^{L}(e^{2\pi iz_{k}})
\mathcal{U}^{R}(e^{-2\pi i\xi_{k}})e^{-\pi iL^{R}(0)}
u_{k}; e^{2\pi iz_{1}}, e^{2\pi i\xi_{1}},\nn
&&\quad\quad \quad\quad\quad \quad\quad \quad
\dots, e^{2\pi iz_{k}}, e^{2\pi i\xi_{k}})
q_{\tau}^{L^{L}(0)-\frac{c^{L}}{24}}q_{\sigma}^{L^{R}(0)-\frac{c^{R}}{24}}
\end{eqnarray}
is absolutely convergent to a multi-valued analytic function 
of $z_{1}, \xi_{1},
\dots, z_{k}, \xi_{k}$, $\tau$ and $\sigma$ in the region 
given by $1>|e^{2\pi iz_{1}}|, \dots, |e^{2\pi iz_{k}} |>|q_{\tau}|>0$,
$1>|e^{2\pi i\xi_{1}}|, \dots, |e^{2\pi i\xi_{k}} |>|q_{\sigma}|>0$,
$z_{i}\ne z_{j}$ for $i\ne j$
and $\xi_{i}\ne \xi_{j}$ for $i\ne j$ and $m, n\in \Z$.
In particular, 
the $q_{\tau}$-$q_{\overline{\tau}}$-trace 
\begin{eqnarray}\label{tr-m-k}
\lefteqn{\tr_{F}m_{k}(\mathcal{U}^{L}(e^{2\pi iz_{1}})
\mathcal{U}^{R}(e^{-2\pi i\overline{z}_{1}})e^{-\pi iL^{R}(0)}u_{1}, }\nn
&&\quad\quad \quad\quad\quad\dots, 
\mathcal{U}^{L}(e^{2\pi iz_{k}})
\mathcal{U}^{R}(e^{-2\pi i\overline{z}_{k}})e^{-\pi iL^{R}(0)} u_{k}; 
e^{2\pi iz_{1}}, e^{-2\pi i\overline{z}_{1}},\nn
&&\quad\quad \quad\quad\quad \quad\quad \quad
\dots, e^{2\pi iz_{k}}, e^{-2\pi i\overline{z}_{k}})
q_{\tau}^{L^{L}(0)-\frac{c^{L}}{24}}
\overline{q}_{\tau}^{L^{R}(0)-\frac{c^{R}}{24}}
\end{eqnarray}
is absolutely convergent when 
$1>|e^{2\pi iz_{1}}|, \dots, |e^{2\pi iz_{k}} |>|q_{\tau}|>0$ and
$z_{i}\ne z_{j}$ for $i\ne j$.
\end{cor}
\pf
The terms in the series (\ref{general-tr}) are the analytic extensions
of the terms in (\ref{q-tr-fvo}). From Proposition \ref{conv-tr-fvo},
(\ref{q-tr-fvo}) is absolutely 
convergent and can be analytically extended to the region 
given by $z_{i}\ne z_{j}+m+n\tau$ for $i\ne j$ and $m, n\in \Z$
and $\xi_{i}\ne \xi_{j}+m+n\sigma$ for $i\ne j$ and $m, n\in \Z$.
We also know that the resulting analytic extension can be 
expanded as a series in powers of $q_{\tau}$ and $q_{\sigma}$ 
in the region given by $1>|e^{2\pi iz_{1}}|, \dots, 
|e^{2\pi iz_{k}} |>|q_{\tau}|>0$,
$1>|e^{2\pi i\xi_{1}}|, \dots, |e^{2\pi i\xi_{k}} |>|q_{\sigma}|>0$,
$z_{i}\ne z_{j}$ for $i\ne j$
and $\xi_{i}\ne \xi_{j}$ for $i\ne j$ and $m, n\in \Z$. Thus
the terms of this expansion must be equal to 
the terms in the series (\ref{general-tr}). In particular, 
(\ref{general-tr}) is absolutely convergent. 
\epfv

We shall denote the sum of the series (\ref{tr-m-k}) by
\begin{equation}\label{genus-1-corr-fn}
m^{(1)}_{k}(u_{1}, \dots, u_{k};
z_{1}, \overline{z}_{1}, \dots, z_{k}, \overline{z}_{k};
\tau, \overline{\tau}),
\end{equation}
where we use the superscript $(1)$ to indicate that this 
corresponds to a genus $1$ surface.
We have:

\begin{prop}
For any $m, n\in \Z$ and $i=1, \dots, k$,
\begin{eqnarray}\label{periods}
\lefteqn{m^{(1)}_{k}(u_{1}, \dots, u_{k};
z_{1}, \overline{z}_{1}, \dots, z_{i-1}, \overline{z}_{i-1},
z_{i}+m+n\tau, \overline{z}_{i}+m+n\overline{\tau},}\nn
&&\quad\quad \quad \quad \quad \quad \quad \quad \quad \quad \quad 
\quad \quad \quad \quad \quad \quad \quad \quad \quad   z_{i+1},
\overline{z}_{i+1}, \dots, 
z_{k}, \overline{z}_{k};
\tau, \overline{\tau})\nn
&&=m^{(1)}_{k}(u_{1}, \dots, u_{k};
z_{1}, \overline{z}_{1}, \dots, z_{k}, \overline{z}_{k};
\tau, \overline{\tau}).
\end{eqnarray}
\end{prop}
\pf
Note that 
\begin{eqnarray}\label{periods-1}
\lefteqn{m_{k}(\mathcal{U}^{L}(e^{2\pi iz_{1}})
\mathcal{U}^{R}(e^{-2\pi i\overline{z}_{1}})e^{-\pi iL^{R}(0)}u_{1}, }\nn
&&\;\dots, 
\mathcal{U}^{L}(e^{2\pi iz_{k}})
\mathcal{U}^{R}(e^{-2\pi i\overline{z}_{k}})e^{-\pi iL^{R}(0)} u_{k}; 
e^{2\pi iz_{1}}, e^{-2\pi i\overline{z}_{1}},
\dots, e^{2\pi iz_{k}}, e^{-2\pi i\overline{z}_{k}})\nn
&&
\end{eqnarray}
is a single-valued function of $e^{2\pi z_{1}}, \dots, e^{2\pi z_{k}}$. Thus 
we see that (\ref{periods}) 
holds when $n=0$. 

Now from the permutation property for conformal full field algebras
over $V^{L}\otimes V^{L}$, 
we know that 
\begin{eqnarray}\label{periods-2}
\lefteqn{m_{k}(\mathcal{U}^{L}(e^{2\pi iz_{1}})
\mathcal{U}^{R}(e^{-2\pi i\overline{z}_{1}})e^{-\pi iL^{R}(0)}u_{1}, }\nn
&&\;\dots, 
\mathcal{U}^{L}(e^{2\pi iz_{k}})
\mathcal{U}^{R}(e^{-2\pi i\overline{z}_{k}})e^{-\pi iL^{R}(0)} u_{k}; 
e^{2\pi iz_{1}}, e^{-2\pi i\overline{z}_{1}},
\dots, e^{2\pi iz_{k}}, e^{-2\pi i\overline{z}_{k}})\nn
&&=m_{k}(\mathcal{U}^{L}(e^{2\pi iz_{1}})
\mathcal{U}^{R}(e^{-2\pi i\overline{z}_{1}})e^{-\pi iL^{R}(0)}u_{1}, \nn
&&\quad\quad\quad\;\dots, \mathcal{U}^{L}(e^{2\pi iz_{i-1}})
\mathcal{U}^{R}(e^{-2\pi i\overline{z}_{i-1}})e^{-\pi iL^{R}(0)}u_{i-1},\nn
&&\quad\quad\quad\;\dots, 
\mathcal{U}^{L}(e^{2\pi iz_{i+1}})
\mathcal{U}^{R}(e^{-2\pi i\overline{z}_{i+1}})e^{-\pi iL^{R}(0)}u_{i+1},\nn
&&\quad\quad\quad\;\dots, \mathcal{U}^{L}(e^{2\pi iz_{k}})
\mathcal{U}^{R}(e^{-2\pi i\overline{z}_{k}})e^{-\pi iL^{R}(0)} u_{k},\nn
&&\quad\quad\quad\; \quad\;\;\; \mathcal{U}^{L}(e^{2\pi iz_{i}})
\mathcal{U}^{R}(e^{-2\pi i\overline{z}_{i}})e^{-\pi iL^{R}(0)}u_{i}; \nn
&&\quad\quad\quad\; \quad\;\;\;e^{2\pi iz_{1}}, e^{-2\pi i\overline{z}_{1}},
\dots, e^{2\pi iz_{i-1}}, e^{-2\pi i\overline{z}_{i-1}},\nn
&&\quad\quad\quad\;  \dots, e^{2\pi iz_{i+1}}, e^{-2\pi i\overline{z}_{i+1}},
\dots, e^{2\pi iz_{k}}, e^{-2\pi i\overline{z}_{k}}, 
e^{2\pi iz_{i}}, e^{-2\pi i\overline{z}_{i}})
\end{eqnarray}
So we can assume that $i=k$. In this case, when 
$1>|e^{2\pi iz_{1}}|>\dots>|e^{2\pi iz_{k}} |>|q_{\tau}|>0$,
we have 
\begin{eqnarray}\label{periods-3}
\lefteqn{m^{(1)}_{k}(u_{1}, \dots, u_{k};
z_{1}, \overline{z}_{1}, \dots, z_{k}, \overline{z}_{k};
\tau, \overline{\tau})}\nn
&&=\tr_{F}\mathbb{Y}(\mathcal{U}^{L}(e^{2\pi iz_{1}})
\overline{\mathcal{U}^{R}(e^{2\pi iz_{1}})}
u_{1}; e^{2\pi iz_{1}}, \overline{e^{2\pi iz_{1}}})\cdot\nn
&&\quad\quad\quad\cdots \mathbb{Y}(\mathcal{U}^{L}(e^{2\pi iz_{k}}) 
\overline{\mathcal{U}^{R}(e^{2\pi iz_{k}})}
u_{k}; e^{2\pi iz_{k}}, \overline{e^{2\pi iz_{k}}})
q_{\tau}^{L^{L}(0)-\frac{c^{L}}{24}}
\overline{q}_{\tau}^{L^{R}(0)-\frac{c^{R}}{24}}\nn
&&=\tr_{F}
\mathbb{Y}(\mathcal{U}^{L}(e^{2\pi iz_{1}})
\overline{\mathcal{U}^{R}(e^{2\pi iz_{1}})}
u_{1}; e^{2\pi iz_{1}}, \overline{e^{2\pi iz_{1}}})\cdots
q_{\tau}^{L^{L}(0)-\frac{c^{L}}{24}}
\overline{q}_{\tau}^{L^{R}(0)-\frac{c^{R}}{24}}\cdot\nn
&&\quad\quad\quad\quad\cdot
\mathbb{Y}(\mathcal{U}^{L}(e^{2\pi i(z_{k}-\tau)}) 
\overline{\mathcal{U}^{R}(e^{2\pi i(z_{k}-\tau)})}
u_{k}; e^{2\pi i(z_{k}-\tau)}, \overline{e^{2\pi i(z_{k}-\tau)}})\nn
&&=\tr_{F}\mathbb{Y}(\mathcal{U}^{L}(e^{2\pi i(z_{k}-\tau)}) 
\overline{\mathcal{U}^{R}(e^{2\pi i(z_{k}-\tau)})}
u_{k}; e^{2\pi i(z_{k}-\tau)}, \overline{e^{2\pi i(z_{k}-\tau)}})\cdot\nn
&&\quad\quad\quad\quad\cdot\mathbb{Y}(\mathcal{U}^{L}(e^{2\pi iz_{1}})
\overline{\mathcal{U}^{R}(e^{2\pi iz_{1}})}
u_{1}; e^{2\pi iz_{1}}, \overline{e^{2\pi iz_{1}}})\cdot\nn
&&\quad\quad\quad\cdots \mathbb{Y}(\mathcal{U}^{L}(e^{2\pi iz_{k-1}})
\overline{\mathcal{U}^{R}(e^{2\pi iz_{k-1}})}
u_{k-11}; e^{2\pi iz_{k-1}}, \overline{e^{2\pi iz_{k-1}}})\cdot\nn
&&\quad\quad\quad\quad\cdot
q_{\tau}^{L^{L}(0)-\frac{c^{L}}{24}}
\overline{q}_{\tau}^{L^{R}(0)-\frac{c^{R}}{24}}.\nn
\end{eqnarray}
When $1>|e^{2\pi i(z_{k}-\tau)} |>
|e^{2\pi iz_{1}}|>\dots> |e^{2\pi iz_{k-1}} |>|q_{\tau}|>0$,
the right-hand side of (\ref{periods-3}) is actually equal to 
\begin{eqnarray}\label{periods-4}
&m^{(1)}_{k}(u_{k}, u_{1}, \dots, u_{k-1};
z_{k}-\tau, \overline{z}_{k}-\overline{\tau}, 
z_{1}, \overline{z}_{1}, \dots, z_{k-1}, \overline{z}_{k-1};
\tau, \overline{\tau})&\nn
&=m^{(1)}_{k}(u_{1}, \dots, u_{k};
z_{1}, \overline{z}_{1}, \dots, 
z_{k}-\tau, \overline{z}_{k}-\overline{\tau};
\tau, \overline{\tau})&
\end{eqnarray}
Since the left-hand side of (\ref{periods-3}) and the right-hand side
of (\ref{periods-4}) are determined uniquely by their 
values in the regions 
$1>|e^{2\pi iz_{1}}|>\dots>|e^{2\pi iz_{k}} |>|q_{\tau}|>0$
and $1>|e^{2\pi i(z_{k}-\tau)} |>
|e^{2\pi iz_{1}}|>\dots> |e^{2\pi iz_{k-1}} |>|q_{\tau}|>0$, respectively,
we see that the left-hand side of (\ref{periods-3}) must be equal to 
the right-hand side
of (\ref{periods-4}), proving (\ref{periods}) in the case of $m=0$.

Combining the cases $m=0$ and $n=0$, the proposition is proved.
\epfv

From this result, we obtain immediately the following:
\begin{cor}
For $u_{1}, \dots, u_{k}\in F$, there exists a unique smooth function of 
$z_{1}, \dots, z_{k}$ and $\tau$ in the region 
$z_{i}-z_{j}\ne m+n\tau$ for $i\ne j$ and 
$\Im{(\tau)}>0$ such that in the region 
given by $1>|e^{2\pi iz_{1}}|, \dots, |e^{2\pi iz_{k}} |>|q_{\tau}|>0$ and
$z_{i}\ne z_{j}$ for $i\ne j$, this function is equal to 
(\ref{genus-1-corr-fn}).
\end{cor}

We shall still denote the smooth function given in this corollary 
by (\ref{genus-1-corr-fn}).

Now we discuss the modular invariance of these 
functions. For any single valued analytic function $f(\tau)$ of $\tau$
and any $r\in R$, 
we choose branches of the multivalued 
analytic functions $f(\tau)^{r}$ and $\left(\overline{f(\tau)}\right)^{r}$
to be $e^{r\log f(\tau)}$ and $e^{r\overline{\log f(\tau)}}$ and 
still denote them by 
$f(\tau)^{r}$ and $\left(\overline{f(\tau)}\right)^{r}$. 
In particular, for $\alpha, \beta, \gamma, \delta\in \R$,
$$\left(\frac{\alpha\tau+\beta}{\gamma\tau+\delta}\right)^{r}
=e^{r\log (\frac{\alpha\tau+\beta}{\gamma\tau+\delta})}$$
and 
\begin{eqnarray*}
\left(\frac{\alpha\overline{\tau}+\beta}
{\gamma\overline{\tau}+\delta}\right)^{r}
&=&\left(\overline{\frac{\alpha\tau+\beta}
{\gamma\tau+\delta}}\right)^{r}\nn
&=&e^{r\overline{\log (\frac{\alpha\tau+\beta}{\gamma\tau+\delta})}}.
\end{eqnarray*}

\begin{defn}
{\rm For $u_{1}, \dots, u_{k}\in V$, the function 
(\ref{genus-1-corr-fn}) is {\it invariant under the action of} 
\begin{equation}\label{sl-2-z}
\left(\begin{array}{cc}\alpha&\beta\\
\gamma&\delta\end{array}\right)\in SL(2, \Z)
\end{equation}
if 
\begin{eqnarray}\label{full-mod-inv}
\lefteqn{m^{(1)}_{k}\Biggl(
\left(\frac{1}{\gamma\tau+\delta}\right)^{L^{L}(0)}
\left(\frac{1}{\gamma\overline{\tau}+\delta}\right)^{L^{R}(0)}u_{1}, \dots, 
\left(\frac{1}{\gamma\tau+\delta}\right)^{L^{L}(0)}
\left(\frac{1}{\gamma\overline{\tau}+\delta}\right)^{L^{R}(0)}u_{k};}\nn
&&\quad\quad\quad\quad\quad\quad\quad\quad
\frac{z_{1}}{\gamma\tau+\delta}, \frac{\overline{z}_{1}}
{\gamma\overline{\tau}+\delta}, \dots, \frac{z_{k}}{\gamma\tau+\delta}, 
\frac{\overline{z}_{k}}{\gamma\overline{\tau}+\delta};
\frac{\alpha\overline{\tau}+\beta}{\gamma\tau+\delta}, 
\frac{\alpha\overline{\tau}+\beta}{\gamma\overline{\tau}+\delta}\Biggl)\nn
&&=m^{(1)}_{k}(u_{1}, \dots, u_{k};
z_{1}, \overline{z}_{1}, \dots, z_{k}, \overline{z}_{k};
\tau, \overline{\tau}).
\end{eqnarray}
A conformal full field algebra over $V^{L}\otimes V^{R}$ is said to 
be {\it invariant under the action of} (\ref{sl-2-z})
if for all $u_{1}, \dots, u_{k}\in V$,
(\ref{full-mod-inv}) holds.
If (\ref{genus-1-corr-fn}) is invariant under the action of all 
elements of $SL(2, \Z)$, we say that (\ref{genus-1-corr-fn}) 
is {\it modular invariant}.
When the function (\ref{genus-1-corr-fn}) is modular invariant,
we call it the {\it genus-one correslation function} associated to 
$u_{1}, \dots, u_{k}\in V$. 
A conformal full field algebra over $V^{L}\otimes V^{R}$ is said to 
be {\it modular invariant} if it is invariant under the action of
all elements of $SL(2, \C)$.}
\end{defn}

Since the modular group $SL(2, \Z)$ is generated by the 
elements 
$$S=\left(\begin{array}{cc}0&1\\-1&0\end{array}\right)$$
and 
$$T=\left(\begin{array}{cc}1&1\\0&1\end{array}\right),$$
to see whether a conformal full field algebra
over $V^{L}\otimes V^{R}$ is modular invariant, we 
need only discuss the invariance under these two
particular elements. 

For the element $T$, we have:

\begin{prop}
A conformal full field algebra over $V^{L}\otimes V^{R}$ 
is invariant under the action of $T$ if and only if
$c^{L}\equiv c^{R} \mod 24$.
\end{prop}
\pf
In the region 
$1>|e^{2\pi iz_{1}}|, \dots, |e^{2\pi iz_{k}} |>|q_{\tau}|>0$ and
$z_{i}\ne z_{j}$ for $i\ne j$, for $u_{1}, \dots, u_{k}\in V$, 
we have
\begin{eqnarray}\label{cent-chg-1}
\lefteqn{m^{(1)}_{k}(u_{1}, \dots, u_{k};
z_{1}, \overline{z}_{1}, \dots, z_{k}, \overline{z}_{k};
\tau+1, \overline{\tau}+1)}\nn
&&=\tr_{F}m_{k}(\mathcal{U}^{L}(e^{2\pi iz_{1}})
\mathcal{U}^{R}(e^{-2\pi i\overline{z}_{1}})e^{-\pi iL^{R}(0)}u_{1}, \nn
&&\quad\quad \quad\quad\quad\quad\quad\dots, 
\mathcal{U}^{L}(e^{2\pi iz_{k}})
\mathcal{U}^{R}(e^{-2\pi i\overline{z}_{k}})e^{-\pi iL^{R}(0)} u_{k}; \nn
&&\quad\quad \quad\quad\quad \quad\quad \quad\quad 
e^{2\pi iz_{1}}, e^{-2\pi i\overline{z}_{1}}, 
\dots, e^{2\pi iz_{k}}, e^{-2\pi i\overline{z}_{k}})
q_{\tau+1}^{L^{L}(0)-\frac{c^{L}}{24}}
\overline{q}_{\tau+1}^{L^{R}(0)-\frac{c^{R}}{24}}\nn
&&=\tr_{F} m_{k}(\mathcal{U}^{L}(e^{2\pi iz_{1}})
\mathcal{U}^{R}(e^{-2\pi i\overline{z}_{1}})e^{-\pi iL^{R}(0)}u_{1}, \nn
&&\quad\quad \quad\quad\quad \quad\quad \dots, 
\mathcal{U}^{L}(e^{2\pi iz_{k}})
\mathcal{U}^{R}(e^{-2\pi i\overline{z}_{k}})e^{-\pi iL^{R}(0)} u_{k}; \nn
&&\quad\quad \quad\quad\quad \quad\quad  \quad\quad 
e^{2\pi iz_{1}}, e^{-2\pi i\overline{z}_{1}}, 
\dots, e^{2\pi iz_{k}}, e^{-2\pi i\overline{z}_{k}})\cdot\nn
&&\quad\quad \quad\quad\quad \quad \quad
\cdot e^{2\pi i(\tau+1)(L^{L}(0)-\frac{c^{L}}{24})}
e^{-2\pi i(\overline{\tau}+1)(L^{R}(0)-\frac{c^{R}}{24})}\nn
&&=\tr_{F} m_{k}(\mathcal{U}^{L}(e^{2\pi iz_{1}})
\mathcal{U}^{R}(e^{-2\pi i\overline{z}_{1}})e^{-\pi iL^{R}(0)}u_{1}, \nn
&&\quad\quad \quad\quad\quad \quad\quad \dots, 
\mathcal{U}^{L}(e^{2\pi iz_{k}})
\mathcal{U}^{R}(e^{-2\pi i\overline{z}_{k}})e^{-\pi iL^{R}(0)} u_{k}; \nn
&&\quad\quad \quad\quad\quad \quad\quad  \quad\quad 
e^{2\pi iz_{1}}, e^{-2\pi i\overline{z}_{1}}, 
\dots, e^{2\pi iz_{k}}, e^{-2\pi i\overline{z}_{k}})\cdot\nn
&&\quad\quad \quad\quad\quad \quad\quad\cdot
e^{2\pi i\tau(L^{L}(0)-\frac{c^{L}}{24})}
e^{-2\pi i\overline{\tau}(L^{R}(0)-\frac{c^{R}}{24})}
e^{2\pi i(L^{L}(0)-L^{R}(0))}e^{2\pi i \frac{c^{L}-c^{R}}{24}}\nn
&&=\tr_{F} m_{k}(\mathcal{U}^{L}(e^{2\pi iz_{1}})
\mathcal{U}^{R}(e^{-2\pi i\overline{z}_{1}})e^{-\pi iL^{R}(0)}u_{1}, \nn
&&\quad\quad \quad\quad\quad \quad\quad \dots, 
\mathcal{U}^{L}(e^{2\pi iz_{k}})
\mathcal{U}^{R}(e^{-2\pi i\overline{z}_{k}})e^{-\pi iL^{R}(0)} u_{k}; \nn
&&\quad\quad \quad\quad\quad \quad\quad  \quad\quad 
e^{2\pi iz_{1}}, e^{-2\pi i\overline{z}_{1}}, 
\dots, e^{2\pi iz_{k}}, e^{-2\pi i\overline{z}_{k}})\cdot\nn
&&\quad\quad \quad\quad\quad \quad\quad
\cdot q_{\tau}^{L^{L}(0)-\frac{c^{L}}{24}}
\overline{q}_{\tau}^{L^{R}(0)-\frac{c^{R}}{24}}
e^{2\pi i \frac{c^{L}-c^{R}}{24}}\nn
&&=m^{(1)}_{k}(u_{1}, \dots, u_{k};
z_{1}, \overline{z}_{1}, \dots, z_{k}, \overline{z}_{k};
\tau, \overline{\tau})e^{2\pi i \frac{c^{L}-c^{R}}{24}},
\end{eqnarray}
where we have used (\ref{sing-val-1}). 
On the other hand, $F$ is invariant under $T$ means
\begin{eqnarray}\label{cent-chg-2}
\lefteqn{m^{(1)}_{k}(u_{1}, \dots, u_{k};
z_{1}, \overline{z}_{1}, \dots, z_{k}, \overline{z}_{k};
\tau+1, \overline{\tau}+1)}\nn
&&=m^{(1)}_{k}(u_{1}, \dots, u_{k};
z_{1}, \overline{z}_{1}, \dots, z_{k}, \overline{z}_{k};
\tau, \overline{\tau})
\end{eqnarray}
{}From (\ref{cent-chg-1}) and (\ref{cent-chg-2}) and the fact that 
(\ref{genus-1-corr-fn}) are not all $0$ (for example, 
it is not $0$ when $u_{1}=\cdots=u_{k}=\one$,
we see that $F$ is invariant under $T$ if and 
only if $e^{2\pi i \frac{c^{L}-c^{R}}{24}}=0$, or equivalently, 
$c^{L}\equiv c^{R} \mod 24$.
\epfv

The invariance under $S$ is the most important. We need to first 
introduce matrix elements associated to the actions of $S$ on the 
chiral genus-one correlation functions. For 
$a^{L}, a^{L}_{1}\in \A^{L}$, let 
$\{\Y_{a^{L}a^{L}_{1}; i}^{a^{L}_{1}; (1)}\; |\; i=1, \dots, 
N_{a^{L}a^{L}_{1}}^{a^{L}_{1}}\}$ and
$\{\Y_{a^{L}a^{L}_{1}; i}^{a^{L}_{1}; (2)}\; |\; i=1, \dots, 
N_{a^{L}a^{L}_{1}}^{a^{L}_{1}}\}$ be basis of the spaces 
$\V_{a^{L}a^{L}_{1}}^{a^{L}_{1}}$ of intertwining operators of types 
${W^{a^{L}_{1}}\choose W^{a^{L}}W^{a^{L}_{1}}}$ 
and for $a^{R}, a^{R}_{1}\in \A^{R}$, let $\{Y_{a^{R}a^{R}_{1}; i}
^{a^{R}_{1}; (1)}\; 
|\; i=1, \dots, 
N_{a^{R}a^{R}_{1}}^{a^{R}_{1}}\}$, $\{Y_{a^{R}a^{R}_{1}; i}
^{a^{R}_{1}; (2)}\; 
|\; i=1, \dots, 
N_{a^{R}a^{R}_{1}}^{a^{R}_{1}}\}$ be  basis of the space
$\V_{a^{R}a^{R}_{1}}^{a^{R}_{1}}$ of intertwining operators of types 
${W^{a^{R}_{1}}\choose W^{a^{R}}W^{a^{R}_{1}}}$. 
From Theorem \ref{mod-inv}, we know that for $a^{L}\in \A^{L}$,
there exist $S(\Y_{a^{L}a^{L}_{1}; i}^{a^{L}_{1}; (1)};
\Y_{a^{L}a^{L}_{2}; j}^{a^{L}_{2}; (2)})$ for $a^{L}_{1}, 
a^{L}_{2}\in \A^{L}$, $i=1, \dots, 
N_{a^{L}a^{L}_{1}}^{a^{L}_{1}}$, $j=1, \dots, 
N_{a^{L}a^{L}_{2}}^{a^{L}_{2}}$  and for
$a^{R}\in \A^{R}$, there exist
$S(\Y_{a^{R}a^{R}_{1}; i}^{a^{R}_{1}; (1)};
\Y_{a^{R}a^{R}_{2}; j}^{a^{R}_{2}; (2)})$ for $a^{R}_{1}, 
a^{R}_{2}\in \A^{R}$, $i=1, \dots, 
N_{a^{R}a^{R}_{1}}^{a^{R}_{1}}$, $j=1, \dots, 
N_{a^{R}a^{R}_{2}}^{a^{R}_{2}}$, such that 
\begin{eqnarray}\label{s-left}
\lefteqn{(\Phi(\Y_{a^{L}a^{L}_{1}; i}^{a^{L}_{1}; (1)}))
\left(\left(-\frac{1}{\tau}\right)^{L^{L}(0)}w_{a^{L}};
-\frac{z}{\tau}; -\frac{1}{\tau}\right)}\nn
&&=\sum_{a^{L}_{2}\in \A^{L}}\sum_{j=1}^{N_{a^{L}a^{L}_{2}}^{a^{L}_{2}}}
S(\Y_{a^{L}a^{L}_{1}; i}^{a^{L}_{1}; (1)};
\Y_{a^{L}a^{L}_{2}; j}^{a^{L}_{2}; (2)})
(\Phi(\Y_{a^{L}a^{L}_{2}; j}^{a^{L}_{2}; (2)}))(w_{a^{L}}, 
z; \tau)
\end{eqnarray}
for $w_{a^{L}}\in W^{a^{L}}$
and 
\begin{eqnarray}\label{s-right}
\lefteqn{(\Phi(\Y_{a^{R}a^{R}_{1}; i}^{a^{R}_{1}; (1)}))
\left(\left(-\frac{1}{\tau}\right)^{L^{R}(0)}w_{a^{R}};
-\frac{z}{\tau}; -\frac{1}{\tau}\right)}\nn
&&=\sum_{a^{R}_{2}\in \A^{R}}\sum_{j=1}^{N_{a^{R}a^{R}_{2}}^{a^{R}_{2}}}
S(\Y_{a^{R}a^{R}_{1}; i}^{a^{R}_{1}; (1)};
\Y_{a^{R}a^{R}_{2}; j}^{a^{R}_{2}; (2)})
(\Phi(\Y_{a^{R}a^{R}_{2}; j}^{a^{R}_{2}; (2)}))(w_{a^{R}};
z; \tau)
\end{eqnarray}

We need:

\begin{lemma}
For $w_{a^{R}}\in W^{a^{R}}$,
\begin{eqnarray}\label{s-inverse}
\lefteqn{(\Phi(\Y_{a^{R}a^{R}_{1}; i}^{a^{R}_{1}; (1)}))
\left(e^{-\pi iL^{R}(0)}\left(-\frac{1}{\overline{\tau}}\right)^{L^{R}(0)}
w_{a^{R}};
\frac{\overline{z}}{\overline{\tau}}; \frac{1}{\overline{\tau}}\right)}\nn
&&=\sum_{a^{R}_{2}\in \A^{R}}\sum_{j=1}^{N_{a^{R}a^{R}_{2}}^{a^{R}_{2}}}
S^{-1}(\Y_{a^{R}a^{R}_{1}; i}^{a^{R}_{1}; (1)};
\Y_{a^{R}a^{R}_{2}; j}^{a^{R}_{2}; (2)})
(\Phi(\Y_{a^{R}a^{R}_{2}; j}^{a^{R}_{2}; (2)}))(e^{-\pi iL^{R}(0)}w_{a^{R}};
-\overline{z}; -\overline{\tau}),\nn
\end{eqnarray}
where $S^{-1}(\Y_{a^{R}a^{R}_{1}; i}^{a^{R}_{1}; (1)};
\Y_{a^{R}a^{R}_{2}; j}^{a^{R}_{2}; (2)})$ are the 
matrix elements of the inverse of the action 
of $S$, that is,
$$\sum_{a^{R}_{2}\in \A^{R}}
\sum_{j=1}^{N_{a^{R}a^{R}_{2}}^{a^{R}_{2}}}
S^{-1}(\Y_{a^{R}a^{R}_{1}; i}^{a^{R}_{1}; (1)};
\Y_{a^{R}a^{R}_{2}; j}^{a^{R}_{2}; (2)})
S(\Y_{a^{R}a^{R}_{2}; j}^{a^{R}_{2}; (2)};
\Y_{a^{R}a^{R}_{3}; k}^{a^{R}_{3}; (1)})
=\delta_{a^{R}_{1}a^{R}_{3}}\delta_{ik}.$$
\end{lemma}
\pf
From (\ref{s-right}), we see that the action of $S$ maps 
$(\Phi(\Y_{a^{R}a^{R}_{1}; j}^{a^{R}_{1}; (1)}))(w_{a^{R}};
z; \tau)$ to 
$$(\Phi(\Y_{a^{R}a^{R}_{1}; i}^{a^{R}_{1}; (1)}))
\left(\left(-\frac{1}{\tau}\right)^{L^{R}(0)}w_{a^{R}};
-\frac{z}{\tau}; -\frac{1}{\tau}\right).$$
So the inverse of the action of $S$ maps 
$\Phi(\Y_{a^{R}a^{R}_{1}; j}^{a^{R}_{1}; (1)})(w_{a^{R}};
z; \tau)$
to 
$$(\Phi(\Y_{a^{R}a^{R}_{1}; i}^{a^{R}_{1}; (1)}))
\left(e^{-\log \tau L^{R}(0)}w_{a^{R}};
\frac{z}{\tau}; -\frac{1}{\tau}\right).$$
Thus we have 
\begin{eqnarray}\label{s-inverse-1}
\lefteqn{(\Phi(\Y_{a^{R}a^{R}_{1}; i}^{a^{R}_{1}; (1)}))
\left(e^{-\log \tau L^{R}(0)}w_{a^{R}};
\frac{z}{\tau}; -\frac{1}{\tau}\right)}\nn
&&=\sum_{a^{R}_{2}\in \A^{R}}\sum_{j=1}^{N_{a^{R}a^{R}_{2}}^{a^{R}_{2}}}
S^{-1}(\Y_{a^{R}a^{R}_{1}; i}^{a^{R}_{1}; (1)};
\Y_{a^{R}a^{R}_{2}; j}^{a^{R}_{2}; (2)})
(\Phi(\Y_{a^{R}a^{R}_{2}; j}^{a^{R}_{2}; (2)}))(w_{a^{R}};
z; \tau),
\end{eqnarray}
Note that when $\tau$ is in the upper  half plane, 
so is $-\overline{\tau}$. So we can substitute 
$-\overline{\tau}$ for $\tau$ in (\ref{s-inverse-1}) above. 
We also substitute  $-\overline{z}$ 
for $z$ and $e^{-\pi iL^{R}(0)}w_{a^{R}}$ for 
$w_{a^{R}}$ in (\ref{s-inverse-1}).
We obtain
\begin{eqnarray}\label{s-inverse-2}
\lefteqn{(\Phi(\Y_{a^{R}a^{R}_{1}; i}^{a^{R}_{1}; (1)}))
\left(e^{-\log (-\overline{\tau}) L^{R}(0)}
e^{-\pi iL^{R}(0)}w_{a^{R}};
\frac{\overline{z}}{\overline{\tau}}; \frac{1}{\overline{\tau}}\right)}\nn
&&=\sum_{a^{R}_{2}\in \A^{R}}\sum_{j=1}^{N_{a^{R}a^{R}_{2}}^{a^{R}_{2}}}
S^{-1}(\Y_{a^{R}a^{R}_{1}; i}^{a^{R}_{1}; (1)};
\Y_{a^{R}a^{R}_{2}; j}^{a^{R}_{2}; (2)})
(\Phi(\Y_{a^{R}a^{R}_{2}; j}^{a^{R}_{2}; (2)}))(e^{-\pi iL^{R}(0)}w_{a^{R}};
-\overline{z}; -\overline{\tau}).\nn
\end{eqnarray}
Note that $e^{-\pi iL^{R}(0)}$ commutes with 
$e^{-\log (-\overline{\tau}) L^{R}(0)}$. 
Since $0< \arg \tau <\pi$, we have $\pi <\arg \overline{\tau}< 2\pi$, 
$0< \arg -\overline{\tau}<\pi$ and $0<\arg -\frac{1}{\tau}
<\pi$. So $\arg -\frac{1}{\tau}=\arg -\overline{\tau}$.
Then by our convention, on $W^{a^{R}}$, we have
\begin{eqnarray}\label{s-inverse-3}
e^{-\log (-\overline{\tau})L^{R}(0)}
&=&e^{(-\log |-\overline{\tau}|
-i\arg (-\overline{\tau}))L^{R}(0)}\nn
&=&e^{(\log |-\frac{1}{\tau}|
-i\arg (-\frac{1}{\tau}))L^{R}(0)}\nn
&=&e^{\overline{\log (-\frac{1}{\tau})}L^{R}(0)}\nn
&=&\left(-\frac{1}{\overline{\tau}}\right)^{L^{R}(0)}.
\end{eqnarray}
Changing the order
of $e^{-\log (-\overline{\tau}) L^{R}(0)}$ and 
$e^{-\pi iL^{R}(0)}$ and then using (\ref{s-inverse-3}), 
we see that (\ref{s-inverse-2})
gives (\ref{s-inverse}).
\epfv

We now have:

\begin{thm}\label{s-mod-inv-thm}
A conformal full field algebra over $V^{L}\otimes V^{R}$ 
is invariant under the action of $S$ if and only if
\begin{eqnarray}\label{s-mod-inv}
\lefteqn{\sum_{n=1}^{N}\sum_{i=1}^{N_{r^{L}(m)r^{L}(n)}^{r^{L}(n)}}
\sum_{j=1}^{N_{r^{R}(m)r^{R}(n)}^{r^{R}(n)}}
d_{mn; i j}^{n; (1, 1)}S(\Y_{r^{L}(m)r^{L}(n); i}^{r^{L}(n); (1)};
\Y_{r^{L}(m)a^{L}; k}^{a^{L}; (2)})\cdot}\nn
&&\quad\quad\quad\quad\quad\quad\quad\quad\quad\quad\quad
\cdot S^{-1}(\Y_{r^{R}(m)r^{R}(n); j}^{r^{R}(n); (1)};
\Y_{r^{R}(m)a^{R}; l}^{a^{R}; (2)})\nn
&&=\sum_{p\in (r^{L})^{-1}(a^{L})\cap (r^{R})^{-1}(a^{R})}
d_{mp; k l}^{p; (2, 2)}
\end{eqnarray}
for $m=1, \dots, N$, $a^{L}\in \A^{L}$, $a^{R}\in \A^{R}$, 
$k=1, \dots, N_{r^{L}(m)a^{L}}^{a^{L}}$, $l=1, \dots, 
N_{r^{R}(m)a^{R}}^{a^{R}}$.
\end{thm}
\pf
From the convergence property  of the correlation function maps 
(see Definition \ref{ffa}), 
we see that all the functions of the form
(\ref{genus-1-corr-fn}) are  invariant under the action of $S$ 
if and only if
all the functions of the form (\ref{genus-1-corr-fn}) with 
$k=1$ are invariant under the action of $S$. We now show that 
all the functions of the form (\ref{genus-1-corr-fn}) with 
$k=1$ are invariant under the action of $S$
if and only if (\ref{s-mod-inv})
holds. 

If (\ref{s-mod-inv}) holds, then from 
(\ref{full-vo-decomp}),  the definition of $q$-$\overline{q}$-trace of 
a map from $F$ to $\overline{F}$, (\ref{s-left}) and (\ref{s-inverse}), we have
\begin{eqnarray}\label{s-mod-inv-1}
\lefteqn{m^{(1)}_{1}\left(\left(-\frac{1}{\tau}\right)^{L^{L}(0)}w_{r^{L}(m)}
\otimes \left(-\frac{1}{\overline{\tau}}\right)^{L^{R}(0)}w_{r^{R}(m)}; 
-\frac{z}{\tau}, -\frac{\overline{z}}{\overline{\tau}}; 
-\frac{1}{\tau}, -\frac{1}{\overline{\tau}}\right)}\nn
&&=\tr_{F}\mathbb{Y}\Biggl(\mathcal{U}^{L}(e^{-2\pi i\frac{z}{\tau}})
\mathcal{U}^{R}(e^{2\pi i\frac{\overline{z}}{\overline{\tau}}})
\left(-\frac{1}{\tau}\right)^{L^{L}(0)}
w_{r^{L}(m)}\nn
&&\quad\quad\quad\quad\quad\quad
\otimes \left(-\frac{1}{\overline{\tau}}\right)^{L^{R}(0)}
w_{r^{R}(m)}; e^{-2\pi i\frac{z}{\tau}}, 
e^{2\pi i\frac{\overline{z}}{\overline{\tau}}}\Biggr)
q_{-\frac{1}{\tau}}^{L^{L}(0)-\frac{c^{L}}{24}}
q_{\frac{1}{\overline{\tau}}}^{L^{R}(0)-\frac{c^{R}}{24}}\nn
&&=\tr_{F}\sum_{l, n=1}^{N}
\sum_{i=1}^{N_{r^{L}(m)r^{L}(n)}^{r^{L}(l)}}
\sum_{j=1}^{N_{r^{R}(m)r^{R}(n)}^{r^{R}(l)}}
d_{mn; i j}^{l; (1, 1)}\cdot\nn
&&\quad\quad\cdot 
\Biggl(\Y_{r^{L}(m)r^{L}(n); i}^{r^{L}(l); (1)}
\Biggl(\mathcal{U}^{L}(e^{-2\pi i\frac{z}{\tau}})
\left(-\frac{1}{\tau}\right)^{L^{L}(0)}w_{r^{L}(m)}, 
e^{-2\pi i\frac{z}{\tau}}\Biggr)\nn
&&\quad\quad\quad  \otimes 
\Y_{r^{R}(m)r^{R}(n); j}^{r^{R}(l); (1)}
\Biggl(\mathcal{U}^{R}(e^{2\pi i\frac{\overline{z}}{\overline{\tau}}})
e^{-\pi iL^{R}(0)}\left(-\frac{1}{\overline{\tau}}\right)^{L^{R}(0)}w_{r^{R}(m)}, 
e^{2\pi i\frac{\overline{z}}{\overline{\tau}}}\Biggr)\Biggr)\cdot\nn
&&\quad\quad\quad\quad\quad\quad\quad\quad
\cdot q_{-\frac{1}{\tau}}^{L^{L}(0)-\frac{c^{L}}{24}}
q_{\frac{1}{\overline{\tau}}}^{L^{R}(0)-\frac{c^{R}}{24}}\nn
&&=\sum_{n=1}^{N}
\sum_{i=1}^{N_{r^{L}(m)r^{L}(n)}^{r^{L}(l)}}
\sum_{j=1}^{N_{r^{R}(m)r^{R}(n)}^{r^{R}(l)}}
d_{mn; i j}^{n; (1, 1)}
\Biggl(\Biggl(\tr_{W^{r^{L}(n)}}\nn
&&\quad\quad
\Y_{r^{L}(m)r^{L}(n); i}^{r^{L}(n); (1)}\Biggl(
\mathcal{U}^{L}(e^{-2\pi i\frac{z}{\tau}})
\left(-\frac{1}{\tau}\right)^{L^{L}(0)}
w_{r^{L}(m)}, e^{-2\pi i\frac{z}{\tau}}\Biggr)
 q_{-\frac{1}{\tau}}^{L^{L}(0)-\frac{c^{L}}{24}}\Biggr)  
\nn
&&\quad\quad
\cdot\otimes 
\Biggl(\tr_{W^{r^{R}(n)}}\nn
&&\quad\quad\quad \Y_{r^{R}(m)r^{R}(n); j}^{r^{R}(n); (1)}
\Biggl(\mathcal{U}^{R}(e^{2\pi i\frac{\overline{z}}{\overline{\tau}}})
e^{-\pi iL^{R}(0)}\left(-\frac{1}{\overline{\tau}}\right)^{L^{R}(0)}
w_{r^{R}(m)}, 
e^{2\pi i\frac{\overline{z}}{\overline{\tau}}}\Biggr)\cdot\nn
&&\quad \quad\quad\quad\quad\quad\quad\quad\quad\quad\quad\quad
\cdot q_{\frac{1}{\overline{\tau}}}^{L^{R}(0)-\frac{c^{R}}{24}}
\Biggr)\Biggr)\nn
&&=\sum_{n=1}^{N}
\sum_{i=1}^{N_{r^{L}(m)r^{L}(n)}^{r^{L}(l)}}
\sum_{j=1}^{N_{r^{R}(m)r^{R}(n)}^{r^{R}(l)}}
d_{mn; i j}^{n; (1, 1)}\nn
&&\quad\quad
\Biggl((\Phi(\Y_{r^{L}(m)r^{L}(n); i}^{r^{L}(n); (1)}))\Biggl(
\left(-\frac{1}{\tau}\right)^{L^{L}(0)}
w_{r^{L}(m)}; -\frac{z}{\tau}; -\frac{1}{\tau}\Biggr)
\nn
&&\quad\quad\quad
\otimes 
(\Phi(\Y_{r^{R}(m)r^{R}(n); j}^{r^{R}(n); (1)}))
\Biggl(e^{-\pi iL^{R}(0)}\left(-\frac{1}{\overline{\tau}}\right)^{L^{R}(0)}
w_{r^{R}(m)}; \frac{\overline{z}}{\overline{\tau}}; \frac{1}{\overline{\tau}}\Biggr)
\Biggr)\nn
&&=\sum_{k=1}^{N_{r^{L}(m)a^{L}}^{a^{L}}}
\sum_{l=1}^{N_{r^{R}(m)a^{R}}^{a^{R}}}
\sum_{a^{L}\in \A^{L}}\sum_{a^{R}\in \A^{R}}\nn
&&\quad\Biggl(\sum_{n=1}^{N}
\sum_{i=1}^{N_{r^{L}(m)r^{L}(n)}^{r^{L}(l)}}
\sum_{j=1}^{N_{r^{R}(m)r^{R}(n)}^{r^{R}(l)}}
\nn
&&\quad\quad
d_{mn; i j}^{n; (1, 1)}S(\Y_{r^{L}(m)r^{L}(n); i}^{r^{L}(n); (1)};
\Y_{r^{L}(m)a^{L}; k}^{a^{L}; (2)}) 
S^{-1}(\Y_{r^{R}(m)r^{R}(n); j}^{r^{R}(n); (1)};
\Y_{r^{R}(m)a^{R}; l}^{a^{R}; (2)})\Biggr)\cdot\nn
&& \quad\quad \cdot (((\Phi(
\Y_{r^{L}(m)a^{L}; k}^{a^{L}; (2)}))(
w_{r^{L}(m)}; z; \tau))  
\nn
&&\quad\quad\quad \quad
\otimes 
((\Phi(\Y_{r^{R}(m)a^{R}; l}^{a^{R}; (2)}))
(e^{-\pi iL^{R}(0)}w_{r^{R}(m)}; -\overline{z}; 
-\overline{\tau})))\nn
&&=\sum_{k=1}^{N_{r^{L}(m)a^{L}}^{a^{L}}}
\sum_{l=1}^{N_{r^{R}(m)a^{R}}^{a^{R}}}
\sum_{a^{L}\in \A^{L}}\sum_{a^{R}\in \A^{R}}\nn
&&\quad\Biggl(\sum_{n=1}^{N}
\sum_{i=1}^{N_{r^{L}(m)r^{L}(n)}^{r^{L}(l)}}
\sum_{j=1}^{N_{r^{R}(m)r^{R}(n)}^{r^{R}(l)}}
\nn
&&\quad\quad
d_{mn; i j}^{n; (1, 1)}S(\Y_{r^{L}(m)r^{L}(n); i}^{r^{L}(n); (1)};
\Y_{r^{L}(m)a^{L}; k}^{a^{L}; (2)}) 
S^{-1}(\Y_{r^{R}(m)r^{R}(n); j}^{r^{R}(n); (1)};
\Y_{r^{R}(m)a^{R}; l}^{a^{R}; (2)})\Biggr)\cdot\nn
&& \quad\quad \cdot \Biggl(\Biggl(\tr_{W^{a^{L}}}
\Y_{r^{L}(m)a^{L}; k}^{a^{L}; (2)}\Biggl(
\mathcal{U}^{L}(e^{2\pi iz})
w_{r^{L}(m)}, e^{2\pi iz}\Biggr)
q_{\tau}^{L^{L}(0)-\frac{c^{L}}{24}}\Biggr)  
\nn
&&\quad\quad\quad \quad
\otimes 
\Biggl(\tr_{W^{a^{R}}} \Y_{r^{R}(m)a^{R}; l}^{a^{R}; (2)}
\Biggl(\mathcal{U}^{R}(e^{-2\pi i\overline{z}})
e^{-\pi iL^{R}(0)}w_{r^{R}(m)}, 
e^{-2\pi i\overline{z}}\Biggr)
q_{-\overline{\tau}}^{L^{R}(0)-\frac{c^{R}}{24}}
\Biggr)\Biggr)\nn
&&=\sum_{k=1}^{N_{r^{L}(m)a^{L}}^{a^{L}}}
\sum_{l=1}^{N_{r^{R}(m)a^{R}}^{a^{R}}}
\sum_{a^{L}\in \A^{L}}\sum_{a^{R}\in \A^{R}}
\sum_{p\in (r^{L})^{-1}(a^{L})\cap (r^{R})^{-1}(a^{R})}
d_{mp; k l}^{p; (2, 2)}\cdot\nn
&& \quad\quad \cdot \Biggl(\Biggl(\tr_{W^{a^{L}}}
\Y_{r^{L}(m)a^{L}; k}^{a^{L}; (2)}\Biggl(
\mathcal{U}^{L}(e^{2\pi iz})
w_{r^{L}(m)}, e^{2\pi iz}\Biggr)
q_{\tau}^{L^{L}(0)-\frac{c^{L}}{24}}\Biggr)  
\nn
&&\quad\quad\quad \quad
\otimes 
\Biggl(\tr_{W^{a^{R}}} \Y_{r^{R}(m)a^{R}; l}^{a^{R}; (2)}
\Biggl(\mathcal{U}^{R}(e^{-2\pi i\overline{z}})
e^{-\pi iL^{R}(0)}w_{r^{R}(m)}, 
e^{-2\pi i\overline{z}}\Biggr)
q_{-\overline{\tau}}^{L^{R}(0)-\frac{c^{R}}{24}}
\Biggr)\Biggr)\nn
&&=\sum_{k=1}^{N_{r^{L}(m)a^{L}}^{a^{L}}}
\sum_{l=1}^{N_{r^{R}(m)a^{R}}^{a^{R}}}
\sum_{p=1}^{N}
d_{mp; k l}^{p; (2, 2)}\cdot\nn
&& \quad\quad \cdot \Biggl(\Biggl(\tr_{W^{a^{L}}}
\Y_{r^{L}(m)r^{L}(p); k}^{r^{L}(p); (2)}\Biggl(
\mathcal{U}^{L}(e^{2\pi iz})
w_{r^{L}(m)}, e^{2\pi iz}\Biggr)
q_{\tau}^{L^{L}(0)-\frac{c^{L}}{24}}\Biggr)  
\nn
&&\quad\quad\quad \quad
\otimes 
\Biggl(\tr_{W^{a^{R}}} \Y_{r^{R}(m)r^{R}(p); l}^{r^{R}(p); (2)}
\Biggl(\mathcal{U}^{R}(e^{-2\pi i\overline{z}})
e^{-\pi iL^{R}(0)}w_{r^{R}(m)}, 
e^{-2\pi i\overline{z}}\Biggr)
q_{-\overline{\tau}}^{L^{R}(0)-\frac{c^{R}}{24}}
\Biggr)\Biggr)\nn
&&=m^{(1)}_{1}(w_{r^{L}(m)}\otimes w_{r^{R}(m)}; 
z, \overline{z}; \tau, \overline{\tau}).
\end{eqnarray}

Conversely, the calculation (\ref{s-mod-inv-1}) also 
shows that if $F$ is modular invariant, 
then (\ref{s-mod-inv}) holds.
\epfv

\renewcommand{\theequation}{\thesection.\arabic{equation}}
\renewcommand{\thethm}{\thesection.\arabic{thm}}
\setcounter{equation}{0}
\setcounter{thm}{0}

\section{$S$-matrices associated to irreducible modules}

We assume that $V$ is a vertex operator algebra satisfying 
Conditions \ref{u-o-v}, \ref{c-red} and \ref{c-2-c-f}.
We shall continue to use the notations in the preceding sections, 
in particular, those in Section 2.

For $a_{1}, a_{2}\in \A$ and $\Y$ an intertwining operator
of type ${W^{a_{1}}\choose W^{a_{2}}W^{a_{1}}}$, we consider
\begin{eqnarray*}
\Psi(\Y): \coprod_{a\in \A}W^{a}&\to &
\mathbb{G}_{1; 1}\nn
w_{a}&\mapsto & 
\Psi(\Y)(w_{a}; \tau)
\end{eqnarray*}
for $a\in \A$ and $w_{a}\in W^{a}$, 
where for $w_{a}\in W^{a}$,
$$\Psi(\Y)(w_{a}; \tau)=0$$
when $a\ne a_{2}$
and
\begin{eqnarray*}
\Psi(\Y)(w_{a_{2}}; \tau)&=&\Phi(\Y)(w_{a_{2}}; z; \tau)\nn
&=&E(\tr_{W^{a_{1}}}
\Y(\mathcal{U}(e^{2\pi iz})w_{a_{2}}, 
e^{2\pi i z})
q_{\tau}^{L(0)-\frac{c}{24}})
\end{eqnarray*}
when $a= a_{2}$, 
and $\mathbb{G}_{1; 1}$ is  the space spanned by
functions 
of $\tau$ of the form 
$\Psi(\Y)(w_{a}; \tau)$ 
for $a, a_{1}, a_{2}\in \A$ and $\Y$ an intertwining operator
of type ${W^{a_{1}}\choose W^{a_{2}}W^{a_{1}}}$.
Here we use the notation $\Psi(\Y)(w_{a}; \tau)$ instead of 
$\Psi(\Y)(w_{a}; z; \tau)$ because 
it was shown in the proof of Theorem 7.3 in \cite{H7}
that $\Psi(\Y)(w_{a_{2}}; \tau)$ is indeed independent of $z$.
Let $\F_{1; 1}$ be the space of all maps of the form 
$\Psi(\Y)$ for $a_{1}, a_{2}\in \A$ and $\Y$ an intertwining operator
of type ${W^{a_{1}}\choose W^{a_{2}}W^{a_{1}}}$.

Let $\mathbb{G}_{1; 2}$ be the space
of all single-valued analytic functions 
on  the universal covering 
$\widetilde{M}_{1}^{2}$ of 
$$M_{1}^{2}=\{(z_{1}, z_{2}, \tau)\in \C^{3}\;|\; 
z_{1}\ne z_{2}+p\tau +q \; \mbox{\rm for}\; 
p, q\in \Z,\; \tau\in \mathbb{H}\}$$
spanned by functions of the form
$$%\begin{equation}\label{1-product}
(\Phi(\Y_{1}\otimes \Y_{2}))(w_{a_{1}}, w_{a_{2}};
z_{1}, z_{2}; \tau)
$$%\end{equation} 
for $a_{1}, a_{2}, a_{3}, a_{4}\in \A$, $w_{a_{1}}\in W^{a_{1}}$,
$w_{a_{2}}\in W^{a_{2}}$, and $\Y_{1}$ and $\Y_{2}$ intertwining operators 
of types ${W^{a_{4}}\choose W^{a_{1}}W^{a_{3}}}$ and 
${W^{a_{3}}\choose W^{a_{2}}W^{a_{4}}}$, respectively. 
This space is also  spanned by the analytic extensions
\begin{equation}\label{1-iterate}
E(\tr_{W^{a_{6}}}
\Y_{3}(\mathcal{U}(e^{2\pi iz_{1}})
\Y_{4}
(w_{a_{1}},  z_{1}-z_{2})w_{a_{2}}, 
e^{2\pi i z_{2}})
q_{\tau}^{L(0)-\frac{c}{24}})
\end{equation}
to the region given by $\Im{\tau}>0$, $z_{1}\ne z_{2}+k\tau+l$
for $k, l\in \Z$, of functions of the form
$$\tr_{W^{a_{6}}}
\Y_{3}(\mathcal{U}(e^{2\pi iz_{1}})
\Y_{4}
(w_{a_{1}},  z_{1}-z_{2})w_{a_{2}}, 
e^{2\pi i z_{2}})
q_{\tau}^{L(0)-\frac{c}{24}}$$
for $a_{1}, a_{2}, a_{3}, a_{4}\in \A$, $w_{a_{1}}\in W^{a_{1}}$,
$w_{a_{2}}\in W^{a_{2}}$, and $\Y_{3}$ and $\Y_{4}$ intertwining operators 
of types ${W^{a_{4}}\choose W^{a}W^{a_{4}}}$ and 
${W^{a}\choose W^{a_{1}}W^{a_{2}}}$, respectively. 
For $a_{3}\in \A$, let 
$\mathbb{G}^{a_{3}}_{1; 2}$ be the space
of all single-valued analytic functions 
on  the universal covering 
$\widetilde{M}_{1}^{2}$ of 
$M_{1}^{2}$
spanned by functions of the form
(\ref{1-iterate}) for $a_{1}, a_{2}, a_{4}\in \A$, 
$w_{a_{1}}\in W^{a_{1}}$,
$w_{a_{2}}\in W^{a_{2}}$, and $\Y_{3}$ and $\Y_{4}$ intertwining operators 
of types ${W^{a_{4}}\choose W^{a_{3}}W^{a_{4}}}$ and 
${W^{a_{3}}\choose W^{a_{1}}W^{a_{2}}}$, respectively.

Now for $a_{1}, a_{2}, a_{3}\in \A$, 
let $\{\Y_{a_{1}a_{2}; i}^{a_{3}; (1)}\;|\; i=1, \dots, 
N_{a_{1}a_{2}}^{a_{2}}\}$ and 
$\{\Y_{a_{1}a_{2}; i}^{a_{3}; (2)}\;|\; i=1, \dots, 
N_{a_{1}a_{2}}^{a_{2}}\}$ be basis of 
the space $\V_{a_{1}a_{2}}^{a_{3}}$ of intertwining operators of 
type ${W^{a_{3}}\choose W^{a_{1}}W^{a_{2}}}$.
We need the following result proved in \cite{H6}: 

\begin{prop}\label{1-independence0}
For $a_{1}, a_{2}\in \A$, the maps from $W^{a_{1}}\otimes W^{a_{2}}$
to $\mathbb{G}_{1; 2}$ given by 
\begin{eqnarray*}
\lefteqn{w_{a_{1}}\otimes w_{a_{2}}\mapsto}\nn 
&& E(\tr_{W^{a_{4}}}
\Y_{a_{1}a_{3}; i}^{a_{4}; (1)}(\mathcal{U}(e^{2\pi iz_{1}})w_{a_{1}}, 
e^{2\pi i z_{1}})\Y_{a_{2}a_{4}; j}^{a_{3}; (2)}
(\mathcal{U}(e^{2\pi iz_{2}})w_{a_{2}},  e^{2\pi i z_{2}})
q_{\tau}^{L(0)-\frac{c}{24}}),
\end{eqnarray*}
$a_{3}, a_{4}\in \A$, $i=1, \dots, N_{a_{1}a_{3}}^{a_{4}}$,
$j=1, \dots, N_{a_{2}a_{4}}^{a_{3}}$,
are linearly independent. Similarly, for $a_{1}, a_{2}\in \A$, 
the maps from $W^{a_{1}}\otimes W^{a_{2}}$
to $\mathbb{G}_{1; 2}$ given by 
$$w_{a_{1}}\otimes w_{a_{2}}\mapsto E(\tr_{W^{a_{4}}}
\Y_{a_{3}a_{4}; k}^{a_{4}; (3)}(\mathcal{U}(e^{2\pi iz_{2}})
\Y_{a_{1}a_{2}; l}^{a_{3}; (4)}
(w_{a_{1}},  z_{1}-z_{2})w_{a_{2}}, 
e^{2\pi i z_{2}})
q_{\tau}^{L(0)-\frac{c}{24}}),$$
$a_{3}, a_{4}\in \A$, $k=1, \dots, N_{a_{3}a_{4}}^{a_{3}}$,
$l=1, \dots, N_{a_{1}a_{2}}^{a_{3}}$,
are linearly independent.
\end{prop}

For $a\in \A$ and intertwining 
operators $\Y_{1}$ and $\Y_{2}$ of types for intertwining 
operators 
of types ${W^{a_{1}}\choose W^{a}W^{a_{1}}}$ and 
${W^{a}\choose W^{a_{2}}W^{a'_{2}}}$, respectively, we consider the maps
\begin{eqnarray*}
\Psi(\Y_{1}\otimes \Y_{2}):
\coprod_{a_{3}\in \A}W^{a_{3}}\otimes W^{a'_{3}}&\to &
\mathbb{G}_{1; 2}\nn
w_{a_{3}}\otimes w_{a'_{3}}&\mapsto & 
(\Psi(\Y_{1}\otimes 
\Y_{2}))(w_{a_{3}}\otimes w_{a'_{3}}; 
z_{1}, z_{2}; \tau)
\end{eqnarray*}
for $a_{1}, a_{2}\in \A$, $i=1, \dots, 
N_{aa_{1}}^{a_{1}}$, $j=1, \dots, N_{a_{2}a_{2}'}^{a}$,
where 
$$(\Psi(\Y_{1}\otimes 
\Y_{2}))(w_{a_{3}}\otimes w_{a'_{3}}; 
z_{1}, z_{2}; \tau)
=0$$ 
when $a_{3}\ne a_{2}$ and 
\begin{eqnarray*}
\lefteqn{(\Psi(\Y_{1}\otimes 
\Y_{2}))(w_{a_{3}}\otimes w_{a'_{3}}; 
z_{1}, z_{2}; \tau)}\nn
&&=
E(\tr_{W^{a_{1}}}
\Y_{1}(\mathcal{U}(e^{2\pi iz_{2}})
\Y_{2}
(w_{a_{2}}, z_{1}-z_{2})w_{a'_{2}}, e^{2\pi i z_{2}})
q_{\tau}^{L(0)-\frac{c}{24}}).
\end{eqnarray*}

For $a\in \A$, 
let $\F^{a}_{1; 2}$ be the space spanned by linear maps of the form
$\Psi(\Y_{1}\otimes \Y_{2})$
for $a_{1}, a_{2}\in \A$ and for intertwining 
operators 
$\Y_{1}$ and $\Y_{2}$  
of types ${W^{a_{1}}\choose W^{a}W^{a_{1}}}$ and 
${W^{a}\choose W^{a_{2}}W^{a'_{2}}}$, respectively. 
For $a\in \A$, we also let $\F^{na}_{1; 2}$ be 
the space 
spanned by linear maps of 
the form $\Psi(\Y_{1}\otimes \Y_{2})$
for $a_{1}, a_{2}, a_{3}\in \A$, $a_{3}\ne a$,  and for intertwining 
operators 
$\Y_{1}$ and $\Y_{2}$  
of types ${W^{a_{1}}\choose W^{a_{3}}W^{a_{1}}}$ and 
${W^{a_{3}}\choose W^{a_{2}}W^{a'_{2}}}$, respectively. 
Let $\F_{1; 2}$ be the sum of $\F^{a}_{1; 2}$ for $a\in \A$. 
We
have:

\begin{prop}\label{projection}
For $a\in \A$, the intersection of $\F^{a}_{1; 2}$ and $\F^{na}_{1; 2}$ 
is $0$. In particular, 
$$\F_{1; 2}=\F^{a}_{1; 2}\oplus \F^{na}_{1; 2}$$
and there exists a projection $\pi: \F_{1; 2}\to \F^{a}_{1; 2}$. 
\end{prop}
\pf
By Proposition \ref{1-independence0}, 
$\Psi(\Y_{aa_{1}; k}^{a_{1}; (1)}\otimes \Y_{a_{2}a_{2}'; l}^{a_{3}; (2)})$
for $a_{1}, a_{2}, a_{3}\in \A$, $k=1, \dots, 
N_{a_{3}a_{1}}^{a_{1}}$ and $l=1, \dots, N_{a_{2}a_{2}'}^{a_{3}}$
are linearly independent. 
Thus the intersection of the space spanned by 
$\Psi(\Y_{aa_{1}; k}^{a_{1}; (1)}\otimes 
\Y_{a_{2}a_{2}'; l}^{a; (2)})$
for 
$a_{1}, a_{2}\in \A$, $k=1, \dots, 
N_{aa_{1}}^{a_{1}}$, $l=1, \dots, N_{a_{2}a_{2}'}^{a}$, 
and the space spanned 
by $\Psi(\Y_{a_{3}a_{1}; k}^{a_{1}; (1)}\otimes
 \Y_{a_{2}a_{2}'; l}^{a_{3}; (2)})$ for
$a_{1}, a_{2}, a_{3}\in \A$, $a_{3}\ne a$, 
$k=1, \dots, 
N_{a_{3}a_{1}}^{a_{1}}$ and 
$l=1, \dots, N_{a_{2}a_{2}'}^{a_{3}}$ are $0$.
\epfv

In the remaining part of this section, 
$\{\Y_{a_{1}a_{2}; i}^{a_{3};(p)}\;|\;i=1, \dots, 
N_{a_{1}a_{2}}^{a_{3}}\}$ for $p=1, 2, 3, 4, 5, 6$ 
are basis of the space of intertwining operators 
of type ${W^{a_{3}}\choose W^{a_{1}}W^{a_{2}}}$, 
$a_{1}, a_{2}, a_{3}\in \A$.

We have the following lemma which is a generalization of 
Lemma 4.2 in \cite{H10}:

\begin{lemma}
For $a_{1}, a_{2}, a\in \A$, $w_{a_{2}}\in W^{a_{2}}$ 
and $w_{a'_{2}}\in W^{a'_{2}}$, we have
\begin{eqnarray}\label{alpha0}
\lefteqn{(\Psi(\Y_{aa_{1}; p}^{a_{1};(1)}\otimes 
\Y_{a_{2}a'_{2}; q}^{a; (2)}))(w_{a_{2}}\otimes w_{a'_{2}}; 
z_{1}, z_{2}-1; \tau)}\nn
&&=\sum_{a\in \mathcal{A}}\sum_{i=1}^{N_{a_{2}a_{3}}^{a_{1}}}
\sum_{j=1}^{N_{a'_{2}a_{1}}^{a_{3}}}
\sum_{a_{4}\in \mathcal{A}}\sum_{k=1}^{N_{a_{4}a_{1}}^{a_{1}}}
\sum_{l=1}^{N_{a_{2}a'_{2}}^{a_{4}}}e^{-2\pi i(h_{a_{3}}-h_{a_{1}})}\cdot\nn
&&\quad\quad \quad \cdot F^{-1}(\Y_{aa_{1}; p}^{a_{1};{(1)}}\otimes 
\Y_{a_{2}a'_{2}; q}^{a;(2)};
\Y_{a_{2}a_{3}; i}^{a_{1};(3)}\otimes \Y_{a'_{2}a_{1}; j}^{a_{3};(4)})
\cdot \nn
&&\quad\quad \quad  \cdot
F(\Y_{a_{2}a_{3}; i}^{a_{1};(3)}\otimes \Y_{a'_{2}a_{1}; j}^{a_{3};(4)};
\Y_{a_{4}a_{1}; k}^{a_{1};(5)}\otimes \Y_{a_{2}a'_{2}; l}^{a_{4};(6)})\cdot \nn
&&\quad\quad\quad \cdot
E\bigg(\tr_{W^{a_{1}}}
\Y_{a_{4}a_{1}; k}^{a_{1};(5)}(\mathcal{U}(e^{2\pi iz_{2}})\cdot \nn
&&\quad\quad\quad \quad\quad\quad\quad\quad\quad\quad\cdot
\Y_{a_{2}a'_{2}; l}^{a_{4};(6)}(w_{a_{2}},z_{1}-z_{2})w_{a'_{2}},  
e^{2\pi i z_{2}})
q_{\tau}^{L(0)-\frac{c}{24}}\bigg)\nn
&&
\end{eqnarray}
and 
\begin{eqnarray}\label{beta0}
\lefteqn{(\Psi(\Y_{aa_{1}; i}^{a_{1};(1)}\otimes 
\Y_{a_{2}a'_{2}; j}^{a; (2)}))(w_{a_{2}}\otimes w_{a'_{2}}; 
z_{1}, z_{2}+\tau; \tau)}\nn
&&=\sum_{a_{3}\in \mathcal{A}}\sum_{i=1}^{N_{a_{2}a_{3}}^{a_{1}}}
\sum_{j=1}^{N_{a'_{2}a_{1}}^{a_{3}}}
\sum_{a_{4}\in \mathcal{A}}\sum_{k=1}^{N_{a_{4}a_{3}}^{a_{3}}}
\sum_{l=1}^{N_{a_{2}a'_{2}}^{a_{4}}}e^{\pi i (-2h_{a_{2}}+h_{a_{4}})}\nn
&&\quad \quad F^{-1}(\Y_{aa_{1}; p}^{a_{1};(p)}\otimes 
\Y_{a_{2}a'_{2}; j}^{a;(2)};
\sigma_{23}(\Y_{a_{2}a'_{1}; i}^{a'_{3};(3)}) \otimes 
\sigma_{13}(\Y_{a'_{3}a_{1}; j}^{a_{2};(4)}))\cdot \nn
&&\quad\quad \cdot 
F(\Y_{a_{2}a'_{1}; i}^{a'_{3};(3)}\otimes 
\sigma_{123}(\Y_{a'_{3}a_{1}; j}^{a_{2};(4)});
\Y_{a_{4}a'_{3}; k}^{a'_{3};(5)}\otimes \Y_{a_{2}a'_{2}; l}^{a_{4};(6)})\cdot \nn
&&\quad\quad \cdot 
E\bigg(\tr_{W^{a_{3}}}
\Y_{a_{4}a_{3}; k}^{a_{3};(5)}(\mathcal{U}(e^{2\pi iz_{2}})
\Y_{a_{2}a'_{2}; l}^{a_{4};(6)}(
w_{a_{2}}, z_{1}-z_{2})
w_{a'_{2}}, e^{2\pi iz_{2}})
q_{\tau}^{L(0)-\frac{c}{24}}
\bigg).\nn
&&
\end{eqnarray}
In particular,
for any $a_{1}, a_{2}\in \A$, the maps from 
$\coprod_{a_{3}\in \A}W^{a_{3}}\otimes W^{a'_{3}}$ to the space of
single-valued analytic functions on $\widetilde{M}_{1}^{2}$
given by 
\begin{eqnarray*}
w_{a_{3}}\otimes w_{a'_{3}}&\mapsto & 
(\Psi(\Y_{aa_{1}; p}^{a_{1};(1)}\otimes 
\Y_{a_{2}a'_{2}; q}^{a; (2)}))(w_{a_{3}}\otimes w_{a'_{3}}; 
z_{1}, z_{2}-1; \tau),\\
w_{a_{3}}\otimes w_{a'_{3}}&\mapsto & 
(\Psi(\Y_{aa_{1}; i}^{a_{1};(1)}\otimes 
\Y_{a_{2}a'_{2}; j}^{a; (2)}))(w_{a_{3}}\otimes w_{a'_{3}}; 
z_{1}, z_{2}+\tau; \tau)
\end{eqnarray*}
for $w_{a_{3}}\in W^{a_{3}}$ and $w_{a'_{3}}\in W^{a'_{3}}$ are
in $\mathcal{F}_{1;2}$.
\end{lemma}
\pf
Using the definition of $\Psi(\Y_{aa_{1}; p}^{a_{1};(1)}\otimes 
\Y_{a_{2}a'_{2}; q}^{a; (2)})$, 
the genus-one associativity properties ((2.9) and (2.10) in \cite{H10})
and
\begin{equation}\label{y-delta}
\Y(w_{a_{1}}, x)w_{a_{2}}\in x^{\Delta(\Y)}W^{a_{3}}[[x, x^{-1}]],
\end{equation}
where 
$$\Delta(\mathcal{Y})=h_{a_{3}}-h_{a_{1}}-h_{a_{2}},$$
we obtain
\begin{eqnarray*}
\lefteqn{(\Psi(\Y_{aa_{1}; p}^{a_{1};(1)}\otimes 
\Y_{a_{2}a'_{2}; q}^{a; (2)}))(w_{a_{2}}\otimes w_{a'_{2}}; 
z_{1}, z_{2}-1; \tau)}\nn
&&=E\bigg(\tr_{W^{a_{1}}}
\Y_{aa_{1}; p}^{a_{1}}(\mathcal{U}(e^{2\pi i(z_{2}-1)})\cdot\nn
&&\quad\quad\quad \quad\quad\quad\quad\cdot
\Y_{a_{2}a'_{2}; q}^{a}
(w_{a_{2}}, z_{1}-(z_{2}-1))w_{a'_{2}}, e^{2\pi i (z_{2}-1)})
q_{\tau}^{L(0)-\frac{c}{24}}\bigg)\nn
&&=\sum_{a_{3}\in \mathcal{A}}
\sum_{i=1}^{N_{a_{2}a_{3}}^{a_{1}}}\sum_{j=1}^{N_{a'_{2}a_{1}}^{a_{3}}}
F^{-1}(\Y_{aa_{1}; p}^{a_{1};(1)}\otimes \Y_{a_{2}a'_{2}; q}^{a;(2)};
\Y_{a_{2}a_{3}; i}^{a_{1};(3)}\otimes \Y_{a'_{2}a_{1}; j}^{a_{3};(4)})\cdot \nn
&&\quad\quad\quad \cdot 
E\bigg(\tr_{W^{a_{1}}}
\Y_{a_{2}a_{3}; i}^{a_{1}; (3)}(\mathcal{U}(e^{2\pi iz_{1}})w_{a_{2}}, e^{2\pi i z_{1}})\cdot \nn
&&\quad\quad\quad\quad\quad\quad\quad\quad\quad\quad\cdot 
\Y_{a'_{2}a_{1}; j}^{a_{3};(4)}(\mathcal{U}(e^{2\pi i(z_{2}-1)})w_{a'_{2}}, e^{2\pi i (z_{2}-1)})
q_{\tau}^{L(0)-\frac{c}{24}}\bigg)\nn
&&=\sum_{a_{3}\in \mathcal{A}}\sum_{i=1}^{N_{a_{2}a_{3}}^{a_{1}}}
\sum_{j=1}^{N_{a'_{2}a_{1}}^{a_{3}}}
F^{-1}(\Y_{aa_{1}; p}^{a_{1};(1)}\otimes \Y_{a_{2}a'_{2}; q}^{a;(2)};
\Y_{a_{2}a_{3}; i}^{a_{1};(3)}\otimes \Y_{a'_{2}a_{1}; j}^{a_{3};(4)})\cdot \nn
&&\quad\quad\quad \cdot 
E\bigg(\tr_{W^{a_{1}}}
\Y_{a_{2}a_{3}; i}^{a_{1}; (3)}(\mathcal{U}(e^{2\pi iz_{1}})w_{a_{2}}, 
e^{2\pi i z_{1}})\cdot \nn
&&\quad\quad\quad \quad\quad\quad\quad\quad\quad\quad\cdot
\Y_{a'_{2}a_{1}; j}^{a_{3}; (4)}(\mathcal{U}(e^{-2\pi i}
e^{2\pi iz_{2}})w_{a'_{2}},  e^{-2\pi i}e^{2\pi i z_{2}})
q_{\tau}^{L(0)-\frac{c}{24}}\bigg)\nn
&&=\sum_{a_{3}\in \mathcal{A}}\sum_{i=1}^{N_{a_{2}a_{3}}^{a_{1}}}
\sum_{j=1}^{N_{a'_{2}a_{1}}^{a_{3}}}e^{-2\pi i(h_{a_{3}}-h_{a_{1}})}
F^{-1}(\Y_{aa_{1}; p}^{a_{1}}\otimes \Y_{a_{2}a'_{2}; q}^{a};
\Y_{a_{2}a_{3}; i}^{a_{1};(3)}\otimes \Y_{a'_{2}a_{1}; j}^{a_{3};(4)})
\cdot \nn
&&\quad\quad\quad \cdot
E\bigg(\tr_{W^{a_{1}}}
\Y_{a_{2}a_{3}; i}^{a_{1};(3)}(\mathcal{U}(e^{2\pi iz_{1}})w_{a_{2}}, e^{2\pi i z_{1}})\cdot \nn
&&\quad\quad\quad \quad\quad\quad\quad\quad\quad\quad\cdot
\Y_{a'_{2}a_{1}; j}^{a_{3};(4)}(\mathcal{U}(
e^{2\pi iz_{2}})w_{a'_{2}},  e^{2\pi i z_{2}})
q_{\tau}^{L(0)-\frac{c}{24}}\bigg)\nn
&&=\sum_{a_{3}\in \mathcal{A}}\sum_{i=1}^{N_{a_{2}a_{3}}^{a_{1}}}
\sum_{j=1}^{N_{a'_{2}a_{1}}^{a_{3}}}e^{-2\pi i(h_{a_{3}}-h_{a_{1}})}
F^{-1}(\Y_{aa_{1}; p}^{a_{1}}\otimes \Y_{a_{2}a'_{2}; q}^{a};
\Y_{a_{2}a_{3}; i}^{a_{1};(3)}\otimes \Y_{a'_{2}a_{1}; j}^{a_{3};(4)})
\cdot \nn
&&\quad \cdot
\sum_{a_{4}\in \mathcal{A}}\sum_{k=1}^{N_{a_{4}a_{1}}^{a_{1}}}
\sum_{l=1}^{N_{a_{2}a'_{2}}^{a_{4}}}
F(\Y_{a_{2}a_{3}; i}^{a_{1};(3)}\otimes \Y_{a'_{2}a_{1}; j}^{a_{3};(4)};
\Y_{a_{4}a_{1}; k}^{a_{1};(5)}\otimes \Y_{a_{2}a'_{2}; l}^{a_{4};(6)})\cdot \nn
&&\quad\quad\quad \cdot
E\bigg(\tr_{W^{a_{1}}}
\Y_{a_{4}a_{1}; k}^{a_{1};(5)}(\mathcal{U}(e^{2\pi iz_{2}})\cdot \nn
&&\quad\quad\quad \quad\quad\quad\quad\quad\quad\quad\cdot
\Y_{a_{2}a'_{2}; l}^{a_{4};(6)}(w_{a_{2}},z_{1}-z_{2})w_{a'_{2}},  
e^{2\pi i z_{2}})
q_{\tau}^{L(0)-\frac{c}{24}}\bigg),
\end{eqnarray*}
proving (\ref{alpha0}).

The proof of (\ref{beta0}) is more complicated. 
Using the definition of $\Psi(\Y_{aa_{1}; p}^{a_{1};(1)}\otimes 
\Y_{a_{2}a'_{2}; q}^{a; (2)})$, 
the genus-one associativity properties ((2.9) and (2.10) in \cite{H10}), 
the $L(0)$-conjugation property, the property of traces,  and
(\ref{y-delta}), 
\begin{eqnarray}\label{beta1}
\lefteqn{(\Psi(\Y_{aa_{1}; p}^{a_{1};(1)}\otimes 
\Y_{a_{2}a'_{2}; q}^{a; (2)}))
(w_{a_{2}}\otimes w_{a'_{2}}; z_{1}, z_{2}+\tau; \tau)}\nn
&&=E\bigg(\tr_{W^{a_{1}}}
\Y_{aa_{1}; p}^{a_{1};(1)}(\mathcal{U}(e^{2\pi i(z_{2}+\tau)})\cdot\nn
&&\quad\quad\quad \quad\quad\quad\quad\cdot
\Y_{a_{2}a'_{2}; q}^{a;(2)}
(w_{a_{2}}, z_{1}-(z_{2}+\tau))w_{a'_{2}}, e^{2\pi i (z_{2}+\tau)})
q_{\tau}^{L(0)-\frac{c}{24}}\bigg)\nn
&&=\sum_{a_{3}\in \mathcal{A}}\sum_{i=1}^{N_{a_{2}a_{3}}^{a_{1}}}
\sum_{j=1}^{N_{a'_{2}a_{1}}^{a_{3}}}
F^{-1}(\Y_{aa_{1}; p}^{a_{1};(1)}\otimes \Y_{a_{2}a'_{2}; q}^{a; (2)};
\sigma_{23}(\Y_{a_{2}a'_{1}; i}^{a'_{3};(3)})\otimes 
\sigma_{13}(\Y_{a'_{3}a_{1}; j}^{a_{2};(4)}))\cdot \nn
&&\quad\quad\quad \cdot 
E\bigg(\tr_{W^{a_{1}}}
\sigma_{23}(\Y_{a_{2}a'_{1}; i}^{a'_{3};(3)})(
\mathcal{U}(e^{2\pi iz_{1}})w_{a_{2}}, e^{2\pi i z_{1}})\cdot \nn
&&\quad\quad\quad\quad\quad\quad\quad\quad\quad\quad\cdot 
\sigma_{13}(\Y_{a'_{3}a_{1}; j}^{a_{2};(4)})
(\mathcal{U}(e^{2\pi i(z_{2}+\tau)})w_{a'_{2}}, e^{2\pi i (z_{2}+\tau)})
q_{\tau}^{L(0)-\frac{c}{24}}\bigg)\nn
&&=\sum_{a_{3}\in \mathcal{A}}\sum_{i=1}^{N_{a_{2}a_{3}}^{a_{1}}}
\sum_{j=1}^{N_{a'_{2}a_{1}}^{a_{3}}}
F^{-1}(\Y_{aa_{1}; p}^{a_{1}; (1)}\otimes \Y_{a_{2}a'_{2}; q}^{a;(2)};
\sigma_{23}(\Y_{a_{2}a'_{1}; i}^{a'_{3};(3)}) \otimes 
\sigma_{13}(\Y_{a'_{3}a_{1}; j}^{a_{2};(4)}))\cdot \nn
&&\quad\quad\quad \cdot 
E\bigg(\tr_{W^{a_{1}}}
\sigma_{23}(\Y_{a_{2}a'_{1}; i}^{a'_{3};(3)})(
\mathcal{U}(e^{2\pi iz_{1}})w_{a_{2}}, e^{2\pi i z_{1}})\cdot \nn
&&\quad\quad\quad\quad\quad\quad\quad\quad\quad\quad\cdot 
\sigma_{13}(\Y_{a'_{3}a_{1}; j}^{a_{2};(4)})
(\mathcal{U}(q_{\tau}e^{2\pi iz_{2}})w_{a'_{2}}, q_{\tau}e^{2\pi iz_{2}})
q_{\tau}^{L(0)-\frac{c}{24}}\bigg)\nn
&&=\sum_{a_{3}\in \mathcal{A}}\sum_{i=1}^{N_{a_{2}a_{3}}^{a_{1}}}
\sum_{j=1}^{N_{a'_{2}a_{1}}^{a_{3}}}
F^{-1}(\Y_{aa_{1}; p}^{a_{1}; (1)}\otimes \Y_{a_{2}a'_{2}; q}^{a;(2)};
\sigma_{23}(\Y_{a_{2}a'_{1}; i}^{a'_{3};(3)}) \otimes 
\sigma_{13}(\Y_{a'_{3}a_{1}; j}^{a_{2};(4)}))\cdot \nn
&&\quad\quad\quad \cdot 
E\bigg(\tr_{W^{a_{1}}}
\sigma_{23}(\Y_{a_{2}a'_{1}; i}^{a'_{3};(3)})(
\mathcal{U}(e^{2\pi iz_{1}})w_{a_{2}}, e^{2\pi i z_{1}})\cdot \nn
&&\quad\quad\quad\quad\quad\quad\quad\quad\quad\quad\cdot 
q_{\tau}^{L(0)-\frac{c}{24}}
\sigma_{13}(\Y_{a'_{3}a_{1}; j}^{a_{2};(4)})(
\mathcal{U}(e^{2\pi iz_{2}})w_{a'_{2}}, e^{2\pi iz_{2}})
\bigg)\nn
&&=\sum_{a_{3}\in \mathcal{A}}\sum_{i=1}^{N_{a_{2}a_{3}}^{a_{1}}}
\sum_{j=1}^{N_{a'_{2}a_{1}}^{a_{3}}}
F^{-1}(\Y_{aa_{1}; p}^{a_{1}; (1)}\otimes \Y_{a_{2}a'_{2}; q}^{a;(2)};
\sigma_{23}(\Y_{a_{2}a'_{1}; i}^{a'_{3};(3)}) \otimes 
\sigma_{13}(\Y_{a'_{3}a_{1}; j}^{a_{2};(4)}))\cdot \nn
&&\quad\quad\quad \cdot 
E\bigg(\tr_{W^{a_{3}}}
\sigma_{13}(\Y_{a'_{3}a_{1}; j}^{a_{2};(4)})(
\mathcal{U}(e^{2\pi iz_{2}})w_{a'_{2}}, e^{2\pi iz_{2}})\cdot \nn
&&\quad\quad\quad\quad\quad\quad\quad\quad\quad\quad\cdot 
\sigma_{23}(\Y_{a_{2}a'_{1}; i}^{a'_{3};(3)})(
\mathcal{U}(e^{2\pi iz_{1}})w_{a_{2}}, e^{2\pi i z_{1}})
q_{\tau}^{L(0)-\frac{c}{24}}
\bigg).
\end{eqnarray}

Using (2.22) in \cite{H10}, the relations $\sigma_{23}^{2}=1$,
$\sigma_{23}\sigma_{13}=\sigma_{123}$ and the genus-one
associativity, we have
\begin{eqnarray}\label{beta2}
\lefteqn{E\bigg(\tr_{W^{a_{3}}}
\sigma_{13}(\Y_{a'_{3}a_{1}; j}^{a_{2};(4)})(
\mathcal{U}(e^{2\pi iz_{2}})w_{a'_{2}}, e^{2\pi iz_{2}})\cdot} \nn
&&\quad\quad\quad\quad\quad\quad\quad\quad\quad\quad\cdot 
\sigma_{23}(\Y_{a_{2}a'_{1}; i}^{a'_{3};(3)})(
\mathcal{U}(e^{2\pi iz_{1}})w_{a_{2}}, e^{2\pi i z_{1}})
q_{\tau}^{L(0)-\frac{c}{24}}
\bigg)\nn
&&=E\bigg(\tr_{W^{a_{3}}}
\sigma_{23}^{2}(\sigma_{13}(\Y_{a'_{3}a_{1}; j}^{a_{2};(4)}))(
\mathcal{U}(e^{2\pi iz_{2}})w_{a'_{2}}, e^{2\pi iz_{2}})\cdot \nn
&&\quad\quad\quad\quad\quad\quad\quad\quad\quad\quad\cdot 
\sigma_{23}(\Y_{a_{2}a'_{1}; i}^{a'_{3};(3)})(
\mathcal{U}(e^{2\pi iz_{1}})w_{a_{2}}, e^{2\pi i z_{1}})
q_{\tau}^{L(0)-\frac{c}{24}}
\bigg)\nn
&&=e^{-2\pi ih_{a_{2}}}
E\bigg(\tr_{W^{a'_{3}}}
\Y_{a_{2}a'_{1}; i}^{a'_{3};(3)}\bigg(
\mathcal{U}(e^{2\pi iz_{1}})e^{\pi iL(0)}
w_{a_{2}}, e^{-2\pi i z_{1}}\bigg)\cdot \nn
&&\quad\quad\quad\quad\quad\quad\quad\cdot 
\sigma_{23}(\sigma_{13}(\Y_{a'_{3}a_{1}; j}^{a_{2};(4)}))
\bigg(\mathcal{U}(e^{2\pi iz_{2}})e^{\pi iL(0)}w_{a'_{2}}, e^{-2\pi iz_{2}}\bigg)
q_{\tau}^{L(0)-\frac{c}{24}}
\bigg)\nn
&&=e^{-2\pi ih_{a_{2}}}
E\bigg(\tr_{W^{a'_{3}}}
\Y_{a_{2}a'_{1}; i}^{a'_{3};(3)}\bigg(\mathcal{U}(e^{-2\pi iz_{1}})e^{\pi i L(0)}
w_{a_{2}}, e^{-2\pi i z_{1}}\bigg)\cdot \nn
&&\quad\quad\quad\quad\quad\quad\quad\cdot 
\sigma_{123}(\Y_{a'_{3}a_{1}; j}^{a_{2};(4)})\bigg(\mathcal{U}(e^{-2\pi iz_{2}})e^{\pi i L(0)}
w_{a'_{2}}, e^{-2\pi iz_{2}}\bigg)
q_{\tau}^{L(0)-\frac{c}{24}}
\bigg)\nn
&&
\end{eqnarray}

We now prove 
\begin{eqnarray}\label{beta3}
\lefteqn{E\bigg(\tr_{W^{a'_{3}}}
\Y_{a_{2}a'_{1}; i}^{a'_{3};(3)}\bigg(\mathcal{U}(e^{-2\pi iz_{1}})e^{\pi i L(0)}
w_{a_{2}}, e^{-2\pi i z_{1}}\bigg)\cdot} \nn
&&\quad\quad\quad\quad\quad\quad\quad\cdot 
\sigma_{123}(\Y_{a'_{3}a_{1}; j}^{a_{2};(4)})
\bigg(\mathcal{U}(e^{-2\pi iz_{2}})e^{\pi i L(0)}
w_{a'_{2}}, e^{-2\pi iz_{2}}\bigg)
q_{\tau}^{L(0)-\frac{c}{24}}
\bigg)\nn
&&=\sum_{a_{4}\in \mathcal{A}}\sum_{k=1}^{N_{a_{4}a_{3}}^{a_{3}}}
\sum_{l=1}^{N_{a_{2}a'_{2}}^{a_{4}}}
F(\Y_{a_{2}a'_{1}; i}^{a'_{3};(3)}\otimes 
\sigma_{123}(\Y_{a'_{3}a_{1}; j}^{a_{2};(4)});
\Y_{a_{4}a'_{3}; k}^{a'_{3};(5)}\otimes \Y_{a_{2}a'_{2}; l}^{a_{4};(6)})\cdot \nn
&&\quad\quad\quad\cdot E\bigg(\tr_{W^{a'_{3}}}
\Y_{a_{4}a'_{3}; k}^{a'_{3};(5)}\bigg(\mathcal{U}(e^{-2\pi iz_{2}})
\Y_{a_{2}a'_{2}; l}^{a_{4};(6)}\bigg(e^{\pi i L(0)}
w_{a_{2}}, e^{\pi i}(z_{1}-z_{2})\bigg)\cdot \nn
&&\quad\quad\quad\quad\quad\quad\quad\quad\quad
\quad\quad\quad\quad\quad\quad\quad\quad\quad\cdot 
e^{\pi i L(0)}
w_{a'_{2}}, e^{-2\pi iz_{2}}\bigg)
q_{\tau}^{L(0)-\frac{c}{24}}
\bigg).\nn
&&
\end{eqnarray}
To prove (\ref{beta3}), we need only prove 
the restrictions of both sides to a subregion of 
$\widetilde{M}_{1}^{2}$ are equal. 
So we need only prove that
\begin{eqnarray}\label{beta3.5}
\lefteqn{\tr_{W^{a'_{3}}}
\Y_{a_{2}a'_{1}; i}^{a'_{3};(3)}\bigg(\mathcal{U}(e^{-2\pi iz_{1}})e^{\pi i L(0)}
w_{a_{2}}, e^{-2\pi i z_{1}}\bigg)\cdot} \nn
&&\quad\quad\quad\quad\quad\quad\quad\cdot 
\sigma_{123}(\Y_{a'_{3}a_{1}; j}^{a_{2};(4)})
\bigg(\mathcal{U}(e^{-2\pi iz_{2}})e^{\pi i L(0)}
w_{a'_{2}}, e^{-2\pi iz_{2}}\bigg)
q_{\tau}^{L(0)-\frac{c}{24}}\nn
&&=\sum_{a_{4}\in \mathcal{A}}\sum_{k=1}^{N_{a_{4}a_{3}}^{a_{3}}}
\sum_{l=1}^{N_{a_{2}a'_{2}}^{a_{4}}}
F(\Y_{a_{2}a'_{1}; i}^{a'_{3};(3)}\otimes 
\sigma_{123}(\Y_{a'_{3}a_{1}; j}^{a_{2}; (4)});
\Y_{a_{4}a'_{3}; k}^{a'_{3};(5)}\otimes \Y_{a_{2}a'_{2}; l}^{a_{4};(6)})\cdot \nn
&&\quad\quad\quad\cdot \tr_{W^{a'_{3}}}
\Y_{a_{4}a'_{3}; k}^{a'_{3};(5)}\bigg(\mathcal{U}(e^{-2\pi iz_{2}})
\Y_{a_{2}a'_{2}; l}^{a_{4};(6)}\bigg(e^{\pi i L(0)}
w_{a_{2}}, e^{\pi i}(z_{1}-z_{2})\bigg)\cdot \nn
&&\quad\quad\quad\quad\quad\quad\quad\quad\quad
\quad\quad\quad\quad\quad\quad\quad\quad\quad\cdot 
e^{\pi i L(0)}
w_{a'_{2}}, e^{-2\pi iz_{2}}\bigg)
q_{\tau}^{L(0)-\frac{c}{24}}
\bigg)\nn
&&
\end{eqnarray}
holds when $|q_{\tau}|<|e^{-2\pi iz_{2}}|<|e^{-2\pi iz_{1}}|<1$
and $0<|e^{2\pi i(-z_{1}+z_{2})}-1|<1$. 
From the genus-one associativity ((2.9) in \cite{H10}), 
we see  that in this region the left-hand side of 
(\ref{beta3.5}) is equal to 
\begin{eqnarray}\label{beta3.7}
\lefteqn{\sum_{a_{4}\in \mathcal{A}}\sum_{k=1}^{N_{a_{4}a_{3}}^{a_{3}}}
\sum_{l=1}^{N_{a_{2}a'_{2}}^{a_{4}}}
F(\Y_{a_{2}a'_{1}; i}^{a'_{3};(3)}\otimes 
\sigma_{123}(\Y_{a'_{3}a_{1}; j}^{a_{2};(4)});
\Y_{a_{4}a'_{3}; k}^{a'_{3};(5)}\otimes 
\Y_{a_{2}a'_{2}; l}^{a_{4};(6)})\cdot} \nn
&&\quad\quad\quad\cdot \tr_{W^{a'_{3}}}
\Y_{a_{4}a'_{3}; k}^{a'_{3};(5)}\bigg(\mathcal{U}(e^{-2\pi iz_{2}})
\Y_{a_{2}a'_{2}; l}^{a_{4};(6)}\bigg(e^{\pi i L(0)}
w_{a_{2}}, (-z_{1}+z_{2})\bigg)\cdot \nn
&&\quad\quad\quad\quad\quad\quad\quad\quad\quad
\quad\quad\quad\quad\quad\quad\quad\quad\quad\cdot 
e^{\pi i L(0)}
w_{a'_{2}}, e^{-2\pi iz_{2}}\bigg)
q_{\tau}^{L(0)-\frac{c}{24}}
\bigg)\nn
&&
\end{eqnarray}
Now in this region, because 
$|e^{-2\pi iz_{2}}|<|e^{-2\pi iz_{1}}|$, 
the imaginary part of $z_{1}$
must be bigger than the imaginary part of $z_{2}$. 
Thus $z_{1}-z_{2}$ is in the upper half plane. 
This means that  $\arg (z_{1}-z_{2})<\pi$ and 
$\arg (z_{1}-z_{2}) +\pi <2\pi$. So we have
$\arg (-(z_{1}-z_{2}))= \arg (z_{1}-z_{2})+\pi$.
Now for any $n\in \C$, by our convention,
\begin{eqnarray*}
(-z_{1}+z_{2})^{n}&=&e^{n\log (-z_{1}+z_{2})}\nn
&=&e^{n\log (-(z_{1}-z_{2}))}\nn
&=&e^{n(\log |-(z_{1}-z_{2})|+i\arg (-(z_{1}-z_{2})))}\nn
&=&e^{n(\log |(z_{1}-z_{2})|+i\arg (z_{1}-z_{2})+\pi i)}\nn
&=&e^{n(\log (z_{1}-z_{2})+\pi i)}\nn
&=&(e^{\pi i}(z_{1}-z_{2}))^{n}.
\end{eqnarray*}
This shows that indeed when 
$|q_{\tau}|<|e^{-2\pi iz_{2}}|<|e^{-2\pi iz_{1}}|<1$
and $0<|e^{2\pi i(-z_{1}+z_{2})}-1|<1$, 
(\ref{beta3.7}) is equal to the right-hand side of 
(\ref{beta3.5}) and  (\ref{beta3.5}) holds. 
Consequently, we obtain (\ref{beta3}).

Using (2.22) in \cite{H10} and the $L(0)$-conjugation formula, 
we have
\begin{eqnarray}\label{beta4}
\lefteqn{E\bigg(\tr_{W^{a'_{3}}}
\Y_{a_{4}a'_{3}; k}^{a'_{3};(5)}\bigg(\mathcal{U}(e^{-2\pi iz_{2}})
\Y_{a_{2}a'_{2}; l}^{a_{4};(6)}\bigg(e^{\pi i L(0)}
w_{a_{2}}, e^{\pi i}(z_{1}-z_{2})\bigg)\cdot} \nn
&&\quad\quad\quad\quad\quad\quad\quad\quad\quad
\quad\quad\quad\quad\quad\quad\quad\quad\quad\cdot 
e^{\pi i L(0)}
w_{a'_{2}}, e^{-2\pi iz_{2}}\bigg)
q_{\tau}^{L(0)-\frac{c}{24}}
\bigg)\nn
&&=e^{\pi i h_{a_{4}}}E\bigg(\tr_{W^{a_{3}}}
\Y_{a_{4}a_{3}; k}^{a_{3};(5)}\bigg(\mathcal{U}(e^{2\pi iz_{2}})e^{-\pi iL(0)}\cdot \nn
&&\quad\quad\quad\quad\quad
\cdot 
\Y_{a_{2}a'_{2}; l}^{a_{4};(6)}\bigg(e^{\pi i L(0)}
w_{a_{2}}, e^{\pi i}(z_{1}-z_{2})\bigg)e^{\pi i L(0)}
w_{a'_{2}}, e^{2\pi iz_{2}}\bigg)
q_{\tau}^{L(0)-\frac{c}{24}}
\bigg)\nn
&&=e^{\pi ih_{a_{4}}}
E\bigg(\tr_{W^{a_{3}}}
\Y_{a_{4}a_{3}; k}^{a_{3};(5)}(\mathcal{U}(e^{2\pi iz_{2}})
\Y_{a_{2}a'_{2}; l}^{a_{4};(6)}(
w_{a_{2}}, z_{1}-z_{2})
w_{a'_{2}}, e^{2\pi iz_{2}})
q_{\tau}^{L(0)-\frac{c}{24}}
\bigg).\nn
&&
\end{eqnarray}
Combining (\ref{beta1}), (\ref{beta2}), (\ref{beta3}) and 
(\ref{beta4}), we obtain (\ref{beta0}).
\epfv

For $a_{1}, a_{2}, a\in \A$, we define 
\begin{eqnarray}
\alpha(\Psi(\Y_{aa_{1}; p}^{a_{1};(1)}\otimes 
\Y_{a_{2}a'_{2}; q}^{a; (2)})):\coprod_{a_{3}\in \A}W^{a_{3}}
\otimes W^{a'_{3}}&\to &\mathbb{G}_{1;2}\label{alpha-def}\\
\beta(\Psi(\Y_{aa_{1}; p}^{a_{1};(1)}\otimes 
\Y_{a_{2}a'_{2}; q}^{a; (2)})):\coprod_{a_{3}\in \A}W^{a_{3}}
\otimes W^{a'_{3}}&\to &\mathbb{G}_{1;2}\label{beta-def}
\end{eqnarray}
by
\begin{eqnarray*}
\lefteqn{(\alpha(\Psi(\Y_{aa_{1}; p}^{a_{1};(1)}\otimes 
\Y_{a_{2}a'_{2}; q}^{a; (2)})))
(w_{a_{3}}\otimes w_{a'_{3}})}\nn
&&=\left\{\begin{array}{ll}(\Psi(\Y_{aa_{1}; p}^{a_{1};(1)}\otimes 
\Y_{a_{2}a'_{2}; q}^{a; (2)}))(w_{a_{2}}\otimes w_{a'_{2}}; 
z_{1}, z_{2}-1; \tau)&a_{3}=a_{2},\\
0&a_{3}\ne a_{2}\end{array}\right.,\\
\lefteqn{(\beta(\Psi(\Y_{aa_{1}; p}^{a_{1};(1)}\otimes 
\Y_{a_{2}a'_{2}; q}^{a; (2)})))
(w_{a_{3}}\otimes w_{a'_{3}})}\nn
&&=\left\{\begin{array}{ll}(\Psi(\Y_{aa_{1}; p}^{a_{1};(1)}\otimes 
\Y_{a_{2}a'_{2}; q}^{a; (2)}))
(w_{a_{2}}\otimes w_{a'_{2}}; z_{1}, z_{2}+\tau; \tau)&a_{3}=a_{2},\\
0&a_{3}\ne a_{2}\end{array}\right.
\end{eqnarray*}
for $a_{3}\in \A$, $w_{a_{3}}\in W^{a}$ and $w_{a'_{3}}\in W^{a'}$.

\begin{prop}
For $a_{1}, a_{2}\in \A$, we have 
\begin{eqnarray}\label{alpha}
\lefteqn{\alpha(\Psi(\Y_{aa_{1}; p}^{a_{1};(1)}\otimes 
\Y_{a_{2}a'_{2}; q}^{a; (2)}))}\nn
&&=e^{-2\pi ih_{a_{2}}}
\sum_{a_{4}\in \mathcal{A}}\sum_{k=1}^{N_{a_{4}a_{1}}^{a_{1}}}
\sum_{l=1}^{N_{a_{2}a'_{2}}^{a_{4}}}\nn
&&\quad\quad \quad  (B^{(-1)})^{2}
(\sigma_{12}(\Y_{aa_{1}; p}^{a_{1};{(1)}})\otimes 
\sigma_{12}(\Y_{a_{2}a'_{2}; q}^{a;(2)});
\sigma_{12}(\Y_{a_{4}a_{1}; k}^{a_{1};(5)})\otimes 
\sigma_{12}(\Y_{a_{2}a'_{2}; l}^{a_{4};(6)}))\cdot \nn
&&\quad\quad\quad \cdot
(\Psi(\Y_{a_{4}a_{1}; k}^{a_{1};(5)}\otimes 
\Y_{a_{2}a'_{2}; l}^{a_{4}; (6)}))
(w_{a_{2}}\otimes w_{a'_{2}}; z_{1}, z_{2}; \tau)
\end{eqnarray}
and 
\begin{eqnarray}\label{beta}
\lefteqn{\beta(\Psi(\Y_{aa_{1}; p}^{a_{1};(1)}\otimes 
\Y_{a_{2}a'_{2}; q}^{a; (2)}))}\nn
&&=e^{-2\pi ih_{a_{2}}}
\sum_{a_{3}\in \mathcal{A}}\sum_{i=1}^{N_{a_{2}a_{3}}^{a_{1}}}
\sum_{j=1}^{N_{a'_{2}a_{1}}^{a_{3}}}
\sum_{a_{4}\in \mathcal{A}}\sum_{k=1}^{N_{a_{4}a_{3}}^{a_{3}}}
\sum_{l=1}^{N_{a_{2}a'_{2}}^{a_{4}}}e^{\pi ih_{a_{4}}}\nn
&&\quad \quad F(\sigma_{12}(\Y_{aa_{1}; p}^{a_{1};(p)})\otimes 
\sigma_{12}(\Y_{a_{2}a'_{2}; j}^{a;(2)});
\sigma_{123}(\Y_{a_{2}a'_{1}; i}^{a'_{3};(3)}) \otimes 
\sigma_{132}(\Y_{a'_{3}a_{1}; j}^{a_{2};(4)}))\cdot \nn
&&\quad\quad \cdot 
F(\Y_{a_{2}a'_{1}; i}^{a'_{3};(3)}\otimes 
\sigma_{123}(\Y_{a'_{3}a_{1}; j}^{a_{2};(4)});
\Y_{a_{4}a'_{3}; k}^{a'_{3};(5)}\otimes \Y_{a_{2}a'_{2}; l}^{a_{4};(6)})\cdot \nn
&&\quad\quad \cdot 
(\Psi(\Y_{a_{4}a_{3}; k}^{a_{3};(5)}\otimes
\Y_{a_{2}a'_{2}; l}^{a_{4};(6)}))(w_{a_{2}}\otimes
w_{a'_{2}}; z_{1}, z_{2}; \tau)
\end{eqnarray}
\end{prop}
\pf
Using the definitions of $\alpha$
and 
$$(\Psi(\Y_{aa_{1}; p}^{a_{1};(1)}\otimes 
\Y_{a_{2}a'_{2}; q}^{a; (2)}))
(w_{a_{2}}\otimes w_{a'_{2}}; z_{1}, z_{2}; \tau),$$
(\ref{alpha0}) and Proposition 3.1 in \cite{H10}, we have
\begin{eqnarray*}
\lefteqn{(\alpha(\Psi(\Y_{aa_{1}; p}^{a_{1};(1)}\otimes 
\Y_{a_{2}a'_{2}; q}^{a; (2)})))
(w_{a_{2}}\otimes w_{a'_{2}}; z_{1}, z_{2}; \tau)}\nn
&&=(\Psi(\Y_{aa_{1}; p}^{a_{1};(1)}\otimes 
\Y_{a_{2}a'_{2}; q}^{a; (2)}))
(w_{a_{2}}\otimes w_{a'_{2}}; z_{1}, z_{2}-1; \tau)\nn
&&=\sum_{a\in \mathcal{A}}\sum_{i=1}^{N_{a_{2}a_{3}}^{a_{1}}}
\sum_{j=1}^{N_{a'_{2}a_{1}}^{a_{3}}}
\sum_{a_{4}\in \mathcal{A}}\sum_{k=1}^{N_{a_{4}a_{1}}^{a_{1}}}
\sum_{l=1}^{N_{a_{2}a'_{2}}^{a_{4}}}e^{-2\pi i(h_{a_{3}}-h_{a_{1}})}\cdot\nn
&&\quad\quad \quad \cdot F^{-1}(\Y_{aa_{1}; p}^{a_{1};{(1)}}\otimes 
\Y_{a_{2}a'_{2}; q}^{a;(2)};
\Y_{a_{2}a_{3}; i}^{a_{1};(3)}\otimes \Y_{a'_{2}a_{1}; j}^{a_{3};(4)})
\cdot \nn
&&\quad\quad \quad  \cdot
F(\Y_{a_{2}a_{3}; i}^{a_{1};(3)}\otimes \Y_{a'_{2}a_{1}; j}^{a_{3};(4)};
\Y_{a_{4}a_{1}; k}^{a_{1};(5)}\otimes \Y_{a_{2}a'_{2}; l}^{a_{4};(6)})\cdot \nn
&&\quad\quad\quad \cdot
E\bigg(\tr_{W^{a_{1}}}
\Y_{a_{4}a_{1}; k}^{a_{1};(5)}(\mathcal{U}(e^{2\pi iz_{2}})\cdot \nn
&&\quad\quad\quad \quad\quad\quad\quad\quad\quad\quad\cdot
\Y_{a_{2}a'_{2}; l}^{a_{4};(6)}(w_{a_{2}},z_{1}-z_{2})w_{a'_{2}},  
e^{2\pi i z_{2}})
q_{\tau}^{L(0)-\frac{c}{24}}\bigg)\nn
&&=e^{-2\pi ih_{a_{2}}}
\sum_{a\in \mathcal{A}}\sum_{i=1}^{N_{a_{2}a_{3}}^{a_{1}}}
\sum_{j=1}^{N_{a'_{2}a_{1}}^{a_{3}}}
\sum_{a_{4}\in \mathcal{A}}\sum_{k=1}^{N_{a_{4}a_{1}}^{a_{1}}}
\sum_{l=1}^{N_{a_{2}a'_{2}}^{a_{4}}}\cdot\nn
&&\quad\quad \quad \cdot F(\sigma_{12}(\Y_{aa_{1}; p}^{a_{1};{(1)}})\otimes 
\sigma_{12}(\Y_{a_{2}a'_{2}; q}^{a;(2)});
\sigma_{12}(\Y_{a_{2}a_{3}; i}^{a_{1};(3)})\otimes 
\sigma_{12}(\Y_{a'_{2}a_{1}; j}^{a_{3};(4)}))
\cdot \nn
&&\quad\quad \quad \cdot e^{-2\pi i(h_{a_{3}}-h_{a_{1}}-h_{a_{2}})}
\cdot \nn
&&\quad\quad \quad  \cdot
F^{-1}(\sigma_{12}(\Y_{a_{2}a_{3}; i}^{a_{1};(3)})\otimes 
\sigma_{12}(\Y_{a'_{2}a_{1}; j}^{a_{3};(4)});
\sigma_{12}(\Y_{a_{4}a_{1}; k}^{a_{1};(5)})\otimes 
\sigma_{12}(\Y_{a_{2}a'_{2}; l}^{a_{4};(6)}))\cdot \nn
&&\quad\quad\quad \cdot
E\bigg(\tr_{W^{a_{1}}}
\Y_{a_{4}a_{1}; k}^{a_{1};(5)}(\mathcal{U}(e^{2\pi iz_{2}})\cdot \nn
&&\quad\quad\quad \quad\quad\quad\quad\quad\quad\quad\cdot
\Y_{a_{2}a'_{2}; l}^{a_{4};(6)}(w_{a_{2}},z_{1}-z_{2})w_{a'_{2}},  
e^{2\pi i z_{2}})
q_{\tau}^{L(0)-\frac{c}{24}}\bigg)\nn
&&=e^{-2\pi ih_{a_{2}}}
\sum_{a_{4}\in \mathcal{A}}\sum_{k=1}^{N_{a_{4}a_{1}}^{a_{1}}}
\sum_{l=1}^{N_{a_{2}a'_{2}}^{a_{4}}}\nn
&&\quad\quad \quad  (B^{(-1)})^{2}
(\sigma_{12}(\Y_{aa_{1}; p}^{a_{1};{(1)}})\otimes 
\sigma_{12}(\Y_{a_{2}a'_{2}; q}^{a;(2)});
\sigma_{12}(\Y_{a_{4}a_{1}; k}^{a_{1};(5)})\otimes 
\sigma_{12}(\Y_{a_{2}a'_{2}; l}^{a_{4};(6)}))\cdot \nn
&&\quad\quad\quad \cdot
(\Psi(\Y_{a_{4}a_{1}; k}^{a_{1};(5)}\otimes 
\Y_{a_{2}a'_{2}; l}^{a_{4}; (6)}))
(w_{a_{2}}\otimes w_{a'_{2}}; z_{1}, z_{2}; \tau),
\end{eqnarray*}
proving (\ref{alpha}). 

Using the definitions of $\beta$, 
$$\Psi(\Y_{aa_{1}; p}^{a_{1};(1)}\otimes 
\Y_{a_{2}a'_{2}; q}^{a; (2)}))
(w_{a_{2}}\otimes w_{a'_{2}}; z_{1}, z_{2}; \tau),$$
and (\ref{beta0}), 
we have
\begin{eqnarray*}
\lefteqn{(\beta(\Psi(\Y_{aa_{1}; p}^{a_{1};(1)}\otimes 
\Y_{a_{2}a'_{2}; q}^{a; (2)}))
(w_{a_{2}}\otimes w_{a'_{2}}; z_{1}, z_{2}; \tau)))
(w_{a_{2}}, w_{a'_{2}}; z_{1}, z_{2}; \tau)}\nn
&&=\Psi(\Y_{aa_{1}; p}^{a_{1};(1)}\otimes 
\Y_{a_{2}a'_{2}; q}^{a; (2)}))
(w_{a_{2}}\otimes w_{a'_{2}}; z_{1}, z_{2}+\tau; \tau)\nn
&&=\sum_{a_{3}\in \mathcal{A}}\sum_{i=1}^{N_{a_{2}a_{3}}^{a_{1}}}
\sum_{j=1}^{N_{a'_{2}a_{1}}^{a_{3}}}
\sum_{a_{4}\in \mathcal{A}}\sum_{k=1}^{N_{a_{4}a_{3}}^{a_{3}}}
\sum_{l=1}^{N_{a_{2}a'_{2}}^{a_{4}}}e^{\pi i (-2h_{a_{2}}+h_{a_{4}})}\nn
&&\quad \quad F^{-1}(\Y_{aa_{1}; p}^{a_{1};(p)}\otimes 
\Y_{a_{2}a'_{2}; j}^{a;(2)};
\sigma_{23}(\Y_{a_{2}a'_{1}; i}^{a'_{3};(3)}) \otimes 
\sigma_{13}(\Y_{a'_{3}a_{1}; j}^{a_{2};(4)}))\cdot \nn
&&\quad\quad \cdot 
F(\Y_{a_{2}a'_{1}; i}^{a'_{3};(3)}\otimes 
\sigma_{123}(\Y_{a'_{3}a_{1}; j}^{a_{2};(4)});
\Y_{a_{4}a'_{3}; k}^{a'_{3};(5)}\otimes \Y_{a_{2}a'_{2}; l}^{a_{4};(6)})\cdot \nn
&&\quad\quad \cdot 
E\bigg(\tr_{W^{a_{3}}}
\Y_{a_{4}a_{3}; k}^{a_{3};(5)}(\mathcal{U}(e^{2\pi iz_{2}})
\Y_{a_{2}a'_{2}; l}^{a_{4};(6)}(
w_{a_{2}}, z_{1}-z_{2})
w_{a'_{2}}, e^{2\pi iz_{2}})
q_{\tau}^{L(0)-\frac{c}{24}}
\bigg)\nn
&&=e^{-2\pi ih_{a_{2}}}
\sum_{a_{3}\in \mathcal{A}}\sum_{i=1}^{N_{a_{2}a_{3}}^{a_{1}}}
\sum_{j=1}^{N_{a'_{2}a_{1}}^{a_{3}}}
\sum_{a_{4}\in \mathcal{A}}\sum_{k=1}^{N_{a_{4}a_{3}}^{a_{3}}}
\sum_{l=1}^{N_{a_{2}a'_{2}}^{a_{4}}}e^{\pi ih_{a_{4}}}\nn
&&\quad \quad F(\sigma_{12}(\Y_{aa_{1}; p}^{a_{1};(p)})\otimes 
\sigma_{12}(\Y_{a_{2}a'_{2}; j}^{a;(2)});
\sigma_{123}(\Y_{a_{2}a'_{1}; i}^{a'_{3};(3)}) \otimes 
\sigma_{132}(\Y_{a'_{3}a_{1}; j}^{a_{2};(4)}))\cdot \nn
&&\quad\quad \cdot 
F(\Y_{a_{2}a'_{1}; i}^{a'_{3};(3)}\otimes 
\sigma_{123}(\Y_{a'_{3}a_{1}; j}^{a_{2};(4)});
\Y_{a_{4}a'_{3}; k}^{a'_{3};(5)}\otimes \Y_{a_{2}a'_{2}; l}^{a_{4};(6)})\cdot \nn
&&\quad\quad \cdot 
(\Psi(\Y_{a_{4}a_{3}; k}^{a_{3};(5)}\otimes
\Y_{a_{2}a'_{2}; l}^{a_{4};(6)}))(w_{a_{2}}\otimes
w_{a'_{2}}; z_{1}, z_{2}; \tau),
\end{eqnarray*}
proving (\ref{beta}). 
\epfv

By Proposition 2.2 in \cite{H10}, fix any basis 
$\{\Y_{aa_{1}; p}^{a_{1};(1)}\;|\;p=1, \dots, N_{aa_{1}}^{a_{1}}\}$
and $\{\Y_{a_{2}a'_{2}; q}^{a; (2)})\;|\;
q=1, \dots, N_{a_{2}a'_{2}}^{a}\}$,
$$(\Psi(\Y_{aa_{1}; p}^{a_{1};(1)}\otimes 
\Y_{a_{2}a'_{2}; q}^{a; (2)}))
(w_{a_{2}}\otimes w_{a'_{2}}; z_{1}, z_{2}; \tau)$$
for $a_{1}, a_{2}, a\in \A$, $p=1, \dots, N_{aa_{1}}^{a_{1}}$
and $q=1, \dots, N_{a_{2}a'_{2}}^{a}$ 
form a basis of $\mathcal{F}_{1;2}$. 
We use 
$$\alpha(\Psi(\Y_{aa_{1}; p}^{a_{1};(1)}\otimes 
\Y_{a_{2}a'_{2}; q}^{a; (2)}); 
\Psi(\Y_{a_{4}a_{3}; k}^{a_{3};(5)}\otimes
\Y_{a_{2}a'_{2}; l}^{a_{4};(6)}))$$
and 
$$\beta(\Psi(\Y_{aa_{1}; p}^{a_{1};(1)}\otimes 
\Y_{a_{2}a'_{2}; q}^{a; (2)}); 
\Psi(\Y_{a_{4}a_{3}; k}^{a_{3};(5)}\otimes
\Y_{a_{2}a'_{2}; l}^{a_{4};(6)}))$$
to denote the matrix elements of $\alpha$ and $\beta$, respectively,
under the basis 
$$\Psi(\Y_{aa_{1}; p}^{a_{1};(1)}\otimes 
\Y_{a_{2}a'_{2}; q}^{a; (2)})$$
and
$$\Psi(\Y_{a_{4}a_{3}; k}^{a_{3};(5)}\otimes
\Y_{a_{2}a'_{2}; l}^{a_{4};(6)}).$$

\begin{cor}
The matrix elements of $\alpha$ and 
$\beta$
are given by
\begin{eqnarray}\label{alpha-1}
\lefteqn{\alpha(\Psi(\Y_{aa_{1}; p}^{a_{1};(1)}\otimes 
\Y_{a_{2}a'_{2}; q}^{a; (2)}); 
\Psi(\Y_{a_{4}a_{3}; k}^{a_{3};(5)}\otimes
\Y_{a_{2}a'_{2}; l}^{a_{4};(6)}))}\nn
&&=\delta_{a_{3}a_{1}}e^{-2\pi ih_{a_{2}}}\cdot\nn
&&\quad\quad \quad \cdot (B^{(-1)})^{2}
(\sigma_{12}(\Y_{aa_{1}; p}^{a_{1};{(1)}})\otimes 
\sigma_{12}(\Y_{a_{2}a'_{2}; q}^{a;(2)});
\sigma_{12}(\Y_{a_{4}a_{1}; k}^{a_{1};(5)})\otimes 
\sigma_{12}(\Y_{a_{2}a'_{2}; l}^{a_{4};(6)}))\nn
&&
\end{eqnarray}
and
\begin{eqnarray}\label{beta-1}
\lefteqn{\beta(\Psi(\Y_{aa_{1}; p}^{a_{1};(1)}\otimes 
\Y_{a_{2}a'_{2}; q}^{a; (2)}); 
\Psi(\Y_{a_{4}a_{3}; k}^{a_{3};(5)}\otimes
\Y_{a_{2}a'_{2}; l}^{a_{4};(6)}))}\nn
&&=e^{-2\pi ih_{a_{2}}}
\sum_{i=1}^{N_{a_{2}a_{3}}^{a_{1}}}
\sum_{j=1}^{N_{a'_{2}a_{1}}^{a_{3}}}
e^{\pi ih_{a_{4}}}\cdot\nn
&&\quad \quad \cdot F(\sigma_{12}(\Y_{aa_{1}; p}^{a_{1};(1)})\otimes 
\sigma_{12}(\Y_{a_{2}a'_{2}; q}^{a;(2)});
\sigma_{123}(\Y_{a_{2}a'_{1}; i}^{a'_{3};(3)}) \otimes 
\sigma_{132}(\Y_{a'_{3}a_{1}; j}^{a_{2};(4)}))\cdot \nn
&&\quad\quad \cdot 
F(\Y_{a_{2}a'_{1}; i}^{a'_{3};(3)}\otimes 
\sigma_{123}(\Y_{a'_{3}a_{1}; j}^{a_{2};(4)});
\Y_{a_{4}a'_{3}; k}^{a'_{3};(5)}\otimes \Y_{a_{2}a'_{2}; l}^{a_{4};(6)})
\end{eqnarray}
\end{cor}
\pf
This corollary follows directly from the definitions of 
the matrix elements and the formulas
(\ref{alpha}) and (\ref{beta}). 
\epfv

For $a\in \A$, we define an action of modular transformation 
$S$ on $\F_{1;1}$ as follows: For 
$a_{1}\in \mathcal{A}$,  $k=1, \dots, N_{aa_{1}}^{a_{1}}$
and $w_{a}\in W^{a}$, 
let
\begin{eqnarray*}
S(\Psi(\Y)
(w_{a}; \tau))&=&
\Psi(\Y)
\left(\left(-\frac{1}{\tau}\right)^{L(0)}w_{a}; 
-\frac{1}{\tau}\right)\nn
&=&\tr_{W^{a_{1}}}
\Y\left(\mathcal{U}(e^{-2\pi i\frac{z}{\tau}})
\left(-\frac{1}{\tau}\right)^{L(0)}w_{a}, 
e^{-2\pi i \frac{z}{\tau}}\right)
q_{-\frac{1}{\tau}}^{L(0)-\frac{c}{24}}.
\end{eqnarray*}
Here we have used our convention 
$$\left(-\frac{1}{\tau}\right)^{L(0)}=e^{(\log (-\frac{1}{\tau}))L(0)}.$$
Note that by the modular invariance of genus-one one-point functions
proved in  \cite{H7} (see Theorem \ref{mod-inv}), 
$S(\Psi(\Y)
(u; \tau))$ is indeed in $\F_{1;1}$. Thus we do obtain maps 
$S: \F_{1;1}\to \F_{1;1}$. 

For $a\in \A$, this action of $S$ on $\F_{1;1}$ can be 
extended to an action of 
$S$ on $\F_{1;2}^{a}$. For  intertwining 
operators $\Y_{1}$ and $\Y_{2}$ of types for intertwining 
operators 
of types ${W^{a_{1}}\choose W^{a}W^{a_{1}}}$ and 
${W^{a}\choose W^{a_{2}}W^{a'_{2}}}$, respectively, 
we define 
$$(S(\Psi(\Y_{1}\otimes \Y_{2})))
(w_{a}, w_{a'}; z_{1}, z_{2}; \tau)=0$$
when $a\ne a_{2}$ and 
\begin{eqnarray*}
\lefteqn{((S(\Psi(\Y_{1}\otimes \Y_{2})))
(w_{a_{2}}, w_{a'_{2}}; z_{1}, z_{2}; \tau)}\nn
&&=E\Bigg(\tr_{W^{a_{1}}}
\Y_{1}\Bigg(\mathcal{U}(e^{-2\pi i\frac{z_{2}}{\tau}})
\left(-\frac{1}{\tau}\right)^{L(0)}\cdot\nn
&&\quad\quad\quad\quad\quad\quad\quad\quad\quad\cdot \Y_{2}
(w_{a_{2}}, z_{1}-z_{2})w_{a'_{2}}, e^{-2\pi i \frac{z_{2}}{\tau}}\Bigg)
q_{-\frac{1}{\tau}}^{L(0)-\frac{c}{24}}\Bigg)\nn
&&=E\Bigg(\tr_{W^{a_{1}}}
\Y_{1}\Bigg(\mathcal{U}(e^{-2\pi i\frac{z_{2}}{\tau}})
\cdot\nn
&&\quad \cdot \Y_{2}
\left(\left(-\frac{1}{\tau}\right)^{L(0)}w_{a_{2}}, 
-\frac{1}{\tau}z_{1}-\left(-\frac{1}{\tau}z_{2}\right)\Bigg)
\left(-\frac{1}{\tau}\right)^{L(0)}w_{a'_{2}}, e^{-2\pi i \frac{z_{2}}{\tau}}\right)
\cdot\nn
&&\quad\quad\quad\quad\cdot
q_{-\frac{1}{\tau}}^{L(0)-\frac{c}{24}}\Bigg)\nn
&&=\Psi(\Y_{1}\otimes \Y_{2})\left(\left(-\frac{1}{\tau}\right)^{L(0)}w_{a},
\left(-\frac{1}{\tau}\right)^{L(0)}w_{a'}; -\frac{1}{\tau}z_{1}, -\frac{1}{\tau}z_{2};
-\frac{1}{\tau}\right).
\end{eqnarray*}

We have:

\begin{prop}
We have the following formula:
\begin{equation}\label{formula2-3}
S\alpha=\beta S.
\end{equation}
\end{prop}
\pf
We have 
\begin{eqnarray*}
\lefteqn{(\beta(S(\Psi(\Y_{aa_{1}; p}^{a_{1};(1)}\otimes 
\Y_{a_{2}a'_{2}; q}^{a; (2)}))))(w_{a_{2}}\otimes 
w_{a_{2}'}; z_{1}, z_{2}; \tau)}\nn
&&=(S(\Psi(\Y_{aa_{1}; p}^{a_{1};(1)}\otimes 
\Y_{a_{2}a'_{2}; q}^{a; (2)})))
(w_{a_{2}}\otimes w_{a_{2}'}; z_{1}, z_{2}+\tau; \tau)\nn
&&=E\Bigg(\tr_{W^{a_{1}}}
\Y_{aa_{1}; p}^{a_{1}}(\mathcal{U}(e^{-2\pi i\frac{z_{2}+\tau}{\tau}})
\left(-\frac{1}{\tau}\right)^{L(0)}\cdot\nn
&&\quad\quad\quad\quad\quad\cdot \Y_{a_{2}a'_{2}; q}^{a}
(w_{a_{2}}, z_{1}-(z_{2}+\tau))w_{a'_{2}}, 
e^{-2\pi i\frac{z_{2}+\tau}{\tau}})
q_{-\frac{1}{\tau}}^{L(0)-\frac{c}{24}}\Bigg)\nn
&&=E\Bigg(\tr_{W^{a_{1}}}
\Y_{aa_{1}; p}^{a_{1}}\Bigg(\mathcal{U}(e^{2\pi i (-\frac{z_{2}}{\tau}-1)})\cdot\nn
&&\quad\quad\quad \quad\cdot
\Y_{a_{2}a'_{2}; q}^{a}
\Bigg(\left(-\frac{1}{\tau}\right)^{L(0)}w_{a_{2}}, 
-\frac{1}{\tau}z_{1}-\left(-\frac{1}{\tau}z_{2}-1\right)\Bigg)\cdot\nn
&&\quad\quad\quad \quad\cdot \left(-\frac{1}{\tau}\right)^{L(0)}
w_{a'_{2}}, e^{2\pi i (-\frac{z_{2}}{\tau}-1)}\Bigg)
q_{-\frac{1}{\tau}}^{L(0)-\frac{c}{24}}\Bigg)\nn
&&=(\Psi(\Y_{aa_{1}; p}^{a_{1};(1)}\otimes 
\Y_{a_{2}a'_{2}; q}^{a; (2)}))
\Biggl(\Biggl(\left(-\frac{1}{\tau}\right)^{L(0)}w_{a_{2}}\Biggr)\otimes
\Biggl(\left(-\frac{1}{\tau}\right)^{L(0)}w_{a_{2}'}\Biggr); \nn
&&\quad\quad\quad\quad\quad\quad\quad\quad
\quad\quad\quad\quad\quad\quad\quad\quad
\quad\quad\quad\quad-\frac{1}{\tau}z_{1}, -\frac{1}{\tau}z_{2}-1; 
-\frac{1}{\tau}\Biggr)\nn
&&=(\alpha(\Psi(\Y_{aa_{1}; p}^{a_{1};(1)}\otimes 
\Y_{a_{2}a'_{2}; q}^{a; (2)})))
\Biggl(\Biggl(\left(-\frac{1}{\tau}\right)^{L(0)}w_{a_{2}}\Biggr)\otimes
\Biggl(\left(-\frac{1}{\tau}\right)^{L(0)}w_{a_{2}'}\Biggr); \nn
&&\quad\quad\quad\quad\quad\quad\quad\quad
\quad\quad\quad\quad\quad\quad\quad\quad
\quad\quad\quad\quad\quad\quad
 -\frac{1}{\tau}z_{1}, -\frac{1}{\tau}z_{2}; 
-\frac{1}{\tau}\Biggr)\nn
&&=(S(\alpha(\Psi_{a_{1}, a_{2}, e}^{1, 1})))
(w_{a_{2}}, 
w_{a_{2}'}; z_{1}, z_{2}; \tau).
\end{eqnarray*}
Thus we obtain
(\ref{formula2-3}).
\epfv

For fixed $a\in \A$, we know that
$\Psi(\Y_{aa_{1}; p}^{a_{1}; (1)})$ for $a_{1}\in \A$
form a basis of $\F^{a}_{1;1}$ and
there exist unique 
$S(\Y_{aa_{1}; p}^{a_{1}; (1)}; 
\Y_{aa_{2}; k}^{a_{2}; (2)})$ in $\C$ for $a_{1}, a_{2}\in \A$,
$p=1, \dots, N_{aa_{1}}^{a_{1}}$ and $k=1, \dots, N_{aa_{2}}^{a_{2}}$
such that 
\begin{eqnarray}
S(\Psi(\Y_{aa_{1}; p}^{a_{1}; (1)}))
&=&\sum_{a_{2}\in \A}\sum_{k=1}^{N_{aa_{2}}^{a_{2}}}
S(\Y_{aa_{1}; p}^{a_{1}; (1)}; 
\Y_{aa_{2}; k}^{a_{2}; (2)})
\Psi(\Y_{aa_{2}; k}^{a_{2}; (2)})\label{formula2-1}.
\end{eqnarray}
Clearly we also have 
$$
S(\Psi(\Y_{aa_{1}; p}^{a_{1};(1)}\otimes 
\Y_{a_{2}a'_{2}; q}^{a; (2)}))
=\sum_{a_{3}\in \A}\sum_{k=1}^{N_{aa_{3}}^{a_{3}}}
S(\Y_{aa_{1}; p}^{a_{1}; (1)}; 
\Y_{aa_{3}; k}^{a_{3}; (3)})
\Psi(\Y_{aa_{3}; k}^{a_{3};(3)}\otimes 
\Y_{a_{2}a'_{2}; q}^{a; (2)})
$$
for $a_{1}, a_{2}, a_{3}\in \A$.

The following theorem is a generalization of Theorem 
4.6 in \cite{H10}:

\begin{thm}
For $a_{1}, a_{2}, a_{3}\in \A$, we have
\begin{eqnarray}\label{formula2}
\lefteqn{\sum_{k=1}^{N_{a_{4}a_{1}}^{a_{1}}}
S(\Y_{a_{4}a_{1}; k}^{a_{1}; (3)}; 
\Y_{a_{4}a_{5}; m}^{a_{5}; (5)})\cdot}\nn
&&\cdot(B^{(-1)})^{2}
(\sigma_{12}(\Y_{a_{3}a_{1}; p}^{a_{1};{(1)}})\otimes 
\sigma_{12}(\Y_{a_{2}a'_{2}; q}^{a_{3};(2)});
\sigma_{12}(\Y_{a_{4}a_{1}; k}^{a_{1};(3)})\otimes 
\sigma_{12}(\Y_{a_{2}a'_{2}; l}^{a_{4};(4)}))\nn
&&=\sum_{a_{6}\in \A}\sum_{r=1}^{N_{a_{3}a_{6}}^{a_{6}}}
\sum_{i=1}^{N_{a_{2}a_{5}}^{a_{6}}}
\sum_{j=1}^{N_{a'_{2}a_{6}}^{a_{5}}}
e^{\pi ih_{a_{4}}}S(\Y_{a_{3}a_{1}; p}^{a_{1}; (1)}; 
\Y_{a_{3}a_{6}; r}^{a_{6}; (6)})\cdot\nn
&&\quad \quad \cdot F(\sigma_{12}(\Y_{a_{3}a_{6}; r}^{a_{6};(6)})\otimes 
\sigma_{12}(\Y_{a_{2}a'_{2}; q}^{a_{3};(2)});
\sigma_{123}(\Y_{a_{2}a'_{6}; i}^{a'_{5};(7)}) \otimes 
\sigma_{132}(\Y_{a'_{5}a_{6}; j}^{a_{2};(8)}))\cdot \nn
&&\quad\quad \cdot 
F(\Y_{a_{2}a'_{6}; i}^{a'_{5};(7)}\otimes 
\sigma_{123}(\Y_{a'_{5}a_{6}; j}^{a_{2};(8)});
\Y_{a_{4}a'_{5}; k}^{a'_{5};(5)}\otimes \Y_{a_{2}a'_{2}; l}^{a_{4};(4)}).
\end{eqnarray}
\end{thm}
\pf
This follows from (\ref{formula2-3}), (\ref{alpha-1}) and (\ref{beta-1})
immediately.
\epfv

\begin{cor}
For $a_{1}, a_{2}, a_{3}\in \A$, we have
\begin{eqnarray}\label{formula2-cor}
\lefteqn{S^{a_{1}}_{e}(B^{(-1)})^{2}
(\sigma_{12}(\Y_{a_{3}a_{1}; p}^{a_{1};{(1)}})\otimes 
\sigma_{12}(\Y_{a_{2}a'_{2}; q}^{a_{3};(2)});
\Y_{a_{1}e; 1}^{a_{1}}\otimes 
\Y_{a'_{2}a_{2}; l}^{e}))}\nn
&&=\sum_{r=1}^{N_{a_{3}a_{2}}^{a_{2}}}
S(\Y_{a_{3}a_{1}; p}^{a_{1}; (1)}; 
\Y_{a_{3}a_{2}; r}^{a_{2}; (6)})\cdot\nn
&&\quad \quad \cdot F(\sigma_{12}(\Y_{a_{3}a_{2}; r}^{a_{2};(6)})\otimes 
\sigma_{12}(\Y_{a_{2}a'_{2}; q}^{a_{3};(2)});
\Y_{ea_{2}; 1}^{a_{2}} \otimes 
\Y_{a_{2}a'_{2}; 1}^{e})\cdot \nn
&&\quad\quad \cdot 
F(\Y_{a_{2}a'_{2}; 1}^{e}\otimes 
\Y_{a'_{2}e; 1}^{a'_{2}};
\Y_{ee; 1}^{e}\otimes \Y_{a_{2}a'_{2}; 1}^{e}).
\end{eqnarray}
\end{cor}
\pf 
This follows immediately from   (\ref{formula2}) by taking 
$a_{4}=a_{5}=e$.
\epfv

\begin{lemma}\label{triv-f-coef}
For any $a\in \A$, 
$$F(\Y_{aa'; 1}^{e}\otimes 
\Y_{a'e; 1}^{a'};
\Y_{ee; 1}^{e}\otimes \Y_{aa'; 1}^{e})=1.$$
\end{lemma}
\pf
Note that $\Y_{ee; 1}^{e}=Y$ is the vertex operator 
map for the vertex operator algebra and 
$\Y_{ea; 1}^{a}=Y_{W^{a}}$ is the vertex operator 
map for the $V$-module $W^{a}$. For $u, v\in V$, 
$w_{a}\in W^{a}$ and $w_{a'}\in W^{a'}$,
we have 
\begin{eqnarray}\label{triv-f-coef-1}
\lefteqn{(u, \Y_{ee; 1}^{e}(\Y_{aa'; 1}^{e}(w_{a}, z_{1}-z_{2})
w_{a'}, z_{2})v)}\nn
&&=(u, Y(\Y_{aa'; 1}^{e}(w_{a}, z_{1}-z_{2})
w_{a'}, z_{2})v)\nn
&&=(e^{z_{2}L(1)}u, Y(v, -z_{2})\Y_{aa'; 1}^{e}(w_{a}, z_{1}-z_{2})
w_{a'})
\end{eqnarray}
in the region $|z_{2}|>|z_{1}-z_{2}|>0$.
Since $\Y_{aa'; 1}^{e}$ is an intertwining operator, we know that 
the product of $Y$ and $\Y_{aa'; 1}^{e}$ satisfy commutativity. 
Thus the analytic extension of 
the right-hand side of (\ref{triv-f-coef-1}) to the region 
$|z_{1}-z_{2}|>|z_{2}|>0$ is equal to
\begin{eqnarray}\label{triv-f-coef-2}
\lefteqn{(e^{z_{2}L(1)}u, \Y_{aa'; 1}^{e}(w_{a}, z_{1}-z_{2})
Y_{W^{a'}}(v, -z_{2})w_{a'})}\nn
&&=(u, e^{z_{2}L(-1)}\Y_{aa'; 1}^{e}(w_{a}, z_{1}-z_{2})
Y_{W^{a'}}(v, -z_{2})w_{a'}).
\end{eqnarray}
Using the $L(-1)$-conjugation formula for intertwining 
operators and the definition of the intertwining operator 
$\Y_{a'e;1}^{a'}$, we see that in the region given by
$|z_{1}-z_{2}|>|z_{2}|>0$ and $|z_{1}|>|z_{2}|>0$,
the right-hand side of 
(\ref{triv-f-coef-2}) is equal to
\begin{eqnarray}\label{triv-f-coef-3}
\lefteqn{(u, \Y_{aa'; 1}^{e}(w_{a}, z_{1})
e^{z_{2}L(-1)}Y_{W^{a'}}(v, -z_{2})w_{a'})}\nn
&&=(u, \Y_{aa'; 1}^{e}(w_{a}, z_{1})
\Y_{a'e;1}^{a'}(w_{a'}, z_{2})v).
\end{eqnarray}
Thus we see that the left-hand side of (\ref{triv-f-coef-1})
can be analytically extended to the right-hand side of 
(\ref{triv-f-coef-3}). Since the left-hand side of (\ref{triv-f-coef-1})
and the right-hand side of 
(\ref{triv-f-coef-3}) are well defined in the regions 
$|z_{2}|>|z_{1}-z_{2}|>0$ and $|z_{1}|>|z_{2}|>0$, respectively,
they are equal in the intersection $|z_{1}|>|z_{2}|>|z_{1}-z_{2}|>0$.
By the definition of the fusing matrix element 
$F(\Y_{aa'; 1}^{e}\otimes 
\Y_{a'e; 1}^{a'};
\Y_{ee; 1}^{e}\otimes \Y_{aa'; 1}^{e})$, we see that it is 
equal to $1$. 
\epfv

We also have:

\begin{lemma}\label{bilinear-f}
For any $a_{2}, a_{3}\in \A$, we have
$$F(\sigma_{12}(\Y_{a_{3}a_{2}; r}^{a_{2};(6)})\otimes 
\sigma_{12}(\Y_{a_{2}a'_{2}; q}^{a_{3};(2)});
\Y_{ea_{2}; 1}^{a_{2}} \otimes 
\Y_{a_{2}a'_{2}; 1}^{e})=\langle 
\sigma_{132}(\Y_{a_{3}a_{2}; r}^{a_{2};(6)}), 
\sigma_{12}(\Y_{a_{2}a'_{2}; q}^{a_{3};(2)})
\rangle_{\V_{a_{2}a'_{2}}^{a'_{3}}}.$$
\end{lemma}
\pf
This follows directly from (3.3) in \cite{HK}. 
\epfv

For $a_{1}, a_{2}, a_{3}\in \A$, let
$$(\cdot, 
\cdot)_{\V_{a_{1}a_{2}}^{a_{3}}}=\frac{\sqrt{F_{a_{3}}}}
{\sqrt{F_{a_{1}}}\sqrt{F_{a_{2}}}}\langle \cdot, 
\cdot\rangle_{\V_{a_{1}a_{2}}^{a_{3}}},$$
where $\langle \cdot, 
\cdot\rangle_{\V_{a_{1}a_{2}}^{a_{3}}}$ is the bilinear form 
constructed in \cite{HK}.
Then Proposition 3.7 in \cite{HK} actually 
says that $(\cdot, \cdot)_{\V_{a_{1}a_{2}}^{a_{3}}}$
is invariant under the action of 
$S_{3}$. 

Now for $a_{2}, a_{3}\in \A$, we choose
$$\Y_{a_{2}a'_{2}; i}^{a_{3}; (2)}
=\sigma_{13}(\Y_{a'_{3}a'_{2}; i}^{'; a'_{2}; (6)})$$
for $i=1, \dots, N_{a'_{3}a'_{2}}^{a'_{2}}$, where 
$\{\Y_{a'_{1}a'_{2}; i}^{'; a'_{3}}\;|\; 
i=1, \dots, N_{a'_{1}a'_{2}}^{a'_{3}}\}$ is the dual basis of 
$\{\Y_{a_{1}a_{2}; i}^{a_{3}}\;|\;i=1, \dots, N_{a_{1}a_{2}}^{a_{3}}\}$
with respect to the bilinear form $\langle \cdot, 
\cdot\rangle_{\V_{a_{1}a_{2}}^{a_{3}}}$.

The following result generalizes (5.15) in \cite{H10}, which, as was shown 
in \cite{H11}, is equivalent to the nondegeneracy of the semisimple
ribbon (tensor) category of modules for $V$:

\begin{thm}
For $a_{1}, a_{2}, a_{3}\in \A$, we have
\begin{eqnarray}\label{s-formula}
\lefteqn{S(\Y_{a_{3}a_{1}; p}^{a_{1}; (1)}; 
\Y_{a_{3}a_{2}; r}^{a_{2}; (6)})}\nn
&&=\frac{F_{a_{3}}S_{e}^{e}}{F_{a_{1}}F_{a_{2}}}(B^{(-1)})^{2}
(\sigma_{12}(\Y_{a_{3}a_{1}; p}^{a_{1};{(1)}})\otimes 
\sigma_{132}(\Y_{a'_{3}a'_{2}; r}^{'; a'_{2};(6)});
\Y_{a_{1}e; 1}^{a_{1}}\otimes 
\Y_{a'_{2}a_{2}; 1}^{e}).\nn
&&
\end{eqnarray}
\end{thm}
\pf
By Lemma \ref{bilinear-f}, we have
\begin{eqnarray}\label{f-coef-2}
\lefteqn{F(\sigma_{12}(\Y_{a_{3}a_{2}; r}^{a_{2};(6)})\otimes 
\sigma_{12}(\Y_{a_{2}a'_{2}; q}^{a_{3};(2)});
\Y_{ea_{2}; 1}^{a_{2}} \otimes 
\Y_{a_{2}a'_{2}; 1}^{e})}\nn
&&=\langle 
\sigma_{132}(\Y_{a_{3}a_{2}; r}^{a_{2};(6)}), 
\sigma_{12}(\Y_{a_{2}a'_{2}; q}^{a_{3};(2)})
\rangle_{\V_{a_{2}a'_{2}}^{a'_{3}}}\nn
&&=\frac{\sqrt{F_{a_{2}}}\sqrt{F_{a'_{2}}}}{\sqrt{F_{a'_{3}}}}
(\sigma_{132}(\Y_{a_{3}a_{2}; r}^{a_{2};(6)}), 
\sigma_{12}(\Y_{a_{2}a'_{2}; q}^{a'_{3};(2)}))
_{\V_{a_{2}a'_{2}}^{a'_{3}}}\nn
&&=\frac{\sqrt{F_{a_{2}}}\sqrt{F_{a'_{2}}}}{\sqrt{F_{a'_{3}}}}
(\sigma_{13}(\Y_{a_{3}a_{2}; r}^{a_{2};(6)}), 
\Y_{a'_{2}a_{2}; q}^{a'_{3};(2)})
_{\V_{a'_{2}a_{2}}^{a'_{3}}}\nn
&&=\frac{\sqrt{F_{a_{2}}}\sqrt{F_{a'_{2}}}}{\sqrt{F_{a'_{3}}}}
(\sigma_{13}(\Y_{a_{3}a_{2}; r}^{a_{2};(6)}), 
\sigma_{13}(\Y_{a'_{3}a'_{2}; q}^{'; a'_{2};(6)}))
_{\V_{a'_{2}a_{2}}^{a'_{3}}}\nn
&&=\frac{\sqrt{F_{a_{2}}}\sqrt{F_{a'_{2}}}}{\sqrt{F_{a'_{3}}}}
(\Y_{a_{3}a_{2}; r}^{a_{2};(6)}, 
\Y_{a'_{3}a'_{2}; q}^{'; a'_{2};(6)})
_{\V_{a_{3}a_{2}}^{a_{2}}}\nn
&&=\frac{\sqrt{F_{a_{2}}}\sqrt{F_{a'_{2}}}}{\sqrt{F_{a'_{3}}}}
\frac{\sqrt{F_{a_{2}}}}{\sqrt{F_{a_{3}}}\sqrt{F_{a_{2}}}}
\langle\Y_{a_{3}a_{2}; r}^{a_{2};(6)}, 
\Y_{a'_{3}a'_{2}; q}^{'; a'_{2};(6)}\rangle
_{\V_{a_{3}a_{2}}^{a_{2}}}\nn
&&=\frac{\sqrt{F_{a_{2}}}\sqrt{F_{a_{2}}}}{\sqrt{F_{a_{3}}}}
\frac{\sqrt{F_{a_{2}}}}{\sqrt{F_{a_{3}}}\sqrt{F_{a_{2}}}}
\langle\Y_{a_{3}a_{2}; r}^{a_{2};(6)}, 
\Y_{a'_{3}a'_{2}; q}^{'; a'_{2};(6)}\rangle
_{\V_{a_{3}a_{2}}^{a_{2}}}\nn
&&=\frac{F_{a_{2}}}{F_{a_{3}}}\delta_{rq}.
\end{eqnarray}
Using Lemma \ref{triv-f-coef}, (\ref{f-coef-2}), the formula
$$S_{e}^{a_{1}}=\frac{S_{e}^{e}}{F_{a_{1}}}$$
and (\ref{formula2-cor}), we obtain (\ref{s-formula}).
\epfv

\begin{lemma}
For $a_{1}, a_{2}, a_{3}\in \A$, we have
\begin{equation}\label{s-iden-lemma}
S(\sigma_{23}(\Y_{a_{3}a_{1}; p}^{a_{1};(1)}); 
\sigma_{23}(\Y_{a_{3}a_{2}; r}^{a_{2};(6)}))
=S(\Y_{a_{3}a_{1}; p}^{a_{1};(1)}; 
\Y_{a_{3}a_{2}; r}^{a_{2};(6)}).
\end{equation}
\end{lemma}
\pf
By (\ref{formula2-1}) and the definition of 
$\Psi_{a_{1}, a_{3}}^{k}$, we have
\begin{eqnarray}\label{s-iden-lemma-1}
\lefteqn{\tr_{W^{a_{1}}}\Y_{a_{3}a_{1}; p}^{a_{1};(1)}
\left(\mathcal{U}(e^{-2\pi i\frac{z}{\tau}})
\left(-\frac{1}{\tau}\right)^{L(0)}
w_{a_{3}}, e^{-2\pi i\frac{z}{\tau}}\right)
q_{-\frac{1}{\tau}}^{L(0)-\frac{c}{24}}}\nn
&&=\sum_{a_{2}\in \A}\sum_{r=1}^{N_{a_{3}a_{2}}^{a_{2}}}
S(\Y_{a_{3}a_{1}; p}^{a_{1};(1)}; 
\Y_{a_{3}a_{2}; r}^{a_{2};(6)})
\tr_{W^{a_{2}}}\Y_{a_{3}a_{2}; r}^{a_{2};(6)}
(\mathcal{U}(e^{2\pi iz})w_{a_{3}}, e^{2\pi iz})
q_{\tau}^{L(0)-\frac{c}{24}}.\nn
\end{eqnarray}
On the other hand, we have
\begin{eqnarray}\label{s-iden-lemma-2}
\lefteqn{\tr_{W^{a_{1}}}\Y_{a_{3}a_{1}; p}^{a_{1};(1)}
\left(\mathcal{U}(e^{-2\pi i\frac{z}{\tau}})\left(-\frac{1}{\tau}\right)^{L(0)}
w_{a_{3}}, e^{-2\pi i\frac{z}{\tau}}\right)
q_{-\frac{1}{\tau}}^{L(0)-\frac{c}{24}}}\nn
&&=\tr_{W^{a_{1}}}\sigma_{23}^{2}(\Y_{a_{3}a_{1}; p}^{a_{1};(1)})
\left(\mathcal{U}(e^{-2\pi i\frac{z}{\tau}})\left(-\frac{1}{\tau}\right)^{L(0)}
w_{a_{3}}, e^{-2\pi i\frac{z}{\tau}}\right)
q_{-\frac{1}{\tau}}^{L(0)-\frac{c}{24}}\nn
&&=\tr_{(W^{a_{1}})'}\sigma_{23}(\Y_{a_{3}a_{1}; p}^{a_{1};(1)})
\left(\mathcal{U}(e^{2\pi i\frac{z}{\tau}})e^{-\pi iL(0)}
\left(-\frac{1}{\tau}\right)^{L(0)}w_{a_{3}}, e^{2\pi i\frac{z}{\tau}}\right)
q_{-\frac{1}{\tau}}^{L(0)-\frac{c}{24}}\nn
&&=\sum_{a_{2}\in \A}\sum_{r=1}^{N_{a_{3}a_{2}}^{a_{2}}}
S(\sigma_{23}(\Y_{a_{3}a_{1}; p}^{a_{1};(1)}); 
\sigma_{23}(\Y_{a_{3}a_{2}; r}^{a_{2};(6)}))\cdot\nn
&&\quad\quad\quad\quad\quad\quad\cdot 
\tr_{(W^{a_{2}})'}\sigma_{23}(\Y_{a_{3}a_{2}; r}^{a_{2};(6)})
(\mathcal{U}(e^{-2\pi iz})e^{-\pi iL(0)}w_{a_{3}}, e^{-2\pi iz})
q_{\tau}^{L(0)-\frac{c}{24}}\nn
&&=\sum_{a_{2}\in \A}\sum_{r=1}^{N_{a_{3}a_{2}}^{a_{2}}}
S(\sigma_{23}(\Y_{a_{3}a_{1}; p}^{a_{1};(1)}); 
\sigma_{23}(\Y_{a_{3}a_{2}; r}^{a_{2};(6)}))\cdot\nn
&&\quad\quad\quad\quad\quad\quad\cdot 
\tr_{W^{a_{2}}}\Y_{a_{3}a_{2}; r}^{a_{2};(6)}
(\mathcal{U}(e^{2\pi iz})w_{a_{3}}, e^{2\pi iz})
q_{\tau}^{L(0)-\frac{c}{24}}.
\end{eqnarray}
From (\ref{s-iden-lemma-1}) and (\ref{s-iden-lemma-2}),
we obtain (\ref{s-iden-lemma}).
\epfv

The following result is 
a generalization of Theorem 5.6 in \cite{H10}:

\begin{thm}\label{s-iden-thm}
For $a_{1}, a_{2}, a_{3}\in \A$, we have
\begin{equation}\label{s-iden}
S(\Y_{a'_{3}a'_{2}; r}^{';a'_{2}; (6)}; 
\Y_{a'_{3}a'_{1}; p}^{'; a'_{1}; (1)})
=S(\Y_{a_{3}a_{1}; p}^{a_{1};(1)}; 
\Y_{a_{3}a_{2}; r}^{a_{2};(6)}).
\end{equation}
\end{thm}
\pf
From (\ref{s-formula}), we obtain
\begin{eqnarray}\label{s-iden-1}
\lefteqn{S(\Y_{a'_{3}a'_{2}; r}^{';a'_{2}; (6)}; 
\Y_{a'_{3}a'_{1}; p}^{'; a'_{1}; (1)})}\nn
&&=\frac{F_{a'_{3}}S_{e}^{e}}{F_{a'_{1}}F_{a'_{2}}}(B^{(-1)})^{2}
(\sigma_{12}(\Y_{a'_{3}a'_{2}; r}^{'; a'_{2};(6)})\otimes 
\sigma_{132}(\Y_{a_{3}a_{1}; p}^{a_{1};(1)});
\Y_{a'_{2}e; 1}^{a'_{2}}\otimes 
\Y_{a_{1}a'_{1}; 1}^{e}).\nn
&&
\end{eqnarray}
Using (3.1), (3.5), (3.12) in \cite{H10}, the relations
$\sigma^{2}_{132}=\sigma_{123}=\sigma_{12}\sigma_{23}$,
$\sigma_{123}\sigma_{12}=\sigma_{132}\sigma_{23}$, 
$\sigma_{12}\sigma_{132}\sigma_{12}=\sigma_{123}$,
$\sigma_{12}\sigma_{123}\sigma_{12}=\sigma_{132}$ and 
$h_{a}=h_{a'}$ for $a\in \A$, we obtain 
\begin{eqnarray}\label{s-iden-2}
\lefteqn{(B^{(-1)})^{2}
(\sigma_{12}(\Y_{a'_{3}a'_{2}; r}^{'; a'_{2};(6)})\otimes 
\sigma_{132}(\Y_{a_{3}a_{1}; p}^{a_{1};(1)});
\Y_{a'_{2}e; 1}^{a'_{2}}\otimes 
\Y_{a_{1}a'_{1}; 1}^{e})}\nn
&&=\sum_{a_{4}\in \A}\sum_{k=1}^{N_{a'_{2}a_{1}}^{a_{4}}}
\sum_{l=1}^{N_{a_{4}a'_{1}}^{a'_{2}}}
F(\sigma_{12}(\Y_{a'_{3}a'_{2}; r}^{'; a'_{2};(6)})\otimes 
\sigma_{132}(\Y_{a_{3}a_{1}; p}^{a_{1};(1)});
\Y_{a_{4}a'_{1}; l}^{a'_{2}; (3)}
\otimes \Y_{a'_{2}a_{1}; k}^{a_{4}; (4)})\cdot\nn
&&\quad\quad\quad\quad\cdot 
e^{-2\pi i (h_{a_{4}}-h_{a'_{2}}-h_{a_{1}})}
F^{-1}(\Y_{a_{4}a'_{1}; l}^{a'_{2}; (3)}\otimes \Y_{a'_{2}a_{1}; k}^{a_{4}; (4)}; 
\Y_{a'_{2}e; 1}^{a'_{2}}\otimes 
\Y_{a_{1}a'_{1}; 1}^{e}) \nn
&&=\sum_{a_{4}\in \A}\sum_{k=1}^{N_{a'_{2}a_{1}}^{a_{4}}}
\sum_{l=1}^{N_{a_{4}a'_{1}}^{a'_{2}}}
e^{-2\pi i (h_{a_{4}}-h_{a'_{2}}-h_{a_{1}})}\cdot\nn
&&\quad\cdot F(\sigma^{2}_{132}(\Y_{a_{3}a_{1}; p}^{a_{1};(1)})
\otimes \sigma_{123}(\sigma_{12}(\Y_{a'_{3}a'_{2}; r}^{'; a'_{2};(6)}));
\sigma_{123}(\Y_{a'_{2}a_{1}; k}^{a_{4}; (4)})
\otimes \sigma_{132}(\Y_{a_{4}a'_{1}; l}^{a'_{2}; (3)}))\cdot\nn
&&\quad\quad\quad\quad\cdot 
F(\sigma_{12}(\Y_{a_{4}a'_{1}; l}^{a'_{2}; (3)})\otimes 
\sigma_{12}(\Y_{a'_{2}a_{1}; k}^{a_{4}; (4)}); 
\sigma_{12}(\Y_{a'_{2}e; 1}^{a'_{2}})\otimes 
\sigma_{12}(\Y_{a_{1}a'_{1}; 1}^{e})) \nn
&&=\sum_{a_{4}\in \A}\sum_{k=1}^{N_{a'_{2}a_{1}}^{a_{4}}}
\sum_{l=1}^{N_{a_{4}a'_{1}}^{a'_{2}}}
e^{-2\pi i (h_{a_{4}}-h_{a'_{2}}-h_{a_{1}})}\cdot\nn
&&\quad \cdot F(\sigma^{2}_{132}(\Y_{a_{3}a_{1}; p}^{a_{1};(1)})
\otimes \sigma_{123}(\sigma_{12}(\Y_{a'_{3}a'_{2}; r}^{'; a'_{2};(6)}));
\sigma_{123}(\Y_{a'_{2}a_{1}; k}^{a_{4}; (4)})
\otimes \sigma_{132}(\Y_{a_{4}a'_{1}; l}^{a'_{2}; (3)}))\cdot\nn
&&\quad\quad\quad\quad\cdot 
F(\sigma_{132}(\sigma_{12}(\Y_{a'_{2}a_{1}; k}^{a_{4}; (4)}))\otimes 
\sigma_{123}(\sigma_{12}(\Y_{a_{4}a'_{1}; l}^{a'_{2}; (3)})); \nn
&&\quad\quad\quad\quad\quad\quad\quad\quad\quad\quad\quad\quad
\quad\quad
\sigma_{123}(\sigma_{12}(\Y_{a_{1}a'_{1}; 1}^{e}))\otimes 
\sigma_{132}(\sigma_{12}(\Y_{a'_{2}e; 1}^{a'_{2}}))) \nn
&&=\sum_{a_{4}\in \A}\sum_{k=1}^{N_{a'_{2}a_{1}}^{a_{4}}}
\sum_{l=1}^{N_{a_{4}a'_{1}}^{a'_{2}}}
e^{-2\pi i (h_{a_{4}}-h_{a'_{2}}-h_{a_{1}})}\cdot\nn
&&\quad \cdot F(\sigma^{2}_{132}(\Y_{a_{3}a_{1}; p}^{a_{1};(1)})
\otimes \sigma_{123}(\sigma_{12}(\Y_{a'_{3}a'_{2}; r}^{'; a'_{2};(6)}));
\sigma_{123}(\Y_{a'_{2}a_{1}; k}^{a_{4}; (4)})
\otimes \sigma_{132}(\Y_{a_{4}a'_{1}; l}^{a'_{2}; (3)}))\cdot\nn
&&\quad\quad\quad\quad\cdot 
F^{-1}(\sigma_{12}(\sigma_{132}(\sigma_{12}
(\Y_{a'_{2}a_{1}; k}^{a_{4}; (4)})))\otimes 
\sigma_{12}(\sigma_{123}(\sigma_{12}(\Y_{a_{4}a'_{1}; l}^{a'_{2}; (3)}))); \nn
&&\quad\quad\quad\quad\quad\quad\quad\quad\quad\quad
\sigma_{12}(\sigma_{123}(\sigma_{12}(\Y_{a_{1}a'_{1}; 1}^{e})))\otimes 
\sigma_{12}(\sigma_{132}(\sigma_{12}(\Y_{a'_{2}e; 1}^{a'_{2}})))) \nn
&&=\sum_{a'_{4}\in \A}\sum_{k=1}^{N_{a'_{2}a_{1}}^{a_{4}}}
\sum_{l=1}^{N_{a_{4}a'_{1}}^{a'_{2}}}
F(\sigma_{12}(\sigma_{23}(\Y_{a_{3}a_{1}; p}^{a_{1};(1)}))
\otimes \sigma_{132}(\sigma_{23}(\Y_{a'_{3}a'_{2}; r}^{'; a'_{2};(6)}));
\nn
&&\quad\quad\quad\quad\quad\quad\quad\quad\quad\quad\quad\quad
\quad\quad\quad\quad\quad\quad
\sigma_{123}(\Y_{a'_{2}a_{1}; k}^{a_{4}; (4)})
\otimes \sigma_{132}(\Y_{a_{4}a'_{1}; l}^{a'_{2}; (3)}))\cdot\nn
&&\quad\quad\quad\cdot 
e^{-2\pi i (h_{a'_{4}}-h_{a'_{1}}-h_{a_{2}})}\cdot\nn
&&\quad\quad\quad\quad\cdot 
F^{-1}(\sigma_{123}(\Y_{a'_{2}a_{1}; k}^{a_{4}; (4)})\otimes 
\sigma_{132}(\Y_{a_{4}a'_{1}; l}^{a'_{2}; (3)});
\Y_{a'_{1}e; 1}^{a'_{1}}\otimes 
\Y_{a_{2}a'_{2}; 1}^{e}) \nn
&&=(B^{(-1)})^{2}
(\sigma_{12}(\sigma_{23}(\Y_{a_{3}a_{1}; p}^{a_{1};(1)}))\otimes 
\sigma_{132}(\sigma_{23}(\Y_{a'_{3}a'_{2}; r}^{'; a'_{2};(6)}));
\Y_{a'_{1}e; 1}^{a'_{1}}\otimes 
\Y_{a_{2}a'_{2}; 1}^{e}).\nn
&&
\end{eqnarray}
Substituting the right-hand side of (\ref{s-iden-2})
into the right-hand side of (\ref{s-iden-1}) and then using 
(\ref{s-formula}) and (\ref{s-iden-lemma}), we obtain
\begin{eqnarray}\label{s-iden-3}
\lefteqn{S(\Y_{a'_{3}a'_{2}; r}^{';a'_{2}; (6)}; 
\Y_{a'_{3}a'_{1}; p}^{'; a'_{1}; (1)})}\nn
&&=\frac{F_{a'_{3}}S_{e}^{e}}{F_{a'_{1}}F_{a'_{2}}}
(B^{(-1)})^{2}
(\sigma_{12}(\sigma_{23}(\Y_{a_{3}a_{1}; p}^{a_{1};(1)}))\otimes 
\sigma_{132}(\sigma_{23}(\Y_{a'_{3}a'_{2}; r}^{'; a'_{2};(6)}));
\Y_{a'_{1}e; 1}^{a'_{1}}\otimes 
\Y_{a_{2}a'_{2}; 1}^{e})\nn
&&=S(\sigma_{23}(\Y_{a_{3}a_{1}; p}^{a_{1};(1)}); 
\sigma_{23}(\Y_{a_{3}a_{2}; r}^{a_{2};(6)}))\nn
&&=S(\Y_{a_{3}a_{1}; p}^{a_{1};(1)}; 
\Y_{a_{3}a_{2}; r}^{a_{2};(6)}).
\end{eqnarray}
\epfv

\renewcommand{\theequation}{\thesection.\arabic{equation}}
\renewcommand{\thethm}{\thesection.\arabic{thm}}
\setcounter{equation}{0}
\setcounter{thm}{0}

\section{Modular invariance}

We now prove the modular invariance of the conformal
full field algebras over $V\otimes V$ 
constructed in \cite{HK} (see Section 2) for $V$ satisfying 
the conditions in the preceding two sections.

\begin{thm}
The conformal full field algebra over $V\otimes V$ 
given in Theorem \ref{full-const}
is modular invariant.
\end{thm}
\pf
By Theorem \ref{s-mod-inv-thm}, we need only verify 
(\ref{s-mod-inv}). In this case, $\A^{L}=\A^{R}=\A$ and 
we can identify 
the set $\{1, \dots, N\}$ with $\A$. The map $r^{L}$ as 
a map from $\A$ to $\A$ is the identity map and the map 
$^{R}$ as a map from $\A$ to $\A$ is the map $'$. 
Using the basis of the intertwining operators that
we choose in the construction of the conformal full field 
algebra over $V\otimes V$ , for $a_{1}, a_{2}\in \A$, we have
$d_{a_{1}, a_{2}; i, j}^{a_{2}; (1, 1)}
=\delta_{ij}$. So in this case 
the left-hand side of (\ref{s-mod-inv}) is equal to 
\begin{eqnarray}
\lefteqn{\sum_{a_{1}\in \A}\sum_{i=1}^{N_{aa_{1}}^{a_{1}}}
\sum_{j=1}^{N_{a'a'_{1}}^{a'_{1}}}
\delta_{ij}S(\Y_{aa_{1}; i}^{a_{1}}; 
\Y_{aa_{2}; k}^{a_{2}})
S^{-1}(\Y_{a'a'_{1}; j}^{a'_{1}}; 
\Y_{a'a'_{3}; l}^{a'_{3}})}\nn
&&=\sum_{a_{1}\in \A}\sum_{i=1}^{N_{aa_{1}}^{a_{1}}}
S(\Y_{a'a'_{2}; k}^{a'_{2}}; 
\Y_{a'a'_{1}; i}^{a'_{1}})
S^{-1}(\Y_{a'a'_{1}; i}^{a'_{1}}; 
\Y_{a'a'_{3}; l}^{a'_{3}})\nn
&&=\delta_{a'_{2}a'_{3}}\delta_{kl}\nn
&&=\delta_{a_{2}a_{3}}\delta_{kl}\nn
&&=\sum_{a_{4}\in (r^{L})^{-1}(a_{2})\cap (r^{R})^{-1}(a'_{3})}
\delta_{kl},
\end{eqnarray}
which is indeed equal to the right-hand side of 
(\ref{s-mod-inv}).
By Theorem \ref{s-mod-inv-thm}, the conformal full field algebra
over $V\otimes V$ 
is modular invariant. 
\epfv

\noindent {\small \sc Department of Mathematics, Rutgers University,
110 Frelinghuysen Rd., Piscataway, NJ 08854-8019}

\noindent {\em E-mail address}: yzhuang@math.rutgers.edu

\vspace{1em}

\noindent {\small \sc Max Planck Institute for Mathematics
in the Sciences, Inselstrasse 22, D-04103, Leipzig, Germany}

\vspace{.2em}

\noindent and

\vspace{.2em}

\noindent {\small \sc
Institut Des Hautes \'{E}tudes Scientifiques,
Le Bois-Marie, 35, Route De Chartres,
F-91440 Bures-sur-Yvette, France} (current address)

\noindent {\em E-mail address}: kong@ihes.fr 

\end{document}